\documentclass[a4paper,11pt]{article}
\usepackage{epsfig}
\usepackage{amsfonts}
\usepackage{amsmath}
\usepackage{tikz}
\usetikzlibrary{positioning}
\usepackage[numbers]{natbib}
\usepackage{multirow}				
\usepackage{booktabs}				
\usepackage{pdflscape}
\usepackage{epstopdf}
\usepackage{subfigure}
\usepackage{varwidth}       
\usepackage{multirow}
\usepackage{rotating}
\usepackage{url}
\usepackage{color}
\usepackage{algorithm}
\usepackage{algpseudocode}	

\newtheorem{corollary}{Corollary}

\newcommand{\minitab}[2][l]{\begin{tabular}{#1}#2\end{tabular}}     

\allowdisplaybreaks[1]

\begin{document}

\title{Service Network Design Problem with Capacity-Demand Balancing}
\author{
Yusuf Se\c{c}erdin$^{a}$,
Murat Erkoc$^{a}$\\
\textit{\small yusufsecerdin@miami.edu, merkoc@miami.edu}}
\date{
\begin{footnotesize}
     \begin{flushleft}
     $^a$ Department of Industrial Engineering, University of Miami, Coral Gables, FL 33146\\
     \end{flushleft}
	\end{footnotesize}
         \today
}

\maketitle

\begin{abstract}
This paper revisits a well-known network design problem with a practical focus and addresses developing cost-effective strategies to respond to excessive demand in the service network design problem (SNDP) for freight carriers in a multi-period setting. The common assumption in SNDP literature states that the current capacity of freight carriers' assets is capable of handling all of the forecasted demand; however, we assume that there are certain periods such as holiday season in which excessive demand is observed.  The demand strictly exceeds the carrier's capacity; even though, the average demand can be still fulfilled throughout the year. In this sense, we let the carriers have three options to respond the demand: Dispersing or shifting the demand through an early and a late delivery with a penalty, leasing additional asset(s) temporarily, and outsourcing some capacity. We propose a modeling and solution approach that jointly incorporates asset management and sizing, outsourcing (3PLs), and earliness/tardiness penalties. Given a set of commodities to be delivered from origin terminals to destination terminals in a network over multiple periods, the objective of the problem is to minimize the overall operational costs by optimally selecting and scheduling the home fleet with respect to 'demand shifting' choices, selecting services from third parties, and routing the commodities on the designed service network. We propose an arc-based formulation as well as valid inequalities for this problem and present a comprehensive computational study including additional analysis from operational and computational perspectives on the randomly generated instances. 
The formulations with valid inequalities (VIs) outperform the regular formulation in obtaining tighter lower bounds. One set of VIs can improve the CPU time elapsed by 25\% on medium-instances that can be solved optimally within the time limit. Furthermore, we develop a custom multi-phase dedicate-merge-and-mix algorithm (DMaM) including a construction phase and three improving phases to solve CSSND problem with an emphasis of obtaining solutions as high-quality as possible practically short time to respond the change in demand in the real world. DMaM has a promising potential to obtain solutions for especially very large instances whereas the commercial solver cannot initialize the branch-and-bound algorithm due to excessive memory usage. 
\end{abstract}

\section{Introduction}
\label{sec-introduction}
As the world becoming a global village, producing and manufacturing have been affected positively since firms are able to centralize production operations at thoroughly designed facilities that are equipped with advanced technology (i.e., automated production systems, robots, etc.) through high-capital investments. These facilities are planned to be operated effectively and efficiently aiming to produce or manufacture a high volume of outputs by achieving consistent quality so that their outputs are supplied to almost all of the 
customers from all around the world. Freight transportation plays an important role in the distribution of globally marketed products and, the freight carriers provide delivery service between all kind of facilities including factories, warehouses, depots
for all stakeholders in the global markets such as manufacturers, wholesalers, and retailers. Briefly, today's world economy mostly relies on the movement of freight efficiently. 

Moreover, the carriers encounter new challenges eventually since the demand for fast, reliable and low-cost transportation service increases gradually. Furthermore, the carriers usually operate with low operating margins (or operating income) and aiming to maximize utilization of assets owned, due to the high ownership and operating costs in the consolidation-based transportation sector. For instance, 
FedEx Freight and XPO Logistics,-the first two of the Journal of Commerce's (JOC) list of Top 50 US and Canadian LTL Trucking Companies-, publicize their operating margins as 6.9\% and 4.6\% in 2016, respectively \cite{JOC:2017,FedExAnnualRep:2016,XPOAnnualRep:2016}. Thus, the carriers have to design and operate their transportation networks efficiently and effectively while satisfying customer expectation perfectly to improve and sustain high operating margin.

Designing transportation networks is divided into three levels including strategic (i.e., facility location), tactical (i.e., distribution planning, allocation of demand points to facilities, resource planning), and operational planning (i.e., vehicle routing and scheduling). In this paper, we specifically focus on resource (asset) planning as well as routing and scheduling of assets in the planning horizon for a consolidation-based freight carrier. This process is known as designing of service network and accepted somewhere between the tactical and operational level of the planning process in transportation systems. From the operations research standpoint, service network design problems are the type of fixed charge capacitated multi-commodity network design problems with static or dynamic MIP formulations \cite{Crainic:2000}. Thereafter, asset positioning and balancing are also considered within service network design problems, since the carriers aim to decrease their operating costs and to increase the utilization of owned assets for operational purposes. This effort is known as "full-asset-utilization" policy in transportation systems, which needs asset circulation continuously throughout the service network. Thus, "design-balanced constraints" are introduced into service network design formulations to satisfy circulation of assets by imposing the number of assets entering and leaving a node must be balanced (e.g., \citealp{Barnhart&Schneur:1996, Crainic&Kim:2007, Andersen&Christiansen:2009}, \citealp*{Andersen&Crainic&Christiansen:2009_TRC,Pedersen&Crainic&Madsen:2009}).

The common assumption in the literature usually states that the current capacity of freight carrier's
assets is capable of handling all of the forecasted demand within the upcoming planning horizon (or it is assumed that an unlimited number of assets are available). As an exception to this, Crainic et al.~\cite{Crainic&Hewitt&Toulouse&Vu:2016_TS} assumed that the number of resources is determined in advance. Besides Barnhart and Schneur~\cite{Barnhart&Schneur:1996}, Andersen and Christiansen~\cite{Andersen&Christiansen:2009} and Lai and Lo~\cite{Lai&Lo:2004} initially solved the problem with fixed resources but also proposed fleet composition as an extension. Since their formulation was able to compose the fleet mix as well when the right-hand side values of the regarding constraints were relaxed. Last, Crainic et al.~\cite{Crainic&Hewitt&Toulouse&Vu:2016_EUROJ} considered resource acquisition as a strategic decision due to the assumption is that resources are permanently acquired. None of the studies explicitly examined short-term resource acquisition in case of an urgent response such as a peak-demand period. 

In contrast to the literature on SND problems, we consider a more realistic situation in which demand strictly exceeds the capacity of freight carrier's assets in a certain period motivated by observations in practice. 
This period might be any peak season observed in pre-Thanksgiving (Black Friday), pre-Christmas and new year's eve, and pre-school periods in the US. Regarding these facts and considered capacity assumption, we assume that the carrier takes three actions to be able to respond to peak demand. First, the carrier may choose to deliver the commodity late (demand shifting), which is observed frequently in practice. In the retailing sector, firms including Walmart, Amazon, Bestbuy, etc. offer to deliver items later than usual for free or may ask additional money to deliver them on time. Besides, we also would like to consider the early case and assume that the carrier would disperse (spread) the excess demand by picking up some of the commodities earlier than their release date (assuming it is possible). The carrier may also deliver some of the commodities later than their due date in exchange for earliness/tardiness fee such that existing assets would be able to handle all demand in the planning horizon.  

Second, the carrier would reserve or immediately buy capacity from another service provider, which refers to outsourcing. The concept of outsourcing is usually observed in express shipment delivery, and air carriers prefer to pay for additional capacity from commercial passenger flights instead of operating an asset for certain routes within their transportation networks. In this paper, outsourcing is an option for the freight carriers, and there is a sufficient number of outsourced services available for any route on the network. 

Third, the carrier would expand the capacity by acquiring (i.e., leasing) additional unit(s) of assets temporarily to compensate for the difference between demand and its capacity. Even though the cost of acquiring an additional unit of an asset is not low, the carriers prefer this option frequently to manage and operate their business with complete flexibility. This option corresponds to fleet sizing in the transportation sector.   

Capacity scaling based on excess demand on a transportation network has not been addressed extensively. In this paper, we propose three possible actions such as demand shifting (early/late delivery), outsourcing and fleet sizing (or capacity expansion) that can be taken into account when demand strictly exceeds capacity for a short-term in the planning horizon for a consolidation-based freight carrier. All of these actions are temporary solutions for a peak demand period to get through an issue of capacity shortage. To emphasize the scope, the problem we study still belongs to the tactical level of planning for a transportation system. Among the considered actions above short-term fleet sizing and outsourcing has not been studied very well, demand shifting has never been considered. In contrast to our study, permanent resource acquisition has been addressed in Crainic et al.~\cite{Crainic&Hewitt&Toulouse&Vu:2016_EUROJ}, however, the authors considered it as a long-term strategic decision and positioned their study such that strategic (resource acquisition) and tactical (vehicle routing and scheduling) level decisions are made jointly. Carriers usually generate a schedule weekly and repeat the same schedule for few consecutive weeks in the same month or quarter. Given that the peak period demand lasts as long as a month or five-six weeks, a carrier might need to acquire an asset temporarily for a short time period. From a decision-making perspective, it can be claimed that this acquisition is still within a tactical level of decisions, although \cite{Crainic&Hewitt&Toulouse&Vu:2016_EUROJ} assumed resource acquisition as a long-term strategic decision in service network design. This paper aims to fill this gap in the literature as well as taking one of the alternative actions (i.e., demand shifting and outsourcing) against peak demand.  


The contribution of this paper can be summarized in threefold. First, it recognizes explicit capacity shortage in case of observing peak demand periods ahead and introduces capacity scaling (or alternatively capacity management) in transportation networks as well as an arc-based formulation for the capacity scaling service network design (CSSND) problem. Second, two sets of valid inequalities (VIs) are generated for the proposed formulation to make the formulation more effective.
The formulations with valid inequalities (VIs) outperform the regular formulation in obtaining tighter lower bounds. One set of VIs can improve the CPU time elapsed by 25\% on medium-instances that can be solved optimally within the time limit. Third, we propose a multi-phase dedicate-merge-and-mix (DMaM) algorithm including a construction phase and three improving phases to solve CSSND problem.

The rest of the paper is organized as follows. Section~\S~\ref{sec-literatureRev} reviews the literature and surveys of previous studies that focus on the service network design problem. In Section~\S~\ref{sec-problem_defintion_formulation}, we introduce the capacity scaling problem in transportation networks within the service network design context and present the developed formulation. Section~\S~\ref{sec-computational_exp} presents all computational experiments including solution of regular formulation, valid inequalities and a set of additional analysis. Proposed Dedicate-merge-and-mix algorithm (DMaM) that is specifically developed for the studied problem is explained and experimented on the test instances in Section~\S~\ref{sec-DMaM}. Finally, we conclude this paper in Section~\S~\ref{sec-conclusion}.

\section{Literature Review}
\label{sec-literatureRev}
Service network design (SND) problems are a type of capacitated (fixed charge) multi-commodity network design problems and belong to the class of NP-Hard problems in terms of complexity. The reader can refer to a series of comprehensive reviews about further details pertaining to SND problems (e.g., \citealp{Crainic:2000,Crainic&Kim:2007,Steadieseifi&Dellaert&Nuijten&VanWoensel&Raoufi:2014}). Besides, the reader should also refer to Magnanti and Wong~\cite{Magnanti&Wong:1984} and Balakrishnan et al.~\cite{Balakrishnan&Magnanti&Mirchandani:1997} for extensive reviews about network design problems. 

In reviewing the literature, we specifically focus on studies in which asset positioning and balancing are considered simultaneously while designing service networks. As stated in the introduction, asset balancing  and "design-balanced constraints" refer to the same problem feature, which are addressed in the literature broadly (e.g., \citealp{Andersen&Christiansen&Crainic&Gronhaug:2011}, \citealp*{Andersen&Crainic&Christiansen:2009_TRC}, \citealp{Crainic&Hewitt&Toulouse&Vu:2016_TS, Lai&Lo:2004, Li&Wei&Aneja&Tian:2017}, \citealp{Li&Wei&Aneja&Tian&Cui:2017}, \citealp*{Pedersen&Crainic&Madsen:2009,Teypaz&Schrenk&Cung:2010,Vu&Crainic&Toulouse:2013}).

The issue of asset management; alternatively, design-balanced constraints are one of the main topics for researchers who seek to increase the utilization of assets in transportation systems. Barnhart and Schneur~\cite{Barnhart&Schneur:1996}, Kim et al.~\cite{Kim&Barnhart&Ware&Reinhardt:1999}, Armacost et al.~\cite{Armacost&Barnhart&Ware:2002}, and Barnhart et al.~\cite{Barnhart&Krishnan&Kim&Ware:2002} are the earlier studies which consider asset management in service network design problems specifically in the applications of express shipment delivery. The scope of asset management in these studies includes balancing a number of fleets (or aircraft) and equipment at each node on the network. Thereafter, management of assets are generalized through design-balanced constraints, and the notation is introduced by Pedersen et al.~\cite{Pedersen&Crainic&Madsen:2009} which studied a capacitated multi-commodity network design (CMND) problem and named it as Designed-balanced capacitated multicommodity network design (DBCMND) problem. 

Among the studies which consider asset management, Pedersen et al.~\cite{Pedersen&Crainic&Madsen:2009} developed generic arc- and cycle-based formulations for the DBCMND problem and proposed a tabu search metaheuristics to address the arc-based model for achieving computational efficiency. Andersen et al.~\cite{Andersen&Crainic&Christiansen:2009_TRC} proposed both arc-based and path-based formulations accompanying with two more formulations by mixing arc/cycle design variables and arc/path flow variables for this problem, compare these formulations computationally. Chouman and Crainic~\cite{Chouman&Crainic:2015} considered the same problem and proposed a cutting plane matheuristic with learning mechanism to combine an exact lower-bound computing method and a variable fixing procedure, which feeds a MIP solver. Bai et al.~\cite{Bai&Kendall&Qu&Atkin:2012} studied service network design problem with asset balancing and proposed a guided local search metaheuristic within a multi-start framework. Bai et al.~\cite{Bai&Woodward&Subramanian&Cartlidge:2018} extended this approach by studying a new neighborhood structure for getting a more effective solution strategy. Hewitt~\cite{Hewitt:2010} dealt with human resources (drivers) in service network design considering daily driving limitations imposed by the Department of Transportation (DoT) or union regulations rather than scheduling equipment. 

Moreover, Andersen et al.~\cite{Andersen&Crainic&Christiansen:2009_EJOR} firstly considered the multiple fleets case and fleet coordination in service network design to improve the integration of vehicle management and service network design concepts and introduced \textit{service network design with asset management and multiple fleet coordination problem} (SNDAM-mFC). The authors addressed the synchronization between the collaborating services such as new (internal) and existing (external) services particularly observed in intermodal transportation where ferry services are assumed external services and fixed in terms of arrival/departure times and terminals. Likewise, Lai et al.~\cite{Li&Wei&Aneja&Tian:2017} considered heterogeneous assets in multi-commodity network design and proposed a Tabu-search metaheuristic based on the decomposition of the problem.    

Real-life size instances of the SND problem with asset management are also specifically addressed in the literature (e.g., \cite{Andersen&Christiansen&Crainic&Gronhaug:2011} ,\citealp{Teypaz&Schrenk&Cung:2010,Vu&Crainic&Toulouse:2013}). Teypaz et al.~\cite{Teypaz&Schrenk&Cung:2010} proposed a three-step decomposition-based heuristic algorithm with the objective of profit maximization for a carrier. The step of the proposed heuristic includes network construction, selection of commodities and, vehicle routing and scheduling. Andersen et al.~\cite{Andersen&Christiansen&Crainic&Gronhaug:2011} studied SNDAM problem and proposed firstly a branch-and-price (B\&P) algorithm which is accompanied with a mechanism to add linear relaxation cuts dynamically and an acceleration technique for updating upper-bound. Vu et al.~\cite{Vu&Crainic&Toulouse:2013} presented a three-phase matheuristic in which an exact solver is combined with two heuristic methods. In the first two phases, the heuristic methods, namely tabu search and path relinking, generate as many feasible solutions as possible and reduce the problem size in the first two phases. Then, the exact solver searches on the restricted solution space and is able to solve large size of problem instances.  

The researchers also specifically focused on variants of service network design problems which differentiate by transportation modes as well as applications such as express shipment delivery network design for air carriers (e.g., \citealp{Kim&Barnhart&Ware&Reinhardt:1999,Barnhart&Schneur:1996}, \citealp*{Armacost&Barnhart&Ware:2002}); rail services network design (e.g., \citealp*{Bauer&Bektas&Crainic:2010, Zhu&Crainic&Gendreau:2014}, \citealp{Andersen&Christiansen:2009}). Lai and Lo~\cite{Lai&Lo:2004} considered ferry service network design and ferry fleet management and applied their formulation and two-phase heuristic algorithm on the case of Hong Kong. Crainic and Sgalambro~\cite{Crainic&Sgalambro:2014} studied urban-vehicle service network design problem as well as its variants for two-tier city logistics and proposed formulations for all considered cases. Zhu et al.~\cite{Zhu&Crainic&Gendreau:2014} addressed the scheduled service network design problem for freight rail transportation and proposed a comprehensive formulation which integrates several core decisions belong to tactical planning process such as service selection and scheduling, car classification and blocking, train to make up and to flow of shipments. 

All of the studies above presumes an unlimited number of resources available for service network planning; however, this leads to failure for the objective of maximizing resource utilization ("full-asset-utilization") as emphasized earlier in \S\ref{sec-introduction}. Andersen and Christiansen~\cite{Andersen&Christiansen:2009} decide how many assets (locomotives) to utilize in the operations of Polcorridor rail service network design. Lai and Lo~\cite{Lai&Lo:2004} determine the number of ferries in operation that is limited by the maximum allowable fleet size for Hong Kong ferry service network design. Crainic et al.~\cite{Crainic&Hewitt&Toulouse&Vu:2016_TS} is one of the papers which considers an SND with resource constraints based on the limited total number of resources that is determined in advance. In contrast to this, a total number of resources (or assets) utilized in the planning horizon is determined simultaneously while designing service network in our paper. Barnhart and Schneur~\cite{Barnhart&Schneur:1996} aimed to solve express shipment delivery network design problem for a fixed fleet, besides the proposed model is able to determine optimum fleet composition and size through releasing right-hand side values of available resource constraints. Even though fleet sizing is not addressed directly in the original problem scope, the proposed model is capable of doing that. The authors also determine the number of shipments delivered by commercial air in the same manner with the outsourcing option in that paper.  

Finally, the most similar study to this paper in the literature is Crainic et al.~\cite{Crainic&Hewitt&Toulouse&Vu:2016_EUROJ} which considers strategic resource acquisition and allocation of resources in service network design. In this study, resources are acquired for permanently, and
the acquisition cost is amortized over a series of periods. However, we assume that resources may be acquired temporarily such a leasing contract; thus the only cost of leasing incurs for the period in which the resource is acquired. The authors also take the outsourcing option into account such a way the resource is temporarily acquired from the third-party only for executing a particular service. On the other hand, the outsourcing option in our case does not include third-party's resource acquisition and is not necessary to occupy the resource fully though, the outsourced volume of demand might be a partial load in third-party carrier's resource. Note that the carrier in our paper cannot interfere with the schedule of third-party carrier's resource, it assumed as given for the carrier. The authors named this problem as \textit{scheduled service network design with resource acquisition and management} and proposed a cycle-based formulation as well as a matheuristic which combines column generation, slope scaling, intensification and diversification procedures, and exact optimization. 

To highlight the main differences between our study and the aforementioned similar studies in the SND literature, we classified those studies based on problem features we consider in Table~\ref{tab:lit_comparison}.   

\begin{table}[htbp]
	\centering
	\scriptsize
	\caption {Problem notation.} \label{tab:lit_comparison}
	{
		\begin{tabular}{l c c c c}
			\midrule
			
\multirow{2}{*}{\bf Paper}	&	\bf Excess demand 	&	\bf Fleet sizing 	&	\multirow{2}{*}{\bf Outsourcing}	&	\bf Demand shifting	\\	
	&	\bf (capacity shortage)	&	\bf (or acquisition)	&		&	\bf (earliness/tardiness)	\\	\midrule
{Barnhart and Schneur~\cite{Barnhart&Schneur:1996}}	&	x	&	Not in problem scope	&	\checkmark	&	x	\\	
{Andersen and Christiansen~\cite{Andersen&Christiansen:2009} }	&	x	&	\checkmark	&	x	&	x	\\	
{Lai and Lo~\cite{Lai&Lo:2004}}	&	x	&	\checkmark	&	x	&	x	\\	
{Crainic et al.~\cite{Crainic&Hewitt&Toulouse&Vu:2016_TS}}	&	x	& \color{red}	Finite resources	&	x	&	x	\\	
{Crainic et al.~\cite{Crainic&Hewitt&Toulouse&Vu:2016_EUROJ}}	&	Implicitly	&	Long-term	&	\checkmark	&	x	\\	
This paper	&	\checkmark	&	\checkmark	&	\checkmark	&	\checkmark	\\		

			\midrule	
		\end{tabular}
	}
\end{table}

To sum up, the contribution of our paper can be summarized such that we introduce capacity-demand balancing problem in service network design and specifically consider a planning horizon in which demand exceeds current capacity and the decision maker has to take action(s) among short-term resource acquisition (leasing), outsourcing or demand shifting  (earliness/tardiness) to respond to the demand. To the best of our knowledge, fleet sizing and outsourcing have not been studied very well, besides demand shifting has never been considered before.

\section{Problem Definition and Formulation}
\label{sec-problem_defintion_formulation}


In service network design problems, a service can be defined as transportation of a commodity between its origin and destination with a determined capacity level and constant speed imposed by the asset (or resource in general) as well as known departure and arrival times. A set of selected services comprises a schedule for a fixed length of time (i.e., day, week, or month) and it is usually assumed that the determined schedule is cyclic meaning that it is being repeated for a certain period (i.e., season or year). To operate a service a resource is required, so an \textit{asset} is assigned to a particular service during the schedule. To sum up, it can be stated that the solution of a service network design problem is briefly a schedule in which services are selected, and assets are assigned to those selected services, by doing that routes of the assets with their corresponding departure and arrival times and flow of commodities through the network within allowed time interval are determined.   

\begin{figure}[htbp]
	\centering
	\fbox{\includegraphics[width=0.5\textwidth]{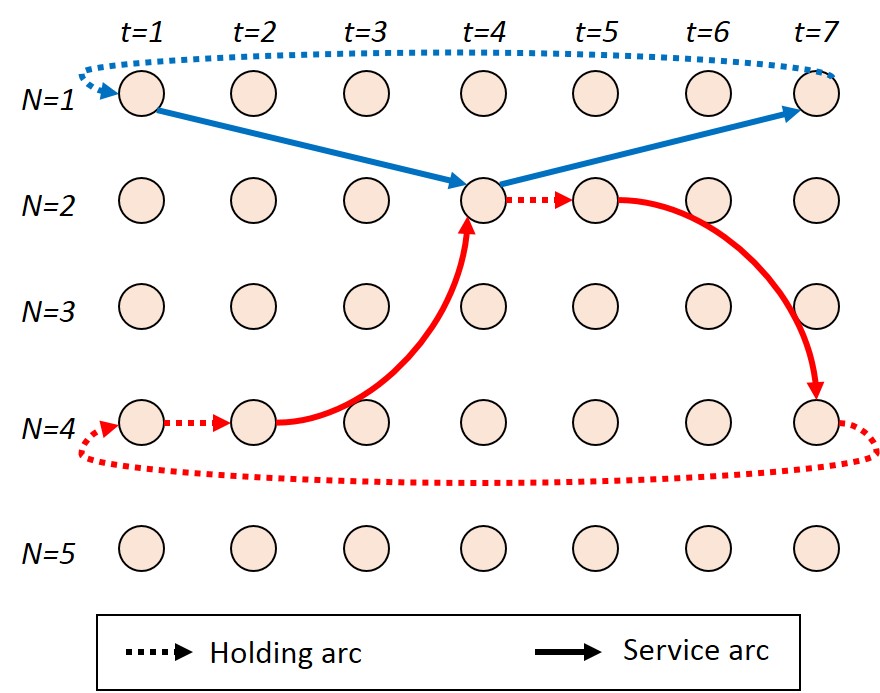}}
	\caption{A sample time-space network with two assets.}
	\label{fig:TSN}
\end{figure}

Network design problems with scheduling decisions are studied on a special network that is called \textit{time-space network} instead of a physical (or static) network to incorporate the time dimension into the problem compared to classical network design problems. A finite planning horizon is divided into the identical time periods, and for each period the nodes of the physical network are duplicated so that any movement on the physical network can be represented both spatially and temporally. Let $N^{'}$ and $N$ represent the set of nodes on the physical and time-space network, respectively, let $T$ denote the set of time periods. An illustration of a time-space network consists of five nodes and seven periods is given in Figure~\ref{fig:TSN}. Note that notation is selected in a way to be consistent with the formulations given in Andersen et al.~\cite{Andersen&Crainic&Christiansen:2009_TRC} and Pedersen et al.~\cite{Pedersen&Crainic&Madsen:2009}.

The arcs on a time-space network represent activity and are classified into two groups, holding arcs and service arcs. A holding arc is directed from a period to another for the same location and represents only time-wise movement (see Figure~\ref{fig:TSN}). The time elapsed during loading/unloading of vehicles, the trailers or cars wait for a transfer to another truck or train might be represented by a holding arc. Thus, any resource is assumed busy on a holding arc within the commodity time-window (when loaded with the commodity) while it is idle when no commodity seizes it (out of commodity time-window). A service arc corresponds to transportation between two locations and the difference between periods of these locations shows the time elapsed during transportation activity. Repositioning of assets also occur via service arcs; however, no flow of commodities takes place in repositioning. In our case, the service arcs are also divided into two categories as arcs of offered services and arcs of outsourced services. Thus, let $A_{h}$ denote holding arcs, let $A_{s}$ and $A_{o}$ symbolize offered services and outsourced services arcs on the time-space network, respectively. We also denote the set of all arcs by $A$ on the time-space network, while the set $A^{'}$ includes all arcs on the physical network. Since the considered problem is capacitated, each arc except holding arcs has a capacity which is denoted by $u_{ij}$. 
  
The distances between nodes on the physical network are defined in terms of periods and assumed to be an integer number denoted by $d_{ij}$. We also assume that distances between nodes on the physical network satisfy the triangle inequality. In addition to this, as it is discussed in Section~\ref{sec-introduction} "full-asset-utilization" policy results with circular routes which means that any asset must return to its starting node at the end of each cycle. To remind that this circulation is enforced through "design-balanced" constraints. This issue has two results on the assumptions, the length of a route (or cycle) for any asset must be equal to the schedule length as \cite{Andersen&Crainic&Christiansen:2009_TRC} stated. We have the following observation for the second result: 

\begin{corollary}
The maximum value for the distance between any pair of nodes on the physical network must be less than or equal to half of the schedule length,; thus any asset is able to make a return trip to its starting node at the beginning of next cycle, or mathematically:	
\begin{eqnarray}
d_{ij} \leq \left\lfloor \frac{|T|}{2}  \right\rfloor \qquad \forall i,j \in N^{'} \label{MaximumDistance} 
\end{eqnarray} 
\end{corollary}    

The reason behind of this observation is that an asset is supposed to be dispatched to its beginning node at the end of each planning horizon; thus any asset should complete any transportation activity in half of the time to be able to make a return-trip route in the worst case. Besides, this observation has also an effect on the frequency of a service provided between any pair of physical nodes. Given that the distance is half of the schedule length for any pair of nodes, then the service between these nodes can only occur once in each schedule.  

As stated above services can be performed by a resource assuming that this is a single type asset and exactly one unit is assigned to any service consistently with the previous studies. Human resources such as crews and workers, transporting vehicles such as trucks, aircraft, rail-cars and ships, loading/unloading vehicles such as forklifts, cranes, and pallet jacks are examples of assets. In this paper, we also assume that assets are divided into two sets, owned assets and assets that can be acquired either through leasing or buying at the beginning of the planning horizon. Let $V_{1}$ represent owned assets assuming that this is a finite set from the realistic standpoint. The assets that can be acquired are symbolized by $V_{2}$ which is limited by a relatively high $\bar{V_{2}}$; even though, there is no limit to acquire additional assets practically except a limited budget. The union of these two sets, $V$, comprises the fleet which is assigned to services and routed through the network during the planning horizon.    

Any amount of flow between any origin-destination pair demanded to transport is called a commodity and let $K$ represent the set of all commodities. Any element of $K$ is associated with a volume $w^{k}$, an origin node $O^{'}_{k}$ and a destination node $D^{'}_{k}$ as well as a release date and a deadline in terms of the period index. Thus, let parameters $O_{k}$ and $D_{k}$ represent the origin and destination of a commodity on the time-space network. The transportation activity for any commodity may start at any time after the beginning of the period of the release date and must be completed before the end of the period of the deadline.

The option of earliness/tardiness affect all parameters regarding the commodities. We assume that any commodity can be delivered as much as one period earlier or later than its original deadline. To be able to incorporate this option into the formulation, we define two more commodities as \textit{dummy commodities} which corresponds to each original commodity within $K$. Dummy commodities represent early and late deliveries of the original commodity, and their origin and destination nodes are denoted by $O^{2}_{k}$ and $D^{2}_{k}$ on the time-space network. The time periods correspond to origin and destination of dummy commodities are $t_{O_{k}}-1$ and $t_{D_{k}}-1$ for the former and $t_{O_{k}}+1$ and $t_{D_{k}}+1$ for the latter. The combination of dummy commodities and the corresponding original commodity is called as \textit{transformed commodity} and set of original commodities, $K$, is replaced with this concept in the proposed problem formulation. Let $L$ represent a set of dummy commodities and let union set of $K \cup L$ represent the aggregated set of original and dummy commodities, namely transformed commodities. As indicated implicitly, the union set of $K \cup L$ includes three times more element than the set of original commodities, $K$, since there are three transformed commodities per original one in the problem.   



The original and transformed commodities are related to each other through an incident matrix, $\alpha_{kl}$ which equals to 1, if transformed commodity $ l \in K \cup L $ is incident to original commodity $ k \in K $; 0, otherwise. With consideration of transformed commodities, the volume of a transformed commodity is denoted by $w_{2}^{k}$, and the corresponding type is represented by $q^{k}$ either as early, original or tardy.


The operating cost of such a service network is considered in three parts. The first part includes fixed cost of operating an asset in the planning horizon, for which we define $f$ and $g$ to represent the fixed cost of operating an owned asset and a leased or an acquired asset, respectively. All costs associated with the crew, depreciation, etc. are assumed to be included in the fixed cost of operating assets. The routing cost of commodities on the network is involved in the second part. Let $c^{k}_{ij}$ denote variable cost of routing transformed commodity $k$ on any service arc $(i,j)$ including offered and outsourced services to cover fuel cost en route and handling at the terminals. The final part consists of penalty for delivering commodities early and tardy, a fee denoted by $r_{e}$ and $r_{l}$ incur for each unit of flow picked up early or delivered late, respectively.    

In addition to routing of assets and flow of commodities, as defined earlier the capacity-demand balancing problem in service network design for a consolidation-based freight forwarder in this paper also includes the following decisions: (1) determine whether to add additional assets into the fleet, if that is the case, decide routes of them, (2) whether to pick up any commodity earlier or deliver any commodity late, in that case which commodities are delivered early/late, and (3) whether to choose an outsourced service in case of not choosing other options to handle the excess demand and decide which commodities are routed through selected outsourced services.

For the sake of readers' convenience, all the notation given so far are also presented in Table~\ref{tab:Notation}. 


\begin{table}[htbp]
	\centering
	\scriptsize
	\caption {Problem notation.} \label{tab:Notation}
	{
		\begin{tabular}{l l}
			\midrule	
			\bf Notation	&	\bf Description	\\	\midrule	
			\multicolumn{2}{l}{\bf Sets:} \\		
			$G=(N,A)$             & Time-space network.\\
			$G^{'}=(N^{'},A^{'})$             & Physical (static) network.\\
			$N$             & Set of nodes on the time-space network.\\
			$N^{'}$             & Set of nodes on the physical network.\\
			$A_{h}$             & Set of holding arcs on the time-space network.\\
			$A_{s}$             & Set of arcs representing services on the time-space network.\\
			$A_{o}$             & Set of arcs representing outsourced services on the time-space network.\\
			$A=A_{h} \cup A_{s} \cup A_{o} $             & Set of all arcs on the time-space network.\\
			$A^{'}$             & Set of all arcs on the physical network.\\
			$T$             & Set of time periods.\\
			$V_{1}=\{1,..,\bar{V_{1}}\}$             & Set of owned assets.\\
			$V_{2}=\{1,..,\bar{V_{2}}\}$             & Set of assets that can be acquired as additional units.\\
			$V=V_{1} \cup V_{2}$             & Set of all assets such that $V=\{1,..,\bar{V_{1}},\bar{V_{1}}+1,..,\bar{V_{1}}+\bar{V_{2}}\}$.\\
			$K$             & Set of original commodities.\\
			\multirow{2}{*}{\bf $L$} 	& Set of dummy commodities representing early and tardy versions of \\
			& each corresponding original commodity. \\
			$K \cup L$             & Joint set of transformed commodities. \\
			\\
			\multicolumn{2}{l}{\bf Parameters:} \\
			\multirow{2}{*}{\bf $d_{ij}$}     & Distance between node $i \in N^{'}$ and node $j \in N^{'}$ in terms of number of periods \\
			& on the physical network. \\
			$w^{k}$    & Volume of commodity $k \in K$ that needs to be transported.\\
			$w_{2}^{k}$    & Volume of transformed commodity $k \in K \cup L$ that needs to be transported.\\
			$t_{i}$		& Period index of node $i \in N$. \\
			$O^{'}_{k}$ & Origin node of commodity $k \in K$ on static network, $G^{'}=(N^{'},A^{'})$.\\
			$O_{k}$ & Origin node of commodity $k \in K$ on time-space network, $G=(N,A)$.\\
			\multirow{2}{*}{\bf $O_{k}^{2}$} & Origin node of transformed commodity $k \in K \cup L$ on time-space network, \\
			& $G=(N,A)$.\\
			$D^{'}_{k}$ & Destination node of commodity $k \in K$ on static network, $G^{'}=(N^{'},A^{'})$. \\
			$D_{k}$ & Destination node of commodity $k \in K$ on time-space network, $G=(N,A)$. \\
			\multirow{2}{*}{\bf $D_{k}^{2}$} & Destination node of transformed commodity $k \in K \cup L$ on time-space network, \\
			& $G=(N,A)$.\\
			$c^{k}_{ij}$   & Cost of transporting one unit of transformed commodity $k \in K \cup L$ on arc $(i,j) \in A$. \\
			$f$       & Fixed cost of operating a unit of asset. \\
			$g$       & Fixed cost of acquiring/leasing and operating an additional unit of asset. \\
			$u_{ij}$   & Capacity of service operated on arc $(i,j) \in A$. \\
			$b^{k}_{ij}=\min\{w^{k},u_{ij}\}$   & Additional parameter defined to obtain stronger formulation. \\
			$q^{k}$   & Type of transformed commodity $k \in K \cup L$ either as early, original or tardy. \\
			$r_{e}/r_{l}$   & Penalty for cost of transporting one unit of dummy commodity early/tardy. \\
			\multirow{2}{*}{\bf $\alpha_{kl}$} & 1, if transformed commodity $ l \in K \cup L $ is incident to original commodity $ k \in K $; \\ 
			& 0, otherwise. \\
            \multirow{2}{*}{\bf $\beta_{kt}$} & 1, if transformed commodity $ k \in K \cup L $ might be in transit in period $ t \in T $; \\ 
			& 0, otherwise. \\
			\midrule	
		\end{tabular}
	}
\end{table}

In order to capture all decisions regarding capacity-demand balancing problem on a service network design, we propose the following sets of decision variables:

\bigskip
\begin{tabular}{lp{13cm}}
$ y^{v}_{ij}  =  \left\{
\begin{array}{ll}
1, & \text{if arc $ (i,j) \in A_{h} \cup A_{s} $ is selected and operated by asset $ v \in V $} \\
0, & \text{otherwise}
\end{array} \right.  $ 
\end{tabular}

\bigskip
\begin{tabular}{lp{13cm}}
$ \delta_{v} =  \left\{
\begin{array}{ll}
1, & \text{if asset $ v \in V $  is utilized on an activity} \\
0, & \text{otherwise}
\end{array} \right.  $ 
\end{tabular}

\bigskip
\begin{tabular}{lp{13cm}}
$ p_{k}  =  \left\{
\begin{array}{ll}
1, & \text{if transformed commodity $ k \in K \cup L $ is selected to be delivered}\\
0, & \text{otherwise}
\end{array} \right.  $ 
\end{tabular}

\bigskip
\begin{tabular}{lp{13cm}}
$ s_{ij}^{k}  =  \left\{
\begin{array}{ll}
1, & \text{if the outsourced service on arc $ (i,j) \in A_{o} $ is selected for delivering}  \\
   &  \text{transformed commodity $ k \in K \cup L$}   \\
0, & \text{otherwise}
\end{array} \right.  $ 
\end{tabular}

\bigskip
\begin{tabular}{lp{13cm}}
$x^{k}_{ij}$ = & Amount of flow of transformed commodity $ k \in K \cup L $ routed on arc $ (i,j) \in A $.
\end{tabular}

\bigskip
Depending on given notation and defined decision variables, capacity-demand balancing problem on service network design can be formulated as follows:
\begin{align}
\mbox{\textbf{($P$)}}
\qquad Min \quad   &      \sum_{v \in V:v \leq V_{1}} f \delta_{v} + \sum_{v \in V:v \geq V_{1}+1} g \delta_{v}
                 + \sum_{(i,j) \in A_{h} \cup A_{s}} \sum_{k \in K \cup L: q_{k}=2} c^{k}_{ij} x^{k}_{ij} 
                 \nonumber \\
            &    \qquad
                 + \sum_{(i,j) \in A_{h} \cup A_{s}} \sum_{k \in K \cup L: q_{k}=1} r_{e} c^{k}_{ij} x^{k}_{ij} 
                 + \sum_{(i,j) \in A_{h} \cup A_{s}} \sum_{k \in K \cup L: q_{k}=3} r_{l} c^{k}_{ij} x^{k}_{ij}
                 \nonumber \\
            &    \qquad
                 + \sum_{(i,j) \in A_{o}}\sum_{k \in K \cup L} c^{k}_{ij} s^{k}_{ij}
                                                          \label{SND-Obj} \\
s. \: t. \quad  & \sum_{(i,j) \in A_{h} \cup A_{s}, \tilde{t_{i}}, \tilde{t_{j}} \in T: \tilde{t_{i}} \leq t < \tilde{t_{j}}} x^{k}_{ij} \leq 0 \qquad \forall \: k \in K \cup L, t \in T: \beta_{kt}=0    \label{SND-Com_not_in_transit} \\
                & \sum_{(i,j) \in A_{h} \cup A_{s}, \tilde{t_{i}}, \tilde{t_{j}} \in T: \tilde{t_{i}} \leq t < \tilde{t_{j}}} y^{v}_{ij} - \delta_{v} = 0 \qquad \forall \: v \in V, t \in T     \label{SND-Asset_for_one_activity} \\
                & \sum_{j \in N:(i,j) \in A_{h} \cup A_{s}} y^{v}_{ij} - \sum_{j \in N:(j,i) \in A_{h} \cup A_{s}} y^{v}_{ji}=0 \qquad \forall \: i \in N, v \in V
                    \label{SND-Design_balance} \\
                & \sum_{v \in V} y^{v}_{ij} \leq 1 \qquad \forall \: (i,j) \in A_{s}
                	\label{SND-One_asset_for_each_activity} \\
                & \sum_{l \in K \cup L: \alpha_{kl}=1} p_{l} \geq 1 \qquad \forall \: k \in K
                \label{SND-One_trans_com_selected} \\	  
                & \sum_{j \in N:(i,j) \in A} x^{k}_{ij} - \sum_{j \in N:(j,i) \in A} x^{k}_{ji} = \begin{cases}
                w^{k} p_{k}, & \text{$i=O^{k}_{2}$}.\\
                -w^{k} p_{k}, & \text{$i=D^{k}_{2}$}.\\
                0, & \text{otherwise}.
                \end{cases} \qquad  \forall \: i \in N, k \in K \cup L \label{SND-FlowBalance} \\             
                & \sum_{k \in K \cup L} x^{k}_{ij} - \sum_{v \in V} u_{ij} y^{v}_{ij}
                    \leq 0 \qquad \forall \: (i,j) \in A_{h} \cup A_{s}
                    \label{SND-Forcing_and_capacity_weak} \\   
                & x^{k}_{ij} - \sum_{v \in V} b^{k}_{ij} y^{v}_{ij}
                \leq 0 \qquad \forall \: (i,j) \in A_{h} \cup A_{s}, k \in K \cup L
                \label{SND-Forcing_and_capacity_strong} \\
                 & x^{k}_{ij} - \sum_{v \in V} u_{ij} s^{k}_{ij}
                \leq 0 \qquad \forall \: (i,j) \in A_{o}, k \in K \cup L
                \label{SND-Forcing_and_Capacity_Outsourced} \\
                & y^{v}_{ij} \in \{0,1\} \qquad \forall \: (i,j) \in A_{h} \cup A_{s}, v \in V
                    \label{SND-Domain-y} \\
                & \delta_{v} \in \{0,1\} \qquad \forall \: v \in V
                \label{SND-Domain-delta} \\
                & x^{k}_{ij} \geq 0 \qquad \forall \: (i,j) \in A,  k \in K \cup L
                \label{SND-Domain-x} \\
                & s^{k}_{ij} \in \{0,1\} \qquad \forall \: (i,j) \in A_{o},  k \in K \cup L
                    \label{SND-Domain-s} \\
                & p_{k} \in \{0,1\} \qquad \forall \: k \in K \cup L
                    \label{SND-Domain-u}
\end{align}

The objective function \eqref{SND-Obj} accounts for the total cost over the planning horizon including the following five terms: (i) fixed cost of using owned assets, (ii) fixed cost of acquiring and using additional assets, (iii) cost of transporting commodities on offered services, (iv) cost of transporting commodities on offered services earlier or later with a penalty for earliness/tardiness and, (v) cost of transporting commodities on outsourced services respectively.

\begin{figure}[htbp]
	\centering
	\fbox{\includegraphics[width=0.8\textwidth]{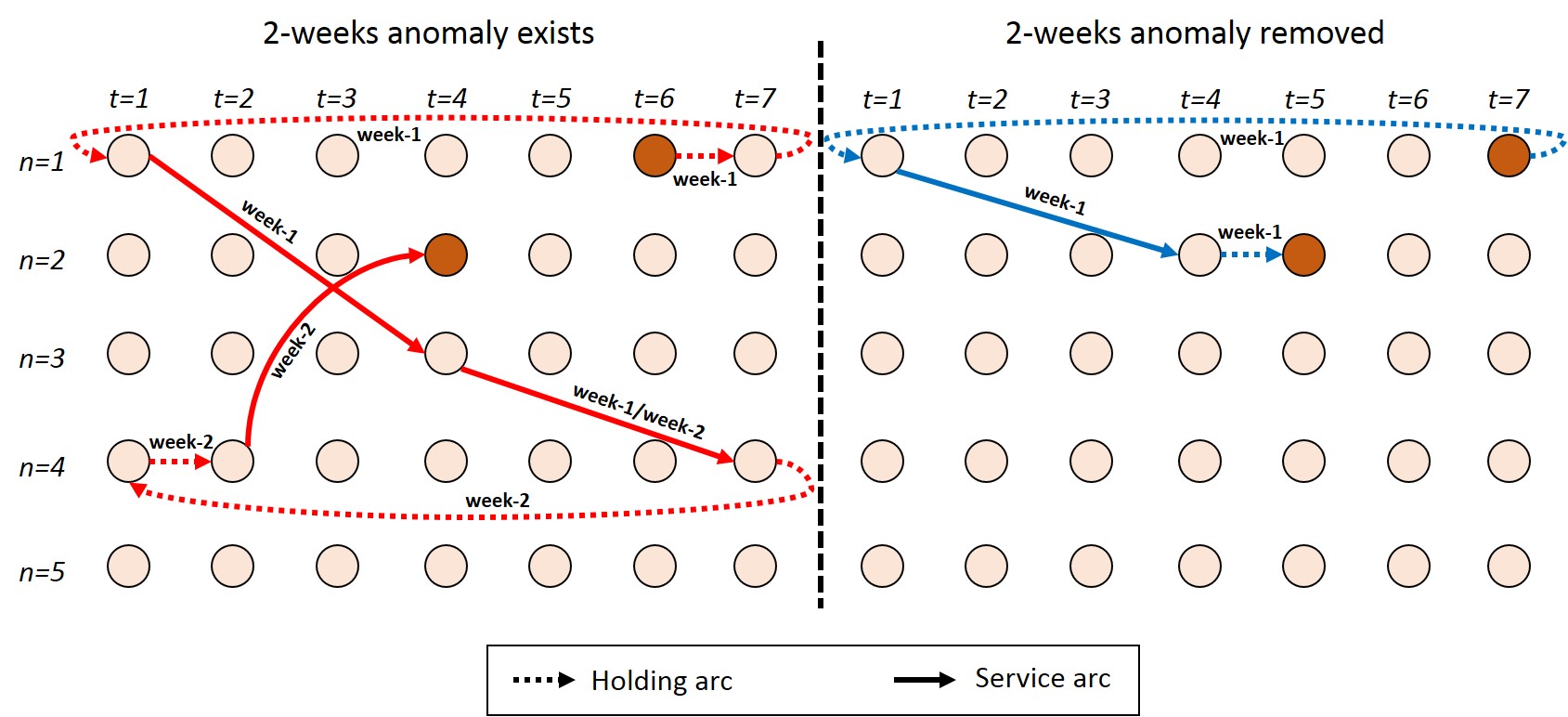}}
	\caption{Illustration of 2-weeks anomaly in commodity flows.}
	\label{fig:anomaly}
\end{figure}

Constraints~\eqref{SND-Com_not_in_transit} are added into the formulation in this paper and do not exist in the formulations proposed in the SND literature so far. The reason for adding these depend on our observation in the preliminary analysis, and we aim to prevent an anomaly that may exist in commodity flows, which is illustrated in Figure~\ref{fig:anomaly}. This anomaly rarely occurs in a few instances for few commodities but violates the overall feasibility of the problem literally. To refer it simply, we prefer to call it as \emph{2-weeks anomaly}. As indicated in Figure~\ref{fig:anomaly}, the transformed commodity (original-type) shown on the left network in the figure (originated at the sixth period at node-1 and destined to the fourth period at node-2) is delivered in almost two weeks from its release date. When we add the constraint set~\eqref{SND-Com_not_in_transit}, 2-weeks anomaly is eliminated from the solution and the tardy-type transformed commodity (originated at the seventh period at node-1 and destined to the fifth period at node-2) that is incident to the same original commodity is delivered instead as shown on the right network in the figure. To be able to formulate this constraint set, we define a new $0-1$ parameter ($\beta_{kt}$) which takes the value of one for the period $t \in T$ when transformed commodity $k \in K \cup L$ would be in transit in a feasible solution based on commodity's time-window. Thus, we are able to restrict the values of decision variable of commodity flows ($x^{k}_{ij}$) to zero on an arc that takes place when the commodity is supposed to be not in transit.  

Constraints \eqref{SND-Asset_for_one_activity} assure that in each period, an owned or added asset must be assigned to only one activity (holding or service) if it is utilized. \cite{Andersen&Crainic&Christiansen:2009_TRC} stated that this constraint imposes a maximum route-length requirement for a schedule. Even though, this constraint set is valid theoretically, the particular condition given under sum notation makes impossible to write it down mathematically due to few arcs wrapping around the planning horizon. We propose to replace constraint set \eqref{SND-Asset_for_one_activity} by \eqref{SND-Asset_for_one_activity_C1}-\eqref{SND-Asset_for_one_activity_C3} based on new sets defined in Table~\ref{tab:Notation2}.

\begin{table}[htbp]
	\centering
	\scriptsize
	\caption {Notation for replacing constraint set \eqref{SND-Asset_for_one_activity}.} \label{tab:Notation2}
	{
		\begin{tabular}{l l}
			\midrule
			\bf Notation	&	\bf Description	\\	\midrule			
			\multicolumn{2}{l}{\bf Additional Sets:} \\
			$A^{'}_{h}/A^{'}_{s}$ & Holding/service arcs that do not wrap around the planning horizon (regular), $ t_{origin} < t_{destination} $.\\
			$A^{''}_{h}/A^{''}_{s}$ & Holding/service arcs that DO wrap around the planning horizon (circular), $ t_{destination} < t_{origin} $. \\
			$A^{'}_{o}$ & Arcs of outsourced services, that do not wrap around the planning horizon (regular), $ t_{origin} < t_{destination} $.\\
			$A^{''}_{o}$ & Arcs of outsourced services, that DO wrap around the planning horizon (circular), $ t_{destination} < t_{origin} $. \\
			$T_{1} $ & Time periods which are only visited by regular holding/service arcs, i.e. $T_{1}={3,4}$.\\
			$T_{2} $ & Time periods which exist in circular arcs as origin nodes, i.e. $T_{2}={5,6,7}$. \\
			$T_{3} $ & Time periods which exist in circular arcs as destination nodes, i.e. $T_{3}={1,2}$. \\
			\midrule	
		\end{tabular}
	}
\end{table}
\begin{align}
& \sum_{(i,j) \in A_{h} \cup A_{s}, \tilde{t_{i}}, \tilde{t_{j}} \in T: \tilde{t_{i}} \leq t < \tilde{t_{j}}} y^{v}_{ij} - \delta_{v} = 0 \qquad \forall \: v \in V, t \in T_{1}  \label{SND-Asset_for_one_activity_C1} \\
& \sum_{(i,j) \in A^{'}_{h} \cup A^{'}_{s}, \tilde{t_{i}}, \tilde{t_{j}} \in T: \tilde{t_{i}} \leq t < \tilde{t_{j}}} y^{v}_{ij} + \sum_{(i,j) \in A^{''}_{h} \cup A^{''}_{s}, \tilde{t_{i}}, \tilde{t_{j}} \in T: \tilde{t_{i}} \leq t < \tilde{t_{j}}+ |T|} y^{v}_{ij} - \delta_{v} = 0 \qquad \forall \: v \in V, t \in T_{2}  \label{SND-Asset_for_one_activity_C2} \\
& \sum_{(i,j) \in A^{'}_{h} \cup A^{'}_{s}, \tilde{t_{i}}, \tilde{t_{j}} \in T: \tilde{t_{i}} \leq t < \tilde{t_{j}}} y^{v}_{ij} + \sum_{(i,j) \in A^{''}_{h} \cup A^{''}_{s}, \tilde{t_{i}}, \tilde{t_{j}} \in T: \tilde{t_{i}} \leq t+|T| , t < \tilde{t_{j}}} y^{v}_{ij} - \delta_{v} = 0 \qquad \forall \: v \in V, t \in T_{3}  \label{SND-Asset_for_one_activity_C3}
\end{align}

The balance of assets are satisfied through constraints \eqref{SND-Design_balance} in which the incoming number of assets must be equal to outgoing assets for each node on the time-space network. 
Constraints \eqref{SND-One_asset_for_each_activity} enforce that only one asset is assigned to a service at most. 
Constraints \eqref{SND-One_trans_com_selected} ensure that at least one of the transformed commodity which is associated with the original correspondent commodity must be selected to be delivered. Since the fact that the objective is minimized, only one of a transformed commodity out of three will be selected for delivery, obviously.
Inequalities \eqref{SND-FlowBalance} are flow balance constraints for the transformed commodities in case of the corresponding one is selected to be delivered. Otherwise, this set of constraints become redundant for two out of three transformed commodities that are incident to a particular original commodity. 
Constraints \eqref{SND-Forcing_and_capacity_weak} and \eqref{SND-Forcing_and_capacity_strong} are the weak and strong capacity and forcing constraints for the arcs of offered services and holding arcs, respectively. The flow on the arcs of outsourced services is restricted through constraints \eqref{SND-Forcing_and_Capacity_Outsourced} if the corresponding arc is selected. Note that, strong version of capacity and forcing constraint \eqref{SND-Forcing_and_capacity_strong} is redundant for MIP formulation as \cite{Andersen&Crainic&Christiansen:2009_TRC} stated as well. 
Finally, \eqref{SND-Domain-y}-\eqref{SND-Domain-u} are domain constraints.

\subsection{Valid Inequalities}
\label{subsec-valid_inequal}
To enhance the performance of the formulation, we also propose few sets of valid inequalities which were discovered during the preliminary analysis of the problem instances. For this reason, we would like to explain the idea behind the generation of valid cuts before presenting them in this section. 

To clarify the idea perfectly, we prefer to apply the procedure on a sample instance for illustration. Table~\ref{tab:SampleComData} presents the values for commodity-related parameters of the ten original commodities (OCs) with five nodes and seven periods. The number of owned and acquirable assets are seven and five, respectively. The other details of the test instances will be explained later in Subsection of \textit{Instance Generation} (\S\ref{subsec-instance-gen}). 


In Table~\ref{tab:SampleComData}, the second and third column show the origin and destination node of the corresponding OC ($k \in K$) on the time-space network. The volume of each commodity is presented in the fourth column. The origin and destination node on the physical network, and release and due date of any commodity in terms of periods are given in the fifth through eighth columns. For instance, commodity-1 is originated at node-2 on the physical network in the second period and destined to node-1 with a due date of the fifth period. The same information for these parameters is also presented regarding transformed commodities (TCs) ($l \in K \cup L$) in Table~\ref{tab:SampleTComData} included with the type parameter of each transformed commodity in the fifth column.

\begin{table}[htbp]
	\begin{minipage}{.5\linewidth}
			\centering
		\scriptsize
		\caption {Original commodities.} \label{tab:SampleComData}
		{
			\renewcommand{\arraystretch}{0.9}
			\begin{tabular}{l c c c c c c c}
				\midrule
				\bf $k$	&	\bf $O_{k}$	&	\bf $D_{k}$	&	\bf $w^{k}$	&	\bf $O^{'}_{k}$	&	\bf $D^{'}_{k}$	&	\bf $t_{O_{k}}$	&	\bf $t_{D_{k}}$	\\	\midrule
				1	&	9	&	5	&	1	&	2	&	1	&	2	&	5	\\	
				2	&	17	&	13	&	1	&	3	&	2	&	3	&	6	\\	
				3	&	19	&	3	&	1	&	3	&	1	&	5	&	3	\\	
				4	&	11	&	28	&	1	&	2	&	4	&	4	&	7	\\	
				5	&	26	&	17	&	1	&	4	&	3	&	5	&	3	\\	
				6	&	7	&	12	&	1	&	1	&	2	&	7	&	5	\\	
				7	&	18	&	35	&	1	&	3	&	5	&	4	&	7	\\	
				8	&	7	&	25	&	1	&	1	&	4	&	7	&	4	\\	
				9	&	31	&	19	&	1	&	5	&	3	&	3	&	5	\\	
				10	&	29	&	24	&	1	&	5	&	4	&	1	&	3	\\	
				
				\hline
			\end{tabular}
		}
				
		\bigskip		
		
		\centering
		\scriptsize
		\caption {\small Transformed commodities.} \label{tab:SampleTComData}
		{
			\renewcommand{\arraystretch}{0.9}
			\begin{tabular}{l c c c c c c c c}
				\midrule
				\bf $l$	&	\bf $O_{k}$	&	\bf $D_{k}$	&	\bf $w^{k}$	&	\bf $q^{k}$	&	\bf $O^{'}_{k}$	&	\bf $D^{'}_{k}$	&	\bf $t_{O_{k}}$	&	\bf $t_{D_{k}}$	\\	\midrule
				1	&	8	&	4	&	1	&	2	&	2	&	1	&	1	&	4	\\	
				2	&	9	&	5	&	1	&	1	&	2	&	1	&	2	&	5	\\	
				3	&	10	&	6	&	1	&	2	&	2	&	1	&	3	&	6	\\	
				4	&	16	&	12	&	1	&	2	&	3	&	2	&	2	&	5	\\	
				5	&	17	&	13	&	1	&	1	&	3	&	2	&	3	&	6	\\	
				6	&	18	&	14	&	1	&	2	&	3	&	2	&	4	&	7	\\	
				7	&	18	&	2	&	1	&	2	&	3	&	1	&	4	&	2	\\	
				8	&	19	&	3	&	1	&	1	&	3	&	1	&	5	&	3	\\	
				9	&	20	&	4	&	1	&	2	&	3	&	1	&	6	&	4	\\	
				10	&	10	&	27	&	1	&	2	&	2	&	4	&	3	&	6	\\	
				11	&	11	&	28	&	1	&	1	&	2	&	4	&	4	&	7	\\	
				12	&	12	&	22	&	1	&	2	&	2	&	4	&	5	&	1	\\	
				13	&	25	&	16	&	1	&	2	&	4	&	3	&	4	&	2	\\	
				14	&	26	&	17	&	1	&	1	&	4	&	3	&	5	&	3	\\	
				15	&	27	&	18	&	1	&	2	&	4	&	3	&	6	&	4	\\	
				16	&	6	&	11	&	1	&	2	&	1	&	2	&	6	&	4	\\	
				17	&	7	&	12	&	1	&	1	&	1	&	2	&	7	&	5	\\	
				18	&	1	&	13	&	1	&	2	&	1	&	2	&	1	&	6	\\	
				19	&	17	&	34	&	1	&	2	&	3	&	5	&	3	&	6	\\	
				20	&	18	&	35	&	1	&	1	&	3	&	5	&	4	&	7	\\	
				21	&	19	&	29	&	1	&	2	&	3	&	5	&	5	&	1	\\	
				22	&	6	&	24	&	1	&	2	&	1	&	4	&	6	&	3	\\	
				23	&	7	&	25	&	1	&	1	&	1	&	4	&	7	&	4	\\	
				24	&	1	&	26	&	1	&	2	&	1	&	4	&	1	&	5	\\	
				25	&	30	&	18	&	1	&	2	&	5	&	3	&	2	&	4	\\	
				26	&	31	&	19	&	1	&	1	&	5	&	3	&	3	&	5	\\	
				27	&	32	&	20	&	1	&	2	&	5	&	3	&	4	&	6	\\	
				28	&	35	&	23	&	1	&	2	&	5	&	4	&	7	&	2	\\	
				29	&	29	&	24	&	1	&	1	&	5	&	4	&	1	&	3	\\	
				30	&	30	&	25	&	1	&	2	&	5	&	4	&	2	&	4	\\	
				
				\hline
			\end{tabular}
		}
	\end{minipage}%
	\begin{minipage}{.5\linewidth}
		\centering
		\scriptsize
		\caption {\small OCs' release times and due dates in tabular format.} \label{tab:ComDataTabular}
		{
			\renewcommand{\arraystretch}{0.9}
			\begin{tabular}{l c c c c c c c}
				\midrule
				\bf $k$	&	\bf $T1$	&	\bf $T2$	&	\bf $T3$	&	\bf $T4$	&	\bf $T5$	&	\bf $T6$	&	\bf $T7$	\\	\midrule
				1	&		&	1	&	1	&	1	&	1	&		&		\\	
				2	&		&		&	1	&	1	&	1	&	1	&		\\	
				3	&	1	&	1	&	1	&		&	1	&	1	&	1	\\	
				4	&		&		&		&	1	&	1	&	1	&	1	\\	
				5	&	1	&	1	&	1	&		&	1	&	1	&	1	\\	
				6	&	1	&	1	&	1	&	1	&	1	&		&	1	\\	
				7	&		&		&		&	1	&	1	&	1	&	1	\\	
				8	&	1	&	1	&	1	&	1	&		&		&	1	\\	
				9	&		&		&	1	&	1	&	1	&		&		\\	
				10	&	1	&	1	&	1	&		&		&		&		\\				
				\hline
			\end{tabular}
		}
	\bigskip
	\bigskip
	
	\centering
	\scriptsize
	\caption {\small TCs' release times and due dates in tabular format.} \label{tab:TComDataTabular}
	{
		\renewcommand{\arraystretch}{0.9}
		\begin{tabular}{l c c c c c c c c}
			\midrule
			\bf $l$	&	\bf $k$	&	\bf $T1$	&	\bf $T2$	&	\bf $T3$	&	\bf $T4$	&	\bf $T5$	&	\bf $T6$	&	\bf $T7$	\\	\midrule
			1	&	\multirow{3}{*}{1}	&	1	&	1	&	\bf 1	&	\bf 1	&		&		&		\\	
			2	&		&		&	1	&	\bf 1	&	\bf 1	&	1	&		&		\\	
			3	&		&		&		&	\bf 1	&	\bf 1	&	1	&	1	&		\\	\midrule
			4	&	\multirow{3}{*}{2}	&		&	1	&	1	&	\bf 1	&	\bf 1	&		&		\\	
			5	&		&		&		&	1	&	\bf 1	&	\bf 1	&	1	&		\\	
			6	&		&		&		&		&	\bf 1	&	\bf 1	&	1	&	1	\\	\midrule
			7	&	\multirow{3}{*}{3}	&	\bf 1	&	\bf 1	&		&	1	&	1	&	\bf 1	&	\bf 1	\\	
			8	&		&	\bf 1	&	\bf 1	&	1	&		&	1	&	\bf 1	&	\bf 1	\\	
			9	&		&	\bf 1	&	\bf 1	&	1	&	1	&		&	\bf 1	&	\bf 1	\\				
			\hline
		\end{tabular}
	}
	
	\bigskip
	\bigskip
	
	\centering
	\scriptsize
	\caption {Asset requirements based on OCs.} \label{tab:SumComDataInTransit}
	{
		\renewcommand{\arraystretch}{0.9}
		\begin{tabular}{l c c c c c c c}
			\midrule
			\bf $k$	&	\bf $T1$	&	\bf $T2$	&	\bf $T3$	&	\bf $T4$	&	\bf $T5$	&	\bf $T6$	&	\bf $T7$	\\	\midrule
			1	&		&		&	1	&	1	&		&		&		\\	
			2	&		&		&		&	1	&	1	&		&		\\	
			3	&	1	&	1	&		&		&		&	1	&	1	\\	
			4	&		&		&		&		&	1	&	1	&		\\	
			5	&	1	&	1	&		&		&		&	1	&	1	\\	
			6	&	1	&	1	&	1	&	1	&		&		&		\\	
			7	&		&		&		&		&	1	&	1	&		\\	
			8	&	1	&	1	&	1	&		&		&		&		\\	
			9	&		&		&		&	1	&		&		&		\\	
			10	&		&	1	&		&		&		&		&		\\	\midrule
			Total	&	4	&	5	&	3	&	4	&	3	&	4	&	2	\\	
			
			\hline
		\end{tabular}
	}

	\end{minipage} 
\end{table}


The time window of each original commodity based on corresponding release times and due dates ($k \in K$) are mapped on a week-based tabular format in Table~\ref{tab:ComDataTabular} for better understanding because the main idea behind the valid inequalities comes from this mapping. Any commodity is expected to be assigned to any asset (resource) in a period where a value of one appears in the corresponding column. The tabular format is enhanced by considering transformed commodities in Table~\ref{tab:TComDataTabular} where only nine transformed commodities which correspond to original commodities 1, 2, and 3, are listed in each three-lines row as a sample. In the wide top row, the first line shows the mapping of time window belongs to transformed commodity-1 which is the early version of original commodity-1 on a week. The second row (transformed commodity-2) corresponds the original commodity itself, and the third row is the tardy version of original commodity-1. The columns in which the ones are typed in bold depict the intersecting time-periods of three transformed commodities; thus the corresponding original commodity occupies a resource in these periods no matter of which one of them is selected for delivery. When we apply the same mapping for all original commodities and add all ones in the same column up, we can predict the number of assets utilized in each period to transport all commodities (see Table~\ref{tab:SumComDataInTransit}).




Based on this discussion, the additional few parameters are defined for the formulation of the valid inequalities and presented in Table~\ref{tab:NotationVI}.

\begin{table}[htbp]
	\centering
	\small
	\caption {Notation for Valid Inequalities.} \label{tab:NotationVI}
	{
		\begin{tabular}{l l}
			\midrule
			\bf Notation	&	\bf Description	\\	\midrule			
			\multicolumn{2}{l}{\bf Additional Parameters:} \\
			$\Phi_{t}$ & Number of assets that needs to be utilized for period $t \in T$ within the planning horizon.\\
			$\Gamma=\min\limits_{t \in T}\{\Phi_{t}\}$ & Minimum number of assets that needs to be utilized for the entire planning horizon.\\
			\midrule	
		\end{tabular}
	}
\end{table}

Regarding the idea of predicting the number of assets utilized from the period-based mapping, we propose two sets of valid inequalities based on binary decision variables and given additional parameters as follows:
\begin{align}
& \sum_{v \in V} \delta_{v} \geq \Gamma \label{SND-VI11-MinNbAsset} \end{align}
\begin{subequations}\label{SND-VI4X-NbAssetperPeriod}
\begin{align}
& \sum_{v \in V} \sum_{(i,j) \in A_{h} \cup A_{s}, \tilde{t_{i}}, \tilde{t_{j}} \in T: 
\tilde{t_{i}} \leq t < \tilde{t_{j}}} y^{v}_{ij} 
+ \sum_{k \in K \cup L} \sum_{(i,j) \in A_{o}, \tilde{t_{i}}, \tilde{t_{j}} \in T: 
\tilde{t_{i}} \leq t < \tilde{t_{j}}} s^{k}_{ij}  \geq \Phi_{t} \qquad \forall \: t \in T_{1} \label{SND-VI4X-NbAssetperPeriod_C1} \\
& \sum_{v \in V} \bigg( \sum_{(i,j) \in A^{'}_{h} \cup A^{'}_{s}, \tilde{t_{i}}, \tilde{t_
{j}} \in T: \tilde{t_{i}} \leq t < \tilde{t_{j}}} y^{v}_{ij} + \sum_{(i,j) \in A^{''}_
{h} \cup A^{''}_{s} , \tilde{t_{i}}, \tilde{t_{j}} \in T: tilde{t_{i}} \leq t < \tilde
{t_{j}}+ |T|} y^{v}_{ij} \bigg) \nonumber \\
& + \sum_{k \in K \cup L} \bigg( \sum_{(i,j) \in A^{'}_{o}, \tilde
{t_{i}}, \tilde{t_{j}} \in T: \tilde{t_{i}} \leq t < \tilde{t_{j}}} s^{k}_{ij} + \sum_
{(i,j) \in A^{''}_{o}, \tilde{t_{i}}, \tilde{t_{j}} \in T:\tilde{t_{i}} \leq t < \tilde{t_{j}}+ |T|} s^{k}_{ij} \bigg) \geq \Phi_{t} \qquad \forall \: t \in T_{2} \label{SND-VI4X-NbAssetperPeriod_C2} \\
& \sum_{v \in V} \bigg( \sum_{(i,j) \in A^{'}_{h} \cup A^{'}_{s}, \tilde{t_{i}}, \tilde{t_
{j}} \in T: \tilde{t_{i}} \leq t < \tilde{t_{j}}} y^{v}_{ij} + \sum_{(i,j) \in A^{''}_
{h} \cup A^{''}_{s}, \tilde{t_{i}}, \tilde{t_{j}} \in T: tilde{t_{i}} \leq t+|T|, t < \tilde
{t_{j}}} y^{v}_{ij} \bigg) \nonumber \\
& + \sum_{k \in K \cup L} \bigg( \sum_{(i,j) \in A^{'}_{o}, \tilde
{t_{i}}, \tilde{t_{j}} \in T: \tilde{t_{i}} \leq t < \tilde{t_{j}}} s^{k}_{ij} + \sum_
{(i,j) \in A^{''}_{o}, \tilde
{t_{i}}, \tilde{t_{j}} \in T: tilde{t_{i}} \leq t+|T|, t < \tilde
{t_{j}}} s^{k}_{ij} \bigg) \geq \Phi_{t} \qquad \forall \: t \in T_{3} \label{SND-VI4X-NbAssetperPeriod_C3}
\end{align}
\end{subequations}

In the first set of inequalities, we aim to put a lower bound on the number of assets required from the perspective of the minimum number of resources. Inequality \eqref{SND-VI11-MinNbAsset} gives a loose lower bound on the number of assets that need to be utilized in the optimal solution by forcing it with the minimum predicted number of assets over all periods. The second set of valid inequalities \eqref{SND-VI4X-NbAssetperPeriod} are focused on obtaining a lower bound for a total number of transportation resources required as either owned assets, acquired assets or outsourcing in each period in the same manner.


\subsection{Additional Constraints for Near-optimal Solutions}
\label{subsec-add_const_near_optimal}

In the preliminary analysis, some of the constraints we generated through the procedure that is discussed in the subsection of \textit{Valid Inequalities} (\S~\ref{subsec-valid_inequal}) are not always valid, however, these constraints can be used to produce near-optimal solutions in comparably shorter computational times. We also propose these time-saving constraints in addition to valid inequalities in this paper. Additional time-saving constraints for near-optimal solutions are as follows:
\begin{align}
& \sum_{v \in V} \delta_{v} \geq \Theta \label{SND-VI12-MaxNbAssetI} \\
& \sum_{v \in V} \delta_{v} \leq \Theta \label{SND-VI13-MaxNbAssetII} \\
& \sum_{v \in V} \delta_{v} + \sum_{(i,j) \in A_{o}} \sum_{k \in K \cup L} s^{k}_{ij} \geq \Theta \label{SND-VI30-MaxNbAssetIII}
\end{align}

Constraints \eqref{SND-VI12-MaxNbAssetI} and \eqref{SND-VI13-MaxNbAssetII} are proposed to obtain a tight lower bound for the total number of assets utilized in the planning horizon based on maximum predicted number of assets. Even though these constraints are valid in more than half of the valid-cut test instances separately, they do not work for few instances. However, we observe that these constraints are able to decrease the CPU time considerably while diverging reasonably from the optimal solution. 

In the third set of additional inequalities \eqref{SND-VI30-MaxNbAssetIII}, we consider two binary variables representing the total number of resources utilized either as own fleet (owned and leased) or as outsourced service and, enforce them to be greater than or equal to maximum predicted number of assets over all time periods. Likewise, this inequality provides a tight lower bound for number of resources and is more successful than the former ones in hitting optimal solutions in average.

\section{Computational Experiments}
\label{sec-computational_exp}

The performance of the proposed arc-based model is tested on randomly generated problem instances for which we dedicate the following subsection (\S~\ref{subsec-instance-gen}) and describe the details of the generation process as well as values of problem parameters. The formulation is coded on IBM ILOG CPLEX Optimization Studio and solved with CPLEX 12.6.1. All computational experiment is run on a 64-bit PC with Intel Core i7 CPU 3.40 GHz and 16GB RAM. 

In the following subsections, we first describe the instance generation and selected values for the cost parameters in \S~\ref{subsec-instance-gen}. Then, we present the optimal solutions of 60 instances as well as our comments on obtained results in \S~\ref{subsec-solRegModel}. Moreover, we implement a particular sensitivity analysis on some selected small-size instances to point out the importance of cost parameters on optimal solutions. We continue testing the problem instances with some operational restrictions that might be observed in practice in terms of demand shifting. We report the outcomes of these experiments in \S~\ref{subsec-impact-costparam} and \S~\ref{subsec-impact-operRest}, respectively. Finally, we report the results of the experiment on the impact of valid inequalities by testing medium-size instances in \S~\ref{subsec-impact-VI}.

\subsection{Instance Generation}
\label{subsec-instance-gen}

The problem instances for considered CBSND problem are described through the size of the physical network ($N^{'}$), number of periods ($T$) in the planning horizon, number of commodities ($K$), number of owned ($V_{1}$) as well as acquirable assets ($V_{2}$), and number of arcs ($A$) in the data. Thus, the generation process begins with determining the size of the physical (or static) network and distances between all pair of static nodes. Besides, selected values for all sets and parameters of the generated instances that are tested in the computational experiments are summarized in Table~\ref{tab:parameters}.

\begin{table}[htbp]
\centering
\tiny
\caption {Values of the parameters on generated problem instances.} \label{tab:parameters}
{
\begin{tabular}{l c c c r r r r}
\midrule
\textbf{Description} &	& \textbf{Notation} &	& \multicolumn{4}{c}{\bf Value} \\
\midrule
 & & & & \textbf{Small} & \textbf{Medium} & \textbf{Large} & \textbf{Very Large} \\
 \midrule
\textbf{Sets:} & & & & & & &\\
\hspace{1 cm} Number of static nodes   & & $|N'|$  	&	&   5 & 6 & 7 & 10+ \\
\hspace{1 cm} Number of time periods & & $|T|$  &  & \multicolumn{4}{c}{7}    \\
\hspace{1 cm} Number of original commodities & & $|K|$  	&   &  10, 15, 20 &  20, 25, 30  &  30, 36, 42 & 72+,81+,90+   \\
\hspace{1 cm} Number of owned assets & & $|V_1|$  	&   &  7 & 12 & 15 & 35+    \\
\hspace{1 cm} Number of acquirable assets & & $|V_2|$  	&   &  5 & 7 & 10 & 15+   \\
\hspace{1 cm} Number of arcs & & $|A|$  &  & \multicolumn{4}{c}{$ |N'|*|N'-1|*|T|+|A_h|+|A_o| $}    \\
\textbf{Distance Parameter:} & & & & &\\
\hspace{1 cm} Duration of service as distance & & $d_{ij}$  &  & \multicolumn{4}{r}{$ d_{ij} \leq \left\lfloor \frac{|T|}{2}  \right\rfloor \quad \forall i,j \in N^{'} $}    \\
\textbf{Capacity and Flow Parameter:} & & & & &\\
\hspace{1 cm} Service arc capacity & & $u_{ij}$  &  & \multicolumn{4}{r}{$ 1  \quad \forall (i,j) \in A_s $} \\
\hspace{1 cm} Holding arc capacity & & $u_{ij}$  &  & \multicolumn{4}{r}{$ \infty \quad \forall (i,j) \in A_h $} \\
\hspace{1 cm} Volume of flow & & $w^k$  &  & \multicolumn{4}{r}{$ 1  \quad \forall k \in K$}    \\
\textbf{Cost Parameters:} & & & & &\\
\hspace{1 cm} Fixed cost of owned assets  &	& $ f $ & & \multicolumn{4}{r}{25}\\
\hspace{1 cm} Fixed cost of acquired assets   &	&   $ g $ &	& \multicolumn{4}{r}{50} \\
\hspace{1 cm} Routing cost on offered services arcs  &	&  $ c^{k}_{ij} $ &	& \multicolumn{4}{r}{$ \sim U[0.6, 1] \quad \forall (i,j) \in A_s  $}     \\
\hspace{1 cm} Routing cost on arcs of outsourced services &	&  $ c^{k}_{ij} $ &	& \multicolumn{4}{r}{$ 25 + \sim U[1.2, 2] \quad \forall (i,j) \in A_o $}     \\
\hspace{1 cm} Routing cost on offered services arcs  &	&  $ c^{k}_{ij} $ &	& \multicolumn{4}{r}{$ 0.15 \quad \forall (i,j) \in A_h $}     \\
\hspace{1 cm} Penalty multiplier for earliness/tardiness   &	&   $ r_{e}=r_{l} $ &	& \multicolumn{4}{r}{1.2} \\
\hline
\end{tabular}
}
\end{table}

Generated instances are classified as the size of small, medium, large and very large based on the number of nodes in $N^{'}$. The instances with five static nodes are categorized as '\textit{small}', the ones with six nodes as '\textit{medium}', and the instances with seven nodes as '\textit{large}' and the instances having more than ten nodes are classified as '\textit{very large}' in testing. The planning horizon is assumed to be week-long; thus it consists of seven identical periods, each corresponds a day in the week. 
As stated earlier, this is a sample period within one of the peak seasons throughout a year. The demand is defined as a number of commodities that needs to be transported in the planning horizon concerning given release dates and deadlines. We tested three different cases in which the number of commodities is increased gradually 
to simulate a peak season demand more effectively. The heaviest demand case in each size of the problem is limited with the assumption of only one commodity between any O-D pair. For instance, there would be 20 commodities ($ |N'|*|N'-1|$) in small-size problems in the maximum demand case. The amount of flow for each commodity is determined as unit flow without loss of generality. To be consistent with this assumption, the capacity of all service arcs is also selected as one unit in all instances. However, to remind that we assume no capacity on holding arcs to be consistent with the literature. The capacity of outsourced services is imposed by the other service provider and able to handle all volume on that route.    

The overall capacity of the carrier is defined through how many owned assets are available for the upcoming planning horizon. There is also an option of acquiring (or leasing) assets to compensate for the capacity shortage in this problem. The number of owned and acquirable assets are determined based on problem size and indicated in Table~\ref{tab:parameters}. The number of arcs depends on the number of nodes in the physical network, the number of periods in the planning horizon, and the number of outsourced services provided by the market.  

The operating cost parameters of the assets are determined in favor of owned assets because the primary objective of the carrier is to maximize the utilization of owned resources. Any other settings for cost parameters would conflict with "full-asset-utilization" policy that is one of the bases of this paper. The cost parameters are selected such that fixed and acquiring (leasing) cost are assumed as $ 25 $ and $ 50 $, respectively, while routing one unit of flow incurs a random transportation cost of $ UNIF(0.60, 1) $ and $ 25+UNIF(1.2, 2) $ for offered and outsourced services for each arc and commodity combination, respectively. The holding cost is assumed constant and chosen as 0.15 for all holding arcs while the penalty of delivering commodities early and tardy are assumed equal and are determined as a constant ratio of transportation cost as $ 1.2 $.

As stated earlier, we assume that a single type of asset is used for transporting commodities, so the speed of all services is the same. For this reason, the duration of each service between any pair of nodes on the physical network is interpreted as the distance between those static nodes. In addition to this, the duration of any service cannot be longer than half-length of the planning horizon based on the "full-asset utilization" policy that requires asset balancing. To remind that, the number of assets leaving and entering a node must be balanced, and this also imposes that each particular asset must return to its starting node at the beginning of each schedule. As proposed in Corollary~\ref{MaximumDistance}, the maximum distance value cannot be more than half of the schedule length. Since we consider a week-long planning horizon, the distance matrix might include values of 1, 2, and 3 in all generated instances. 

\begin{figure}[htbp]
\centering
\fbox{\begin{varwidth}{1\textwidth}
\subfigure[Close-range network.]{
\includegraphics[clip,width=0.40\textwidth]{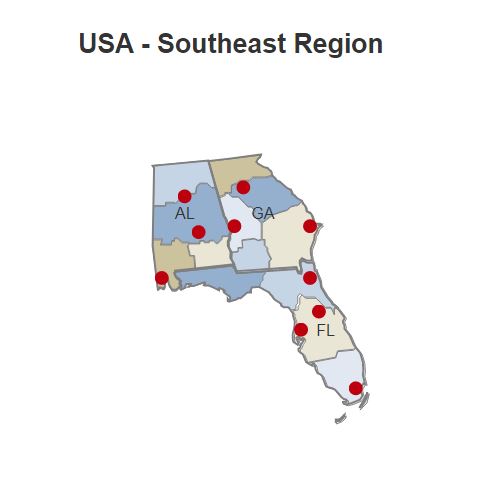}\label{fig:close-range}}
\quad\subfigure[Mid-range network.]{
\includegraphics[clip,width=0.40\textwidth]{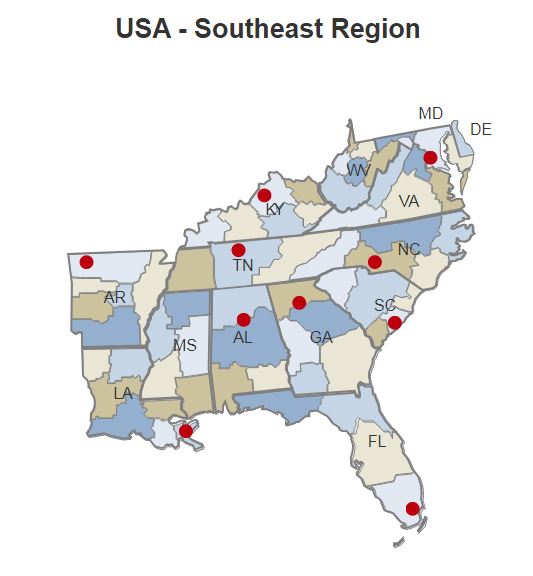}\label{fig:mid-range}}
\end{varwidth}}
\caption{Illustration of close-range and mid-range networks on US regions.}
\label{fig:maps_CR-MR}
\end{figure}

\begin{figure}[htbp]
	\centering
	\fbox{\includegraphics[width=0.7\textwidth]{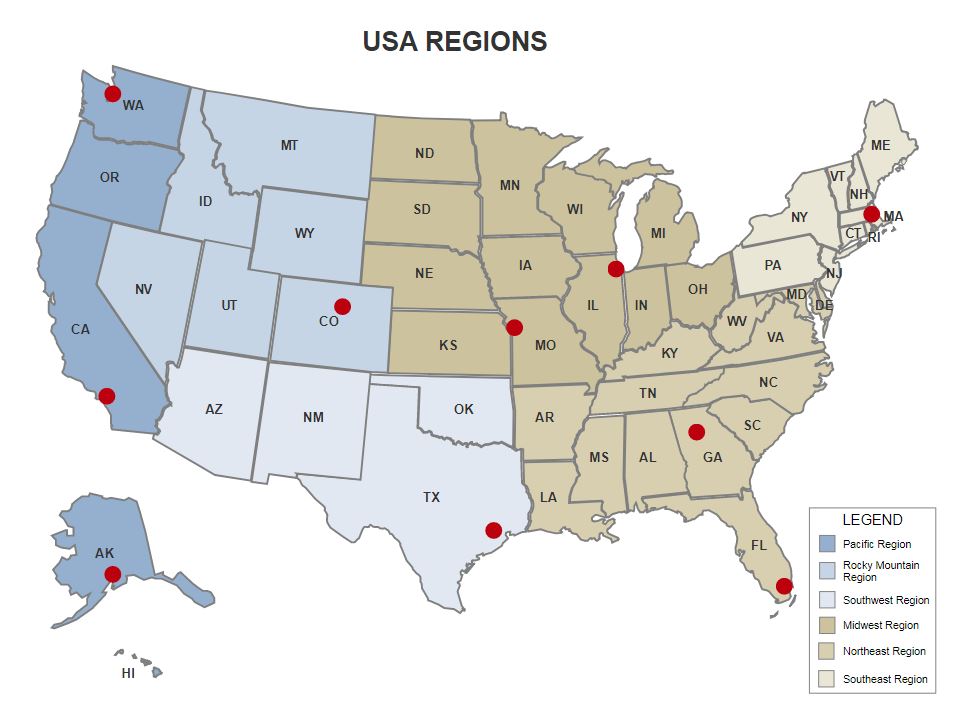}}
	\caption{Illustration of a long-range network on US map.}
	\label{fig:long-range}
\end{figure}

Moreover, we also consider the overall network topology in terms of distance in our analysis. We aim to provide additional insights about how the proximity of nodes on a network affects the MIP formulation both computationally (CPU time) and operationally (total operating cost). Thus, we propose the concept of \textit{distance index} and categorize generated instances based on this value either close-range (CR), medium-range (MR), or long-range (LR) instances. To illustrate this concept, we arbitrarily select ten points which spread over three distinct regions as presented in Figure~\ref{fig:maps_CR-MR} and \ref{fig:long-range}. The ten nodes from three states in the US Southeast Region displayed in Figure\ref{fig:close-range} (Miami (FL), Tampa (FL), Orlando (FL), Atlanta (GA), Columbus (GA), Savannah (GA), Birmingham (AL), Mobile (AL), Montgomery (AL)) compose a close-range network. The network in Figure~\ref{fig:mid-range} corresponds a mid-range network which includes another ten nodes spread over entire US Southeast Region (Miami (FL), Atlanta (GA), New Orleans (LA), Washington (DC), Louisville (KY), Charlotte (NC), Nashville (TN), Birmingham (AL), Rogers (AR), and Charleston (SC)). A long-range network may comprise of ten nodes from all around US (Miami (FL), Seattle (WA), Los Angeles (CA), Chicago (IL), Houston (TX), Boston (MA), Kansas City (MO), Atlanta (GA), Denver (CO), Anchorage (AK)) as depicted in Figure~\ref{fig:long-range}.

The categorization itself is a comparative fact; indeed, a close-range network in a context might be an example of a long-range network in another one, i.e., nodes selected around South Florida might be a close-range network, nodes around Florida might be a mid-range network and nodes displayed in Figure~\ref{fig:close-range} might be a long-range network. However, the comparative feature of this concept does not prevent distance index approach from analyzing in our experiments.  

Without loss of generality, we apply the following rule to categorize the generated instances properly. The range between the minimum possible total distance value (i.e., 20 for small-size instances and 42 for medium-size instances) and the maximum possible total distance value (i.e., 60 for small-size instances and 90 for medium-size instances) is divided into three intervals, each of them approximately takes one-third of the total range. Instances fall within the first interval are assumed as close-range, the ones fall within the second interval are presumed medium-range, and the instances fall within the third interval is assumed as long-range instances. The total distance value and the distance index of each instance belong to the set of small and medium size are presented in a two-fold table below (Table~\ref{tab:instances}) in columns three and four. For the sake of simplicity in referring any instance across our paper, we assign a global ID number to each instance which is also shown in column two in each corresponding part.

\begin{table}[htbp]
	\centering
	\tiny
	\caption {Test instances and corresponding total distance in the static network.} \label{tab:instances}
	{
		\renewcommand{\arraystretch}{0.85}
		\begin{tabular}{l c c c | c | l c c c}
			\midrule
\multirow{2}{*}{\bf Instance}	&	\bf Global	&	\bf Total	&	\bf Distance	&		&	\multirow{2}{*}{\bf Instance}	&	\bf Global	&	\bf Total	&	\bf Distance	\\		
	&	\bf ID	&	\bf Distance	&	\bf Index	&		&		&	\bf ID	&	\bf Distance	&	\bf Index	\\	\cmidrule{1-4}	\cmidrule{6-9}
inst1.n5.c10	&	1	&	\multirow{3}{*}{50}	&	Long	&		&	inst1.n6.c20	&	31	&	\multirow{3}{*}{72}	&	Long	\\		
inst1.n5.c15	&	2	&		&	Range	&		&	inst1.n6.c25	&	32	&		&	Range	\\		
inst1.n5.c20	&	3	&		&	(LR)	&		&	inst1.n6.c30	&	33	&		&	(LR)	\\	\cmidrule{1-4}	\cmidrule{6-9}
inst2.n5.c10	&	4	&	\multirow{3}{*}{48}	&	Long	&		&	inst2.n6.c20	&	34	&	\multirow{3}{*}{52}	&	Medium	\\		
inst2.n5.c15	&	5	&		&	Range	&		&	inst2.n6.c25	&	35	&		&	Range	\\		
inst2.n5.c20	&	6	&		&	(LR)	&		&	inst2.n6.c30	&	36	&		&	(MR)	\\	\cmidrule{1-4}	\cmidrule{6-9}
inst3.n5.c10	&	7	&	\multirow{3}{*}{38}	&	Medium	&		&	inst3.n6.c20	&	37	&	\multirow{3}{*}{62}	&	Medium	\\		
inst3.n5.c15	&	8	&		&	Range	&		&	inst3.n6.c25	&	38	&		&	Range	\\		
inst3.n5.c20	&	9	&		&	(MR)	&		&	inst3.n6.c30	&	39	&		&	(MR)	\\	\cmidrule{1-4}	\cmidrule{6-9}
inst4.n5.c10	&	10	&	\multirow{3}{*}{36}	&	Medium	&		&	inst4.n6.c20	&	40	&	\multirow{3}{*}{68}	&	Medium	\\		
inst4.n5.c15	&	11	&		&	Range	&		&	inst4.n6.c25	&	41	&		&	Range	\\		
inst4.n5.c20	&	12	&		&	(MR)	&		&	inst4.n6.c30	&	42	&		&	(MR)	\\	\cmidrule{1-4}	\cmidrule{6-9}
inst5.n5.c10	&	13	&	\multirow{3}{*}{34}	&	Close	&		&	inst5.n6.c20	&	43	&	\multirow{3}{*}{54}	&	Medium	\\		
inst5.n5.c15	&	14	&		&	Range	&		&	inst5.n6.c25	&	44	&		&	Range	\\		
inst5.n5.c20	&	15	&		&	(CR)	&		&	inst5.n6.c30	&	45	&		&	(MR)	\\	\cmidrule{1-4}	\cmidrule{6-9}
inst6.n5.c10	&	16	&	\multirow{3}{*}{56}	&	Long	&		&	inst6.n6.c20	&	46	&	\multirow{3}{*}{50}	&	Close	\\		
inst6.n5.c15	&	17	&		&	Range	&		&	inst6.n6.c25	&	47	&		&	Range	\\		
inst6.n5.c20	&	18	&		&	(LR)	&		&	inst6.n6.c30	&	48	&		&	(CR)	\\	\cmidrule{1-4}	\cmidrule{6-9}
inst7.n5.c10	&	19	&	\multirow{3}{*}{40}	&	Medium	&		&	inst7.n6.c20	&	49	&	\multirow{3}{*}{66}	&	Medium	\\		
inst7.n5.c15	&	20	&		&	Range	&		&	inst7.n6.c25	&	50	&		&	Range	\\		
inst7.n5.c20	&	21	&		&	(MR)	&		&	inst7.n6.c30	&	51	&		&	(MR)	\\	\cmidrule{1-4}	\cmidrule{6-9}
inst8.n5.c10	&	22	&	\multirow{3}{*}{26}	&	Close	&		&	inst8.n6.c20	&	52	&	\multirow{3}{*}{78}	&	Long	\\		
inst8.n5.c15	&	23	&		&	Range	&		&	inst8.n6.c25	&	53	&		&	Range	\\		
inst8.n5.c20	&	24	&		&	(CR)	&		&	inst8.n6.c30	&	54	&		&	(LR)	\\	\cmidrule{1-4}	\cmidrule{6-9}
inst9.n5.c10	&	25	&	\multirow{3}{*}{46}	&	Medium	&		&	inst9.n6.c20	&	55	&	\multirow{3}{*}{46}	&	Close	\\		
inst9.n5.c15	&	26	&		&	Range	&		&	inst9.n6.c25	&	56	&		&	Range	\\		
inst9.n5.c20	&	27	&		&	(MR)	&		&	inst9.n6.c30	&	57	&		&	(CR)	\\	\cmidrule{1-4}	\cmidrule{6-9}
inst10.n5.c10	&	28	&	\multirow{3}{*}{42}	&	Medium	&		&	inst10.n6.c20	&	58	&	\multirow{3}{*}{74}	&	Long	\\		
inst10.n5.c15	&	29	&		&	Range	&		&	inst10.n6.c25	&	59	&		&	Range	\\		
inst10.n5.c20	&	30	&		&	(MR)	&		&	inst10.n6.c30	&	60	&		&	(LR)	\\

			\hline
		\end{tabular}
	}
\end{table}

Last, we determine a time limit of one hour for small-size instances and four hours for medium-size or above in regular model runs as well as runs with valid inequalities in all experiments.

\subsection{Solutions of the Regular Model}
\label{subsec-solRegModel}

In the first part of the computational experiment, we tested small (instances with 5-nodes physically) and medium (instances with 6-nodes physically) size instances with the regular model and presented the optimal solutions in Table~\ref{tab:n5_results} and Table~\ref{tab:n6_results}, respectively. The first five columns provide information about instances' characteristics including instance name, global instance ID across the paper, number of nodes in the static network, number of original commodities, and number of owned and acquirable (leasable) assets. The optimal solution for each instance is given starting from column-7 through column-12. The number of owned assets utilized and the number of assets leased in each corresponding instance are shown in columns 7 and 8, respectively. In columns 9 through 11, we present the number of commodities delivered on time, that of demand shifted to an earlier/a later date, and that of commodities delivered via an outsourced service on time, respectively. Finally, we present total operating cost obtained for each problem in column 12.

\begin{table}[htbp]
	\centering
	\scriptsize
	\caption {Optimal solutions for instances with $|N|=5$.} \label{tab:n5_results}
	{
		\begin{tabular}{l c c c c c c c c c c c c c}
			\midrule
		
\multirow{3}{*}{\bf Instance}	&	\multirow{3}{*}{\bf ID}	&	\multicolumn{4}{c}{\bf Problem Size}							&		&	\multicolumn{2}{c}{\bf \# of Assets}			&		&	\multicolumn{3}{c}{\bf \# of Commodities}					&	\multirow{1.5}{*}{\bf Total}	\\	\cmidrule{3-6}	\cmidrule{8-9}	\cmidrule{11-13}
	&		&	\multirow{2}{*}{\bf $|N'| $}	&	\multirow{2}{*}{\bf $|K| $}	&	\multirow{2}{*}{\bf $|V_1| $}	&	\multirow{2}{*}{\bf $|V_2| $}	&		&	\multirow{2}{*}{\bf Owned}	&	\multirow{2}{*}{\bf Leased}	&		&	\bf On	&	\bf Early/	&	\multirow{2}{*}{\bf Outsourced}	&	\multirow{1.5}{*}{\bf Cost}	\\			
	&		&		&		&		&		&		&		&		&		&	\bf Time	&	\bf Tardy	&		&		\\	\midrule		
inst1.n5.c10	&	1	&	5	&	10	&	7	&	5	&		&	6	&	-	&		&	3	&	2/5	&	-	&	162.37	\\			
inst1.n5.c15	&	2	&	5	&	15	&	7	&	5	&		&	7	&	-	&		&	6	&	6/2	&	1	&	215.51	\\			
inst1.n5.c20	&	3	&	5	&	20	&	7	&	5	&		&	7	&	1	&		&	8	&	4/6	&	2	&	297.86	\\	\midrule		
inst2.n5.c10	&	4	&	5	&	10	&	7	&	5	&		&	6	&	-	&		&	4	&	3/3	&	-	&	161.42	\\			
inst2.n5.c15	&	5	&	5	&	15	&	7	&	5	&		&	7	&	-	&		&	1	&	7/7	&	-	&	191.48	\\			
inst2.n5.c20	&	6	&	5	&	20	&	7	&	5	&		&	7	&	-	&		&	3	&	7/7	&	3	&	270.76	\\	\midrule		
inst3.n5.c10	&	7	&	5	&	10	&	7	&	5	&		&	5	&	-	&		&	2	&	4/4	&	-	&	135.64	\\			
inst3.n5.c15	&	8	&	5	&	15	&	7	&	5	&		&	7	&	-	&		&	7	&	4/4	&	-	&	189.82	\\			
inst3.n5.c20	&	9	&	5	&	20	&	7	&	5	&		&	7	&	-	&		&	8	&	5/6	&	1	&	221.84	\\	\midrule		
inst4.n5.c10	&	10	&	5	&	10	&	7	&	5	&		&	4	&	-	&		&	3	&	2/5	&	-	&	109.9	\\			
inst4.n5.c15	&	11	&	5	&	15	&	7	&	5	&		&	6	&	-	&		&	9	&	3/3	&	-	&	164.28	\\			
inst4.n5.c20	&	12	&	5	&	20	&	7	&	5	&		&	7	&	-	&		&	11	&	7/2	&	-	&	195.2	\\	\midrule		
inst5.n5.c10	&	13	&	5	&	10	&	7	&	5	&		&	5	&	-	&		&	5	&	3/2	&	-	&	135.52	\\			
inst5.n5.c15	&	14	&	5	&	15	&	7	&	5	&		&	6	&	-	&		&	7	&	6/2	&	-	&	164.89	\\			
inst5.n5.c20	&	15	&	5	&	20	&	7	&	5	&		&	7	&	-	&		&	12	&	3/5	&	-	&	194.23	\\	\midrule		
inst6.n5.c10	&	16	&	5	&	10	&	7	&	5	&		&	7	&	-	&		&	6	&	1/3	&	-	&	184.7	\\			
inst6.n5.c15	&	17	&	5	&	15	&	7	&	5	&		&	7	&	-	&		&	5	&	5/3	&	2	&	240.81	\\			
inst6.n5.c20	&	18	&	5	&	20	&	7	&	5	&		&	7	&	2	&		&	4	&	7/7	&	2	&	348.63	\\	\midrule		
inst7.n5.c10	&	19	&	5	&	10	&	7	&	5	&		&	5	&	-	&		&	5	&	2/3	&	-	&	134.3	\\			
inst7.n5.c15	&	20	&	5	&	15	&	7	&	5	&		&	6	&	-	&		&	4	&	6/5	&	-	&	167.23	\\			
inst7.n5.c20	&	21	&	5	&	20	&	7	&	5	&		&	7	&	1	&		&	8	&	6/6	&	-	&	243.7	\\	\midrule		
inst8.n5.c10	&	22	&	5	&	10	&	7	&	5	&		&	4	&	-	&		&	4	&	3/3	&	-	&	109.26	\\			
inst8.n5.c15	&	23	&	5	&	15	&	7	&	5	&		&	5	&	-	&		&	9	&	4/2	&	-	&	139.03	\\			
inst8.n5.c20	&	24	&	5	&	20	&	7	&	5	&		&	6	&	-	&		&	9	&	8/3	&	-	&	169.64	\\	\midrule		
inst9.n5.c10	&	25	&	5	&	10	&	7	&	5	&		&	7	&	-	&		&	3	&	4/3	&	-	&	184.42	\\			
inst9.n5.c15	&	26	&	5	&	15	&	7	&	5	&		&	7	&	-	&		&	4	&	5/4	&	2	&	242.25	\\			
inst9.n5.c20	&	27	&	5	&	20	&	7	&	5	&		&	7	&	1	&		&	7	&	6/5	&	2	&	294.88	\\	\midrule		
inst10.n5.c10	&	28	&	5	&	10	&	7	&	5	&		&	5	&	-	&		&	6	&	4/0	&	-	&	134.61	\\			
inst10.n5.c15	&	29	&	5	&	15	&	7	&	5	&		&	6	&	-	&		&	8	&	5/2	&	-	&	165.62	\\			
inst10.n5.c20	&	30	&	5	&	20	&	7	&	5	&		&	7	&	1	&		&	10	&	2/8	&	-	&	244.76	\\			
			\hline
		\end{tabular}
	}
\end{table}

In the first place, the number of assets utilized increases as the demand increases in both small and medium size problems, that is an expected outcome, obviously. Although increasing trend can be observed in each group of the instance, the number of assets utilized varies across the group of instances due to differences in topology on the static networks. For example, all of the owned assets are utilized in 16 out of 30 small-size instances in which the capacity shortage is overcame either by outsourcing (in instances 2, 6, 9, 17, and 26), or by leasing additional asset (in instances 21 and 30), or by both (in instances 3, 18, and 27), or only by demand shifting (in instances 5, 8, 12, 15, 16, 25). All these instances but one (instance-15) are classified either as long-range or medium-range cases based on physical networks. In one-third of medium-size instances all owned assets are completely utilized (instances 32, 33, 42, 51, 52, 53, 54, 58, 59, and 60), outsourcing (instances 53 and 56) and leasing additional asset(s) (instances 33, 53, and 54) are chosen as an option for handling excess demand. Again all of these medium-size instances belong to either long-range or medium-range categories in terms of physical network distance index.     

To compare the asset need between long-range and close-range physical networks, the total number of assets utilized is as low as four assets (instance-22) in close-range cases while this value goes up to eight and nine assets in instances 3 and 18 in the long-range counterparts. The maximum total number of assets utilized is doubled, and the minimum number is not lower than six assets across all long-range small-size instances. Thus, the results confirm that asset utilization is directly affected by the distances between nodes on the physical network. Relatively less number of assets are utilized when the static nodes are getting closer, a substantial amount of resource is required to meet the transportation demand when the distances between nodes are getting higher.  

To examine the optimal results in terms of response actions, the model chooses to add an asset(s) into the fleet in five out of 30 small-size and three out of 30 medium-size instances. Besides, the outsourcing option is selected in eight out of 30 small-size and three out of 30 medium-size instances for 15 and five commodities in total, respectively. The outsourcing option is mostly selected to deliver commodities in long-range problems (ten commodities); however, a few commodities are outsourced in mid-range cases (five commodities) in the set of small-size instances, as well. This outcome is consistent in medium-size instances, four outsourced commodities are included in demand of long-range networks.

\begin{table}[htbp]
	\centering
	\scriptsize
	\caption {Optimal solutions for instances with $|N|=6$.} \label{tab:n6_results}
	{
		\begin{tabular}{l c c c c c c c c c c c c r}
			\midrule
\multirow{3}{*}{\bf Instance}	&	\multirow{3}{*}{\bf ID}	&	\multicolumn{4}{c}{\bf Problem Size}							&		&	\multicolumn{2}{c}{\bf \# of Assets}			&		&	\multicolumn{3}{c}{\bf \# of Commodities}					&	\multirow{1.5}{*}{\bf Total}	\\	\cmidrule{3-6}	\cmidrule{8-9}	\cmidrule{11-13}
	&		&	\multirow{2}{*}{\bf $|N'| $}	&	\multirow{2}{*}{\bf $|K| $}	&	\multirow{2}{*}{\bf $|V_1| $}	&	\multirow{2}{*}{\bf $|V_2| $}	&		&	\multirow{2}{*}{\bf Owned}	&	\multirow{2}{*}{\bf Leased}	&		&	\bf On	&	\bf Early/	&	\multirow{2}{*}{\bf Outsourced}	&	\multirow{1.5}{*}{\bf Cost}	\\			
	&		&		&		&		&		&		&		&		&		&	\bf Time	&	\bf Tardy	&		&		\\	\midrule		
inst1.n6.c20	&	31	&	6	&	20	&	12	&	7	&		&	10	&	-	&		&	7	&	6/7	&	-	&	270.54	\\			
inst1.n6.c25	&	32	&	6	&	25	&	12	&	7	&		&	12	&	-	&		&	15	&	5/5	&	-	&	325.624	\\			
inst1.n6.c30	&	33	&	6	&	30	&	12	&	7	&		&	12	&	1	&		&	15	&	5/10	&	-	&	380.984	\\	\midrule		
inst2.n6.c20	&	34	&	6	&	20	&	12	&	7	&		&	7	&	-	&		&	10	&	4/6	&	-	&	200.816	\\			
inst2.n6.c25	&	35	&	6	&	25	&	12	&	7	&		&	9	&	-	&		&	11	&	6/8	&	-	&	252.796	\\			
inst2.n6.c30	&	36	&	6	&	30	&	12	&	7	&		&	10	&	-	&		&	14	&	10/6	&	-	&	283.342	\\	\midrule		
inst3.n6.c20	&	37	&	6	&	20	&	12	&	7	&		&	8	&	-	&		&	7	&	8/5	&	-	&	221.202	\\			
inst3.n6.c25	&	38	&	6	&	25	&	12	&	7	&		&	10	&	-	&		&	13	&	4/8	&	-	&	273.25	\\			
inst3.n6.c30	&	39	&	6	&	30	&	12	&	7	&		&	10	&	-	&		&	12	&	9/9	&	-	&	285.23	\\	\midrule		
inst4.n6.c20	&	40	&	6	&	20	&	12	&	7	&		&	9	&	-	&		&	8	&	8/3	&	1	&	271.176	\\			
inst4.n6.c25	&	41	&	6	&	25	&	12	&	7	&		&	11	&	-	&		&	13	&	6/6	&	-	&	303.052	\\			
inst4.n6.c30	&	42	&	6	&	30	&	12	&	7	&		&	12	&	-	&		&	12	&	6/12	&	-	&	334.638	\\	\midrule		
inst5.n6.c20	&	43	&	6	&	20	&	12	&	7	&		&	8	&	-	&		&	8	&	7/5	&	-	&	221.326	\\			
inst5.n6.c25	&	44	&	6	&	25	&	12	&	7	&		&	9	&	-	&		&	12	&	4/9	&	-	&	250.714	\\			
inst5.n6.c30	&	45	&	6	&	30	&	12	&	7	&		&	10	&	-	&		&	14	&	6/10	&	-	&	283.134	\\	\midrule		
inst6.n6.c20	&	46	&	6	&	20	&	12	&	7	&		&	8	&	-	&		&	9	&	6/5	&	-	&	220.764	\\			
inst6.n6.c25	&	47	&	6	&	25	&	12	&	7	&		&	9	&	-	&		&	11	&	6/8	&	-	&	252.612	\\			
inst6.n6.c30	&	48	&	6	&	30	&	12	&	7	&		&	9	&	-	&		&	10	&	9/11	&	-	&	258.958	\\	\midrule		
inst7.n6.c20	&	49	&	6	&	20	&	12	&	7	&		&	10	&	-	&		&	10	&	4/6	&	-	&	270.074	\\			
inst7.n6.c25	&	50	&	6	&	25	&	12	&	7	&		&	10	&	-	&		&	7	&	11/7	&	-	&	279.326	\\			
inst7.n6.c30	&	51	&	6	&	30	&	12	&	7	&		&	12	&	-	&		&	12	&	10/8	&	-	&	331.228	\\	\midrule		
inst8.n6.c20	&	52	&	6	&	20	&	12	&	7	&		&	12	&	-	&		&	7	&	10/3	&	-	&	320.778	\\			
inst8.n6.c25	&	53	&	6	&	25	&	12	&	7	&		&	12	&	1	&		&	12	&	6/7	&	1	&	402.98	\\			
inst8.n6.c30	&	54	&	6	&	30	&	12	&	7	&		&	12	&	3	&		&	14	&	8/8	&	-	&	482.826	\\	\midrule		
inst9.n6.c20	&	55	&	6	&	20	&	12	&	7	&		&	7	&	-	&		&	7	&	9/4	&	-	&	193.49	\\			
inst9.n6.c25	&	56	&	6	&	25	&	12	&	7	&		&	8	&	-	&		&	9	&	8/8	&	-	&	233.668	\\			
inst9.n6.c30	&	57	&	6	&	30	&	12	&	7	&		&	9	&	-	&		&	12	&	10/8	&	-	&	254.224	\\	\midrule		
inst10.n6.c20	&	58	&	6	&	20	&	12	&	7	&		&	11	&	-	&		&	11	&	4/5	&	-	&	294.856	\\			
inst10.n6.c25	&	59	&	6	&	25	&	12	&	7	&		&	12	&	-	&		&	16	&	6/3	&	-	&	327.904	\\			
inst10.n6.c30	&	60	&	6	&	30	&	12	&	7	&		&	12	&	-	&		&	9	&	9/9	&	3	&	409.11	\\			
			\hline
		\end{tabular}
	}
\end{table}

Moreover, the option of demand shifting seems very attractive based on current cost values, some portion of demand is shifted either at an earlier or a late period in all of the small and medium size instances. The total number of commodities shifted is greater than or equal to that of delivered on time in 23 out of 30 small size instances, whereas this statistics is increased to 25 among medium-size instances.

The similar trends are observed in total operating cost values of the test problems. It increases as the number of commodities increases as well in each instance group. The second observation is that total operating cost is strongly affected by the topology of the physical network, it is more likely to increase when the distances are getting longer since longer distances mean more resources to utilize in transportation activities.   

The reason why the mathematical model mostly chooses the outsourcing option rather than adding a leased resource into the fleet is that the balance between values of the cost parameters of options in our opinion. Thus, we decided to conduct a particular test with a portion of small-size instances in \S~\ref{subsec-impact-costparam} and \S~\ref{subsec-impact-operRest} in the next two subsections to be able to illustrate that how cost parameters are effective on the choice of options to meet the excess demand. We intentionally left the discussion of computational time in \S\ref{subsec-impact-VI} where we compare the performance of regular formulation against valid inequalities.

\subsection{Impact of Cost Parameters on Solutions}
\label{subsec-impact-costparam}

As mentioned earlier, it is observed that cost parameters of the CSSND problem are very effective on optimal solutions in the way of choosing the option(s) to handle capacity shortage (or excess demand) issue in our computational experiment. 

To emphasize the importance of selecting cost parameters properly as well as to prove the flexibility of the mathematical model, we conduct a particular cost analysis on some selected small-size instances in which owned assets are fully utilized. 
This experiment is designed as a \textit{ceteris paribus experiment} and consists of two parts: (1) change in leasing (acquiring) cost and (2) change in penalty ratio of earliness/tardiness. The selected instances with global IDs as well as optimal solutions for each instance that belong to corresponding part of the experiment are shown in
Table~\ref{tab:costAnalysis_g} and \ref{tab:costAnalysis_r2}, respectively. 

First, we tested three values of $g$ around of the fixed cost of utilizing owned assets and the cost of outsourcing option such as $25, 26 $, and $27$, and compared the results with the default case. 
The optimal solutions are presented in 
Table~\ref{tab:costAnalysis_g} 
and it is structured in 2x2 wide rows and columns in which optimal decisions of the problem are summarized with respect to chosen value of $g$ given on top of that column including the default settings. 
The first two columns identify the instance name and ID, the next three columns (columns 3 through 5) show number of assets utilized in the corresponding instance as owned, leased and sum of them in order. Columns 6 through 9 display the number of commodities delivered on time, that of demand delivered either as early/late and via an outsourced service. Finally, the total operating cost of each instance is presented in the last column. This structure is repeated in the second wide column and row, as well. 

When the leasing cost is as low as utilizing an owned asset ($f=g=25$), the model acts in favor of leasing additional assets among the three response actions. Additional assets are leased in all instances but instances 21 and 30, the number of utilized assets increases by 17\% in average across ten instances compared to default case. The total number of commodities delivered by an outsourced service decreases from 15 to zero, while the number of on time deliveries increases by 29\% (from 63 to 81). The total number of shifted commodities has a small change and decreases by 2.8\% (from 107 to 104). 

\begin{landscape} 
 \begin{table}[htbp]
	\centering
	\scriptsize
	\caption {Sensitivity analysis of selected instances with respect to changes in cost, $g$.} \label{tab:costAnalysis_g} 
	{
		\renewcommand{\arraystretch}{1}
		\begin{tabular}{l c c c c c c c c r c c c c c c c c r}
			\midrule
\multirow{4}{*}{\bf Instance}	&	\multirow{4}{*}{\bf ID}	&			\multicolumn{8}{c}{\bf $g=50 $ (original setting)}																	&		&	\multicolumn{8}{c}{\bf $g=25 \: (50\% \: decrease, f=g)$}																			\\	\cmidrule{3-10}	\cmidrule{12-19}				
	&		&			\multicolumn{3}{c}{\bf \# of Assets}					&		&	\multicolumn{3}{c}{\bf \# of Commodities}					&			\multirow{3}{*}{\minitab[c]{\bf Total \\ \bf Cost}}	&		&	\multicolumn{3}{c}{\bf \# of Assets}					&		&	\multicolumn{3}{c}{\bf \# of Commodities}					&			\multirow{3}{*}{\minitab[c]{\bf Total \\ \bf Cost}}			\\	\cmidrule{3-5}	\cmidrule{7-9}	\cmidrule{12-14}	\cmidrule{16-18}		
	&		&			\multirow{2}{*}{\bf Owned}	&	\multirow{2}{*}{\bf Leased}	&	\multirow{2}{*}{\bf Total}	&		&	\bf On	&	\bf Early/	&	\multirow{2}{*}{\bf Out.}	&				&		&	\multirow{2}{*}{\bf Owned}	&	\multirow{2}{*}{\bf Leased}	&	\multirow{2}{*}{\bf Total}	&		&	\bf On	&	\bf Early/	&	\multirow{2}{*}{\bf Out.}	&						\\						
	&		&				&		&		&		&	\bf Time	&	\bf Tardy	&		&				&		&		&		&		&		&	\bf Time	&	\bf Tardy	&		&						\\	\midrule					
inst1.n5.c15	&	2	&			7	&	-	&	7	&		&	6	&	6/2	&	1	&			215.51	&		&	7	&	1	&	8	&		&	7	&	4/4	&	-	&			214.582			\\						
inst1.n5.c20	&	3	&			7	&	1	&	8	&		&	8	&	4/6	&	2	&			297.86	&		&	7	&	3	&	10	&		&	9	&	5/6	&	-	&			270.566			\\						
inst2.n5.c20	&	6	&			7	&	-	&	7	&		&	3	&	7/7	&	3	&			270.76	&		&	7	&	2	&	9	&		&	7	&	7/6	&	-	&			244.256			\\						
inst3.n5.c20	&	9	&			7	&	-	&	7	&		&	8	&	5/6	&	1	&			221.84	&		&	7	&	1	&	8	&		&	6	&	8/6	&	-	&			220.354			\\						
inst6.n5.c15	&	17	&			7	&	-	&	7	&		&	5	&	5/3	&	2	&			240.81	&		&	7	&	2	&	9	&		&	6	&	5/4	&	-	&			240.266			\\						
inst6.n5.c20	&	18	&			7	&	2	&	9	&		&	4	&	7/7	&	2	&			348.63	&		&	7	&	4	&	11	&		&	11	&	4/5	&	-	&			294.704			\\						
inst7.n5.c20	&	21	&			7	&	1	&	8	&		&	8	&	6/6	&	-	&			243.7	&		&	7	&	1	&	8	&		&	8	&	6/6	&	-	&			218.698			\\						
inst9.n5.c15	&	26	&			7	&	-	&	7	&		&	4	&	5/4	&	2	&			242.25	&		&	7	&	2	&	9	&		&	8	&	4/3	&	-	&			238.85			\\						
inst9.n5.c20	&	27	&			7	&	1	&	8	&		&	7	&	6/5	&	2	&			294.88	&		&	7	&	2	&	9	&		&	9	&	5/6	&	-	&			245.416			\\						
inst10.n5.c20	&	30	&			7	&	1	&	8	&		&	10	&	2/8	&	-	&			244.76	&		&	7	&	1	&	8	&		&	10	&	2/8	&	-	&			219.76			\\						
																																													\\	\midrule					
\multirow{4}{*}{\bf Instance}	&	\multirow{4}{*}{\bf ID}	&			\multicolumn{8}{c}{\bf $g=26 \: ( f < g, g < c^{k}_{ij} \: \forall \: (i,j) \in A_o)$}																	&		&	\multicolumn{8}{c}{\bf  $g=27  \: ( f \le g)$}																			\\	\cmidrule{3-10}	\cmidrule{12-19}				
	&		&			\multicolumn{3}{c}{\bf \# of Assets}					&		&	\multicolumn{3}{c}{\bf \# of Commodities}					&			\multirow{3}{*}{\minitab[c]{\bf Total \\ \bf Cost}}	&		&	\multicolumn{3}{c}{\bf \# of Assets}					&		&	\multicolumn{3}{c}{\bf \# of Commodities}					&			\multirow{3}{*}{\minitab[c]{\bf Total \\ \bf Cost}}			\\	\cmidrule{3-5}	\cmidrule{7-9}	\cmidrule{12-14}	\cmidrule{16-18}		
	&		&			\multirow{2}{*}{\bf Owned}	&	\multirow{2}{*}{\bf Leased}	&	\multirow{2}{*}{\bf Total}	&		&	\bf On	&	\bf Early/	&	\multirow{2}{*}{\bf Out.}	&				&		&	\multirow{2}{*}{\bf Owned}	&	\multirow{2}{*}{\bf Leased}	&	\multirow{2}{*}{\bf Total}	&		&	\bf On	&	\bf Early/	&	\multirow{2}{*}{\bf Out.}	&						\\						
	&		&				&		&		&		&	\bf Time	&	\bf Tardy	&		&				&		&		&		&		&		&	\bf Time	&	\bf Tardy	&		&						\\	\midrule					
inst1.n5.c15	&	2	&			7	&	-	&	7	&		&	6	&	6/2	&	1	&			215.514	&		&	7	&	-	&	7	&		&	7	&	2/5	&	1	&			215.68			\\						
inst1.n5.c20	&	3	&			7	&	3	&	10	&		&	6	&	6/8	&	-	&			273.696	&		&	7	&	2	&	9	&		&	9	&	5/5	&	1	&			275.284			\\						
inst2.n5.c20	&	6	&			7	&	2	&	9	&		&	7	&	7/6	&	-	&			246.256	&		&	7	&	2	&	9	&		&	7	&	7/6	&	-	&			248.256			\\						
inst3.n5.c20	&	9	&			7	&	1	&	8	&		&	7	&	7/6	&	-	&			222.016	&		&	7	&	1	&	8	&		&	11	&	5/4	&	-	&			222.608			\\						
inst6.n5.c15	&	17	&			7	&	-	&	7	&		&	5	&	5/3	&	2	&			240.81	&		&	7	&	-	&	7	&		&	5	&	5/3	&	2	&			240.81			\\						
inst6.n5.c20	&	18	&			7	&	3	&	10	&		&	7	&	6/6	&	1	&			298.526	&		&	7	&	3	&	10	&		&	7	&	6/6	&	1	&			301.526			\\						
inst7.n5.c20	&	21	&			7	&	1	&	8	&		&	8	&	6/6	&	-	&			219.698	&		&	7	&	1	&	8	&		&	8	&	6/6	&	-	&			220.698			\\						
inst9.n5.c15	&	26	&			7	&	1	&	8	&		&	5	&	5/4	&	1	&			240.704	&		&	7	&	1	&	8	&		&	5	&	5/4	&	1	&			241.704			\\						
inst9.n5.c20	&	27	&			7	&	2	&	9	&		&	9	&	5/6	&	-	&			247.416	&		&	7	&	2	&	9	&		&	9	&	5/6	&	-	&			249.416			\\						
inst10.n5.c20	&	30	&			7	&	1	&	8	&		&	10	&	2/8	&	-	&			220.76	&		&	7	&	1	&	8	&		&	10	&	2/8	&	-	&			221.76			\\						
																																																			
			\hline			
		\end{tabular}
	}
\end{table} 
\end{landscape}

The change in acquiring cost does not affect some of the instances at all. Instances 21 and 30 are insensitive to changes in acquiring cost, since the optimal solutions are all same across the test for them. Instances 6, 17, 26 and 27 are insensitive to changes in $g$ for some specific values. Instances 6, 26 and 27 are unaffected for  $ 26 \leq g < 50 $, instance 17 is indifferent for $ 25 < g $.

In the second part of this experiment, we conducted a test to understand the impact of changes in the penalty ratio on the optimal solutions. We tested four values of $r$ such that $1.5, 2.5, 5 $ and, $10$. The results of the same ten test instances are presented in Table~\ref{tab:costAnalysis_r2} and it is structured in 2x3 wide rows and columns where optimal decisions of the problem are summarized with respect to chosen value of $r$ given on top of that wide-column. The first column identifies the instance ID, the next three columns (columns 2 through 4) show number of assets utilized in the corresponding instance as owned (Own), leased (Lea) and sum of both (Tot) in this order. Columns 5 through 7 display number of commodities delivered on time (OT), that of demand delivered either as early/late (E/T) and commodities delivered via an outsourced service (Out). Finally, total operating cost of each instance is presented in the last column. This structure is repeated in the second and third wide columns and second wide-row, as well.

The consequences of the second part can be summarized in such a general way as penalty ratio ($r$) increases: (1) number of commodities delivered on time by offered and outsourced services increases while number of shifted commodities (delivered early/late) decreases, which is expected obviously, (2) the model utilizes all owned assets to deliver as many commodities as possible on time, however, it chooses to outsource instead of leasing additional asset(s) since delivering by a leased asset (on-time delivery even) is more expensive compared to outsourcing. 

The most dramatic change in choosing a response action to excess demand for the carrier is seen in instance-18 when going from $r=1.2$ to $r=1.5$. The number of assets utilized decreases by two, the number of commodities shifted decreases from 14 to five, the number of commodities delivered on-time increases from four to nine and the number of commodities outsourced increases from two to six. Instance-18, itself, shows the impact of earliness/tardiness penalty parameter on optimal solutions in choosing the most appropriate response action. 

A careful reader may notice the abnormality existing in instance-9 when going from $r=2.5$ to $r=5$. The number of commodities delivered on-time decreases from 12 to 11 while the number of commodities shifted increases from seven to eight. This is basically in conflict with general trend in this test, the reason is that the solver reached the one hour of time limit and branch-and-bound ended with a optimality gap of 3.61\%. There might be another integer solution which is not discovered yet in this specific instance.

\begin{landscape} 
 \begin{table}[htbp]
	\centering
	\tiny
	\caption {Sensitivity analysis of selected instances with respect to changes in penalty ratio, $r$.} \label{tab:costAnalysis_r2}
	{
		\renewcommand{\arraystretch}{1.1}
		\begin{tabular}{l c c c c c c c r c c c c c c c c r c c c c c c c c r}
			\midrule
\multirow{4}{*}{\bf Ins.}	&	\multicolumn{8}{c}{\bf $r=1.2$ (original setting)}															&		&	\multicolumn{8}{c}{\bf $r=1.5$}															&		&	\multicolumn{8}{c}{\bf $r=2.5$}															\\	\cmidrule{2-9}	\cmidrule{11-18}	\cmidrule{20-27}			
	&	\multicolumn{3}{c}{\bf \# of Assets}					&		&	\multicolumn{3}{c}{\bf \# of Commodities}					&	\multirow{3}{*}{\minitab[c]{\bf Total \\ \bf Cost}}	&		&	\multicolumn{3}{c}{\bf \# of Assets}					&		&	\multicolumn{3}{c}{\bf \# of Commodities}					&	\multirow{3}{*}{\minitab[c]{\bf Total \\ \bf Cost}}	&		&	\multicolumn{3}{c}{\bf \# of Assets}					&		&	\multicolumn{3}{c}{\bf \# of Commodities}					&	\multirow{3}{*}{\minitab[c]{\bf Total \\ \bf Cost}}	\\	\cmidrule{2-4}	\cmidrule{6-8}	\cmidrule{11-13}	\cmidrule{15-17}	\cmidrule{20-22}	\cmidrule{24-26}
	&	\multirow{2}{*}{\bf Own}	&	\multirow{2}{*}{\bf Lea}	&	\multirow{2}{*}{\bf Tot}	&		&	\multirow{2}{*}{\bf OT}	&	\multirow{2}{*}{\bf E/T}	&	\multirow{2}{*}{\bf Out}	&		&		&	\multirow{2}{*}{\bf Own}	&	\multirow{2}{*}{\bf Lea}	&	\multirow{2}{*}{\bf Tot}	&		&	\multirow{2}{*}{\bf OT}	&	\multirow{2}{*}{\bf E/T}	&	\multirow{2}{*}{\bf Out}	&		&		&	\multirow{2}{*}{\bf Own}	&	\multirow{2}{*}{\bf Lea}	&	\multirow{2}{*}{\bf Tot}	&		&	\multirow{2}{*}{\bf OT}	&	\multirow{2}{*}{\bf E/T}	&	\multirow{2}{*}{\bf Out}	&		\\						
	&		&		&		&		&		&		&		&		&		&		&		&		&		&		&		&		&		&		&		&		&		&		&		&		&		&		\\	\midrule					
2	&	7	&	-	&	7	&		&	6	&	6/2	&	1	&	215.51	&		&	7	&	-	&	7	&		&	9	&	2/3	&	1	&	217.415	&		&	7	&	-	&	7	&		&	11	&	1/2	&	1	&	221.075	\\						
3	&	7	&	1	&	8	&		&	8	&	4/6	&	2	&	297.86	&		&	7	&	1	&	8	&		&	9	&	3/6	&	2	&	300.665	&		&	7	&	-	&	7	&		&	12	&	2/2	&	4	&	305.565	\\						
6	&	7	&	-	&	7	&		&	3	&	7/7	&	3	&	270.76	&		&	7	&	-	&	7	&		&	3	&	7/7	&	3	&	274.19	&		&	7	&	-	&	7	&		&	3	&	7/7	&	3	&	285.63	\\						
9	&	7	&	-	&	7	&		&	8	&	5/6	&	1	&	221.84	&		&	7	&	-	&	7	&		&	9	&	4/6	&	1	&	224.705	&		&	7	&	-	&	7	&		&	12	&	4/3	&	1	&	233.185	\\						
17	&	7	&	-	&	7	&		&	5	&	5/3	&	2	&	240.81	&		&	7	&	-	&	7	&		&	5	&	5/3	&	2	&	242.85	&		&	7	&	-	&	7	&		&	5	&	5/3	&	2	&	249.65	\\						
18	&	7	&	2	&	9	&		&	4	&	7/7	&	2	&	348.63	&		&	7	&	-	&	7	&		&	9	&	2/3	&	6	&	350.3	&		&	7	&	-	&	7	&		&	9	&	2/3	&	6	&	354.68	\\						
21	&	7	&	1	&	8	&		&	8	&	6/6	&	-	&	243.7	&		&	7	&	-	&	7	&		&	10	&	2/6	&	2	&	246.39	&		&	7	&	-	&	7	&		&	10	&	2/6	&	2	&	252.81	\\						
26	&	7	&	-	&	7	&		&	4	&	5/4	&	2	&	242.25	&		&	7	&	-	&	7	&		&	5	&	4/4	&	2	&	244.785	&		&	7	&	-	&	7	&		&	5	&	4/4	&	2	&	252.675	\\						
27	&	7	&	1	&	8	&		&	7	&	6/5	&	2	&	294.88	&		&	7	&	1	&	8	&		&	8	&	6/4	&	2	&	297.505	&		&	7	&	-	&	7	&		&	11	&	2/3	&	4	&	301.85	\\						
30	&	7	&	1	&	8	&		&	10	&	2/8	&	-	&	244.76	&		&	7	&	-	&	7	&		&	11	&	2/5	&	2	&	246.8	&		&	7	&	-	&	7	&		&	12	&	3/3	&	2	&	252.26	\\						
																																																					\\	\midrule					
\multirow{4}{*}{\bf Ins.}	&	\multicolumn{8}{c}{\bf $r=5$}															&		&	\multicolumn{8}{c}{\bf $r=10$}															&																		\\	\cmidrule{2-9}	\cmidrule{11-18}				
	&	\multicolumn{3}{c}{\bf \# of Assets}					&		&	\multicolumn{3}{c}{\bf \# of Commodities}					&	\multirow{3}{*}{\minitab[c]{\bf Total \\ \bf Cost}}	&		&	\multicolumn{3}{c}{\bf \# of Assets}					&		&	\multicolumn{3}{c}{\bf \# of Commodities}					&	\multirow{3}{*}{\minitab[c]{\bf Total \\ \bf Cost}}	&		&		&		&												\\	\cmidrule{2-4}	\cmidrule{6-8}	\cmidrule{11-13}	\cmidrule{15-17}		
	&	\multirow{2}{*}{\bf Own}	&	\multirow{2}{*}{\bf Lea}	&	\multirow{2}{*}{\bf Tot}	&		&	\multirow{2}{*}{\bf OT}	&	\multirow{2}{*}{\bf E/T}	&	\multirow{2}{*}{\bf Out}	&		&		&	\multirow{2}{*}{\bf Own}	&	\multirow{2}{*}{\bf Lea}	&	\multirow{2}{*}{\bf Tot}	&		&	\multirow{2}{*}{\bf OT}	&	\multirow{2}{*}{\bf E/T}	&	\multirow{2}{*}{\bf Out}	&		&		&		&		&		&		&								\\						
	&		&		&		&		&		&		&		&		&		&		&		&		&		&		&		&		&		&		&		&		&		&		&								\\	\midrule					
2	&	7	&	-	&	7	&		&	12	&	1/1	&	1	&	225.82	&		&	7	&	-	&	7	&		&	12	&	1/1	&	1	&	234.67	&		&		&		&		&		&								\\						
3	&	7	&	-	&	7	&		&	12	&	2/2	&	4	&	315.84	&		&	7	&	-	&	7	&		&	12	&	2/2	&	4	&	336.39	&		&		&		&		&		&								\\						
6	&	7	&	-	&	7	&		&	6	&	7/7	&	3	&	314.23	&		&	7	&	-	&	7	&		&	9	&	4/3	&	4	&	346.83	&		&		&		&		&		&								\\						
9	&	7	&	-	&	7	&		&	11	&	3/5	&	1	&	250.34	&		&	7	&	-	&	7	&		&	15	&	2/1	&	2	&	265.96	&		&		&		&		&		&								\\						
17	&	7	&	-	&	7	&		&	5	&	5/3	&	2	&	266.65	&		&	7	&	-	&	7	&		&	9	&	2/1	&	3	&	285.86	&		&		&		&		&		&								\\						
18	&	7	&	-	&	7	&		&	9	&	2/3	&	6	&	365.16	&		&	7	&	-	&	7	&		&	9	&	2/3	&	6	&	386.11	&		&		&		&		&		&								\\						
21	&	7	&	-	&	7	&		&	10	&	3/5	&	2	&	268.68	&		&	7	&	-	&	7	&		&	13	&	2/2	&	3	&	295.37	&		&		&		&		&		&								\\						
26	&	7	&	-	&	7	&		&	5	&	4/4	&	2	&	272.4	&		&	7	&	-	&	7	&		&	9	&	1/2	&	3	&	286.2	&		&		&		&		&		&								\\						
27	&	7	&	-	&	7	&		&	11	&	1/4	&	4	&	311.91	&		&	7	&	-	&	7	&		&	11	&	1/4	&	4	&	331.12	&		&		&		&		&		&								\\						
30	&	7	&	-	&	7	&		&	12	&	3/3	&	2	&	264.86	&		&	7	&	-	&	7	&		&	12	&	4/2	&	2	&	289.71	&		&		&		&		&		&								\\						
\hline
		\end{tabular}
	}
\end{table} 
\end{landscape}

The optimal results are also categorized by distance index of the test instances and presented in Table~\ref{tab:Cost-r_summary} as a summary of the lowest and the highest values of the penalty ratio ($r$). Five out of ten instances have the long-range index while the rest of the instances belongs to mid-range category. A total of 90 and 95 commodities are delivered in the former and latter group, respectively. In case of $r=1.2$, the big portion of commodities are delivered as shifted demand (early/tardy) while this dramatically changes in favor of on-time deliveries when $r=10$. The portion of on-time deliveries which is 56.67\% in long-range instances is less than that (63.16\%) of mid-range instances. Thus more commodities are outsourced in the former group of instances. On the network of the mid-range instances, assets can be utilized on shorter distances, thus more commodities can be delivered on-time with the same amount of asset capacity. 
This basically illustrates the effect of network topology on the asset utilization and response action in service network design. 

\begin{table}[htbp]
	\centering
	\scriptsize
	\caption {Summary of sensitivity analysis on penalty ratio, $r$, in terms of network topology.} \label{tab:Cost-r_summary}
	{
		\begin{tabular}{l r r r r r r r r r r r r}
			\midrule
	
\multicolumn{2}{l}{\multirow{3}{*}{\bf Network Topology}}
  		&	\multicolumn{3}{c}{\bf On-time (OT)}					&		&	\multicolumn{3}{c}{\bf Shifted (E/T)}					&		&	\multicolumn{3}{c}{\bf Outsourced (Out.)}					\\	\cmidrule{3-5}	\cmidrule{7-9}	\cmidrule{11-13}
	&		&	\bf $r=1.2$	&	\bf $r=10$	&	\bf $\Delta$	&		&	\bf $r=1.2$	&	\bf $r=10$	&	\bf $\Delta$	&		&	\bf $r=1.2$	&	\bf $r=10$	&	\bf $\Delta$	\\	\midrule		
\multirow{2}{*}{\bf Long-range (LR)}	&	Values	&	26	&	51	&	25	&		&	54	&	21	&	-23	&		&	10	&	18	&	8	\\			
	&	\%	&	28.89\%	&	56.67\%	&	96.15\%	&		&	60\%	&	23.33\%	&	-61.11\%	&		&	11.11\%	&	20\%	&	80\%	\\	\midrule		
\multirow{2}{*}{\bf Mid-range (MR)}	&	Values	&	37	&	60	&	23	&		&	53	&	21	&	-22	&		&	5	&	14	&	9	\\			
	&	\%	&	38.95\%	&	63.16\%	&	62.16\%	&		&	55.79\%	&	22.10\%	&	-60.37\%	&		&	5.26\%	&	14.74\%	&	180\%	\\				
			\hline
		\end{tabular}
	}
\end{table}

\subsection{Impact of Operational Restrictions on Solutions}
\label{subsec-impact-operRest}

Since demand shifting has not been discussed yet in the service network design literature, we conduct a particular test with the same test instances to reveal and to quantify the impact of demand shifting on the optimal solutions. We propose some operational restrictions to limit earliness and tardiness (demand shifting) at a specific portion of the total number of commodities in the problem. 
We created three more problem configurations along with the regular problem setting based on these restrictions such that: (1) no restriction on shifting demand (default case), (2) 50\% restriction on shifting demand, (3) 25\% restriction on shifting demand, and (4) full restriction on shifting demand. 

To formulate the problem configurations stated above, constraint \eqref{SND-One_type_tCom_restricted} is incorporated into the regular formulation such that: (1) $\lambda=1$, (2) $\lambda=0.5$, (3) $\lambda=0.25$, and (4) $\lambda=0$.
\begin{align}
& \sum_{l \in K \cup L: q^{l}\neq 1} p_{l} \leq \lambda |K \cup L|
  \label{SND-One_type_tCom_restricted}
\end{align}

To refine the impacts of other two options, outsourcing and asset leasing, we decreased the leasing cost to $25$, which makes all assets available for upcoming planning horizon at the same fixed operating cost value for the carrier. 

Results of the analysis are summarized in Table~\ref{tab:operational_restrictions} in which each large cell corresponds to a problem configuration on a 2x2 tabular structure. Within each cell, instance name and ID are shown in the first two columns, the number of assets utilized is presented in columns three through five, the number of commodities is given in columns six through eight, and the total cost is displayed in the last column.

As the restriction increases gradually from 0\% to 100\%, commodities that are delivered on-time and via an outsourced service increase as expected. To make more on-time delivery, the carrier has to lease more assets to compensate the lack of demand shifting option in all of the instances after issuing a restriction over 50\%. 

Again, the optimal results are categorized by distance index of the test instances and reported in Table~\ref{tab:Cost-OpeRest_summary} as a summary of the lowest and the highest values of the restriction. The statistics about the long-range and mid-range instances are exactly same with Table~\ref{tab:Cost-r_summary} above. To generalize, the number of assets utilized dramatically increases (by 22.5\%) from going case 1 to case 4 -20 more assets are utilized in the last case - to deliver as many commodities as possible across all test instances. This leads to an increase at 28.03\% in average in the total operating cost for the carrier. In the breakdown of the summarized results by distance index, it is clearly observed that the network topology affects the model in responding excess demand by making a choice among two actions. 98.95\% of all commodities can be delivered on-time in mid-range instances while this value decreases to 94.45\% in long-range problems. The outsourcing option is needed more (5.55\% versus 1.05\%) in the latter. Asset utilization is also more effective on the network of mid-range instances even though the total number of assets utilized in each group are almost equal. Total of 54 assets are utilized to deliver 94 commodities on-time on mid-range networks, in contrast, 55 assets delivers only 85 commodities on-time on long-range networks. This result also supports the argument that network topology is one of the factor that needs to be considered in network design problems.

\begin{table}[htbp]
	\centering
	\scriptsize
	\caption {Summary of sensitivity analysis on operational restrictions in terms of network topology.} \label{tab:Cost-OpeRest_summary}
	{
		\begin{tabular}{l r r r r r r r r r}
			\midrule
\multicolumn{2}{l}{\multirow{3}{*}{\bf Network Topology}}			&	\multicolumn{2}{c}{\bf On-time (OT)}					&		&	\multicolumn{2}{c}{\bf Outsourced (Out.)}					&		&	\multicolumn{2}{c}{\bf \# of Assets}					\\	\cmidrule{3-4}	\cmidrule{6-7}	\cmidrule{9-10}
	&		&	\bf Case 1	&	\bf Case 4			&		&	\bf Case 1	&	\bf Case 4			&		&	\bf Case 1	&	\bf Case 4			\\	\midrule		
\multirow{2}{*}{\bf Long-range (LR)}	&	Values	&	40	&	85			&		&	0	&	5			&		&	47	&	55			\\			
	&	\%	&	44.44\%	&	94.45\%			&		&	N/A	&	5.55\%			&		&	N/A	&	N/A			\\	\midrule		
\multirow{2}{*}{\bf Mid-range (MR)}	&	Values	&	41	&	94			&		&	0	&	1			&		&	42	&	54			\\			
	&	\%	&	43.16\%	&	98.95\%			&		&	N/A	&	1.05\%			&		&	N/A	&	N/A			\\	\midrule		
\bf Total	&	-	&	\bf 81	&	\bf 179			&		&	\bf 0	&	\bf 6			&		&	\bf 89	&	\bf 109			\\			
	
			\hline
		\end{tabular}
	}
\end{table}

\begin{landscape} 
 \begin{table}[htbp]
	\centering
	\scriptsize
	\caption {Sensitivity analysis of selected instances with respect to operational restrictions.} \label{tab:operational_restrictions}
	{
		\renewcommand{\arraystretch}{1}
		\begin{tabular}{l c c c c c c c c r c c c c c c c c r}
			\midrule
\multirow{4}{*}{\bf Instance}	&	\multirow{4}{*}{\bf ID}	&	\multicolumn{8}{c}{\bf $g=25, \text{no restriction}$}															&		&	\multicolumn{8}{c}{\bf $g=25, \text{rest. 50\%}$}															\\	\cmidrule{3-10}	\cmidrule{12-19}				
	&		&	\multicolumn{3}{c}{\bf \# of Assets}					&		&	\multicolumn{3}{c}{\bf \# of Commodities}					&	\multirow{3}{*}{\minitab[c]{\bf Total \\ \bf Cost}}	&		&	\multicolumn{3}{c}{\bf \# of Assets}					&		&	\multicolumn{3}{c}{\bf \# of Commodities}					&	\multirow{3}{*}{\minitab[c]{\bf Total \\ \bf Cost}}	\\	\cmidrule{3-5}	\cmidrule{7-9}	\cmidrule{12-14}	\cmidrule{16-18}		
	&		&	\multirow{2}{*}{\bf Owned}	&	\multirow{2}{*}{\bf Leased}	&	\multirow{2}{*}{\bf Total}	&		&	\bf On	&	\bf Early/	&	\multirow{2}{*}{\bf Out.}	&		&		&	\multirow{2}{*}{\bf Owned}	&	\multirow{2}{*}{\bf Leased}	&	\multirow{2}{*}{\bf Total}	&		&	\bf On	&	\bf Early/	&	\multirow{2}{*}{\bf Out.}	&		\\						
	&		&		&		&		&		&	\bf Time	&	\bf Tardy	&		&		&		&		&		&		&		&	\bf Time	&	\bf Tardy	&		&		\\	\midrule					
inst1.n5.c15	&	2	&	7	&	1	&	8	&		&	7	&	4/4	&	-	&	214.582	&		&	7	&	1	&	8	&		&	8	&	4/3	&	-	&	214.728	\\						
inst1.n5.c20	&	3	&	7	&	3	&	10	&		&	9	&	5/6	&	-	&	270.566	&		&	7	&	3	&	10	&		&	10	&	5/5	&	-	&	270.694	\\						
inst2.n5.c20	&	6	&	7	&	2	&	9	&		&	7	&	7/6	&	-	&	244.256	&		&	7	&	2	&	9	&		&	10	&	5/5	&	-	&	244.734	\\						
inst3.n5.c20	&	9	&	7	&	1	&	8	&		&	6	&	8/6	&	-	&	220.354	&		&	7	&	1	&	8	&		&	12	&	3/5	&	-	&	220.378	\\						
inst6.n5.c15	&	17	&	7	&	2	&	9	&		&	6	&	5/4	&	-	&	240.266	&		&	7	&	2	&	9	&		&	8	&	4/3	&	-	&	240.364	\\						
inst6.n5.c20	&	18	&	7	&	4	&	11	&		&	11	&	4/5	&	-	&	294.704	&		&	7	&	4	&	11	&		&	11	&	4/5	&	-	&	294.704	\\						
inst7.n5.c20	&	21	&	7	&	1	&	8	&		&	8	&	6/6	&	-	&	218.698	&		&	7	&	2	&	9	&		&	12	&	4/4	&	-	&	242.818	\\						
inst9.n5.c15	&	26	&	7	&	2	&	9	&		&	8	&	4/3	&	-	&	238.85	&		&	7	&	2	&	9	&		&	8	&	4/3	&	-	&	238.85	\\						
inst9.n5.c20	&	27	&	7	&	2	&	9	&		&	9	&	5/6	&	-	&	245.416	&		&	7	&	2	&	9	&		&	10	&	4/6	&	-	&	245.868	\\						
inst10.n5.c20	&	30	&	7	&	1	&	8	&		&	10	&	2/8	&	-	&	219.76	&		&	7	&	1	&	8	&		&	10	&	2/8	&	-	&	219.76	\\						
																																					\\	\midrule					
\multirow{4}{*}{\bf Instance}	&	\multirow{4}{*}{\bf ID}	&	\multicolumn{8}{c}{\bf $g=25, \text{rest. 25\%}$}															&		&	\multicolumn{8}{c}{\bf $g=25, \text{rest. 0\%}$}															\\	\cmidrule{3-10}	\cmidrule{12-19}				
	&		&	\multicolumn{3}{c}{\bf \# of Assets}					&		&	\multicolumn{3}{c}{\bf \# of Commodities}					&	\multirow{3}{*}{\minitab[c]{\bf Total \\ \bf Cost}}	&		&	\multicolumn{3}{c}{\bf \# of Assets}					&		&	\multicolumn{3}{c}{\bf \# of Commodities}					&	\multirow{3}{*}{\minitab[c]{\bf Total \\ \bf Cost}}	\\	\cmidrule{3-5}	\cmidrule{7-9}	\cmidrule{12-14}	\cmidrule{16-18}		
	&		&	\multirow{2}{*}{\bf Owned}	&	\multirow{2}{*}{\bf Leased}	&	\multirow{2}{*}{\bf Total}	&		&	\bf On	&	\bf Early/	&	\multirow{2}{*}{\bf Out.}	&		&		&	\multirow{2}{*}{\bf Owned}	&	\multirow{2}{*}{\bf Leased}	&	\multirow{2}{*}{\bf Total}	&		&	\bf On	&	\bf Early/	&	\multirow{2}{*}{\bf Out.}	&		\\						
	&		&		&		&		&		&	\bf Time	&	\bf Tardy	&		&		&		&		&		&		&		&	\bf Time	&	\bf Tardy	&		&		\\	\midrule					
inst1.n5.c15	&	2	&	7	&	1	&	8	&		&	12	&	2/1	&	-	&	216.368	&		&	7	&	2	&	9	&		&	15	&	-	&	-	&	239.45	\\						
inst1.n5.c20	&	3	&	7	&	4	&	11	&		&	15	&	3/2	&	-	&	294.612	&		&	7	&	5	&	12	&		&	20	&	-	&	-	&	320.27	\\						
inst2.n5.c20	&	6	&	7	&	4	&	11	&		&	15	&	3/2	&	-	&	292.272	&		&	7	&	5	&	12	&		&	20	&	-	&	-	&	342.44	\\						
inst3.n5.c20	&	9	&	7	&	1	&	8	&		&	15	&	3/2	&	-	&	221.094	&		&	7	&	2	&	9	&		&	20	&	-	&	-	&	245.58	\\						
inst6.n5.c15	&	17	&	7	&	3	&	10	&		&	12	&	1/2	&	-	&	264.45	&		&	7	&	3	&	10	&		&	14	&	-	&	1	&	289.34	\\						
inst6.n5.c20	&	18	&	7	&	5	&	12	&		&	15	&	3/2	&	-	&	318.972	&		&	7	&	5	&	12	&		&	16	&	-	&	4	&	420.27	\\						
inst7.n5.c20	&	21	&	7	&	2	&	9	&		&	15	&	2/3	&	-	&	244.588	&		&	7	&	4	&	11	&		&	20	&	-	&	-	&	293.1	\\						
inst9.n5.c15	&	26	&	7	&	3	&	10	&		&	12	&	1/2	&	-	&	264.622	&		&	7	&	4	&	11	&		&	15	&	-	&	-	&	290.31	\\						
inst9.n5.c20	&	27	&	7	&	3	&	10	&		&	15	&	1/4	&	-	&	272.412	&		&	7	&	5	&	12	&		&	19	&	-	&	1	&	347.22	\\						
inst10.n5.c20	&	30	&	7	&	2	&	9	&		&	15	&	1/4	&	-	&	244.258	&		&	7	&	4	&	11	&		&	20	&	-	&	-	&	294.19	\\						
\hline
		\end{tabular}
	}
\end{table} 
\end{landscape}

\subsection{Impact of Valid Inequalities}
\label{subsec-impact-VI}


To be able to gain insights about the impact of valid inequalities (VIs) that are described in \S~\ref{subsec-valid_inequal}, we run a particular experiment and test them on all of the medium-size instances. In this analysis, VIs are added into the proposed arc-based formulation once at a time, and three formulation configurations are compared including the regular one when default settings (presolve) of CPLEX are not changed. The performance of VIs are evaluated in terms of CPU time for all instances tested, besides optimality gap as well as best bound obtained for that of instances cannot be solved within the time limit (4 hours) are also compared. All of the results of the VIs' performance test are presented in
Table~\ref{tab:n6_VI_results_landscape} for the set of 
medium-size instances.


 
In Table~\ref{tab:n6_VI_results_landscape}, the first two columns identify instance name and ID in order, then CPU time elapsed for solving each configuration are presented in the next three columns. Optimality gaps for each configuration follow the CPU time in the next three columns and best bound obtained for those of problems that cannot be solved within the time limit are given in the last three columns. Average CPU time and savings in solution time are provided at the bottom of the table, as well. Besides, we give the number of instances that can be solved optimally within the time limit and the number of instances for which the best bound is obtained through the corresponding formulation in the last two rows of the table, respectively.

All three configurations are able to solve only 9 out of 30 problems optimally within the time limit; thus they perform equally in this aspect. 
When it comes to CPU time performance, the third configuration $(P)+\eqref{SND-VI4X-NbAssetperPeriod}$ is the best, which improves the average solution time by 3\% that corresponds to six minutes, approximately. When the average is calculated over those of nine instances solved optimally, the third configuration achieves a 25\% improvement on average (marked with a * at the bottom of the Table~\ref{tab:n6_VI_results_landscape}). The impact of constraints~\eqref{SND-VI4X-NbAssetperPeriod} in some cases is tremendous, for example, formulation $(P)+\eqref{SND-VI4X-NbAssetperPeriod}$ solved instance-50 in one-fifth of the time needed for regular formulation. Likewise, the amount of savings for instance-57, 43 and 34 are 3000, 2000, and 810 seconds, respectively. On the other hand, formulation $(P)+\eqref{SND-VI4X-NbAssetperPeriod}$ is behind of regular formulation only in two (instances-55 and 37) out of nine optimally solved instances in which the differences are about 50 and 200 seconds, respectively.

The second configuration's performance is almost equal to that of basic formulation's performance with a slight difference about 20 seconds in overall average solution time. In four out of nine instances, $(P)+\eqref{SND-VI11-MinNbAsset}$ is faster than $(P)$, while $(P)$ reaches optimal solution sooner in other five instances compared to the second configuration.

In the comparison of obtaining the best (lower) bound in this experiment, the configuration which takes the maximum lower bound for each instance beats the others. The maximum lower bound in each row is emphasized by the bold font in Table~\ref{tab:n6_VI_results_landscape}, besides the number of best bound for each configuration is given at the bottom row. The third configuration outperforms the others for reaching the best bound in more than half of the instances that cannot be solved optimally within the time limit (eleven out of 21 instances). The second configuration comes second with six wins out of 21 instances, and the regular formulation comes last with only five out of 21 instances. The results claim that $(P)+\eqref{SND-VI4X-NbAssetperPeriod}$ is the tightest formulation compared to the other two configurations to get the best lower bound. 

Furthermore, four out of nine problems that are solved optimally within the time limit belong to the close-range category, five out of them are of mid-range category. None of the long-range instances can be solved optimally within the time limit. The highest four optimality gaps, 22.36\%, 17.94\%, 14.22\% and 11.83\%, obtained with regular formulation are observed in solving long-range instances. These outcomes prove that the network topology is among one of the factors on the complexity of CSSND problems and long-range cases are the most complex class within these instances.

 \begin{table}[htbp]
	\centering
	\tiny
	\caption {Comparison of CPU time and optimality gap with valid inequalities on instances with $|N|=6$.} \label{tab:n6_VI_results_landscape}
	{
		\renewcommand{\arraystretch}{1.01}
		\begin{tabular}{l c r r r c r r r r r r r}
			\midrule
\multirow{2}{*}{\bf Instance}	&	\multirow{2}{*}{\bf ID}	&	\multicolumn{3}{c}{\bf CPU Time (sec.)}					&		&	\multicolumn{3}{c}{\bf Optimality Gap (\%)}					&		&	\multicolumn{3}{c}{\bf Best Bound}					\\	\cmidrule{3-5}	\cmidrule{7-9}	\cmidrule{11-13}
	&		&	\multirow{2}{*}{\bf $(P)$}	&	\multirow{2}{*}{\bf $(P)+\eqref{SND-VI11-MinNbAsset}$}	&	\multirow{2}{*}{\bf $(P)+\eqref{SND-VI4X-NbAssetperPeriod}$}	&		&	\multirow{2}{*}{\bf $(P)$}	&	\multirow{2}{*}{\bf $(P)+\eqref{SND-VI11-MinNbAsset}$}	&	\multirow{2}{*}{\bf $(P)+\eqref{SND-VI4X-NbAssetperPeriod}$}	&		&	\multirow{2}{*}{\bf $(P)$}	&	\multirow{2}{*}{\bf $(P)+\eqref{SND-VI11-MinNbAsset}$}	&	\multirow{2}{*}{\bf $(P)+\eqref{SND-VI4X-NbAssetperPeriod}$}	\\			
	&		&		&		&		&		&		&		&		&		&		&		&		\\	\midrule		
inst1.n6.c20	&	31	&	14400	&	14400	&	14400	&		&	4.43	&	4.85	&	3.58	&		&	258.565	&	257.373	&	\bf 260.634	\\			
inst1.n6.c25	&	32	&	14400	&	14400	&	14400	&		&	6.24	&	6.06	&	6.17	&		&	305.295	&	\bf 305.708	&	305.467	\\			
inst1.n6.c30	&	33	&	14400	&	14400	&	14400	&		&	11.83	&	12.92	&	11.99	&		&	\bf 335.929	&	332.491	&	335.39	\\	\midrule		
inst2.n6.c20	&	34	&	2130.62	&	1245.89	&	1320.35	&		&	0	&	0	&	0	&		&	-	&	-	&	-	\\			
inst2.n6.c25	&	35	&	14400	&	14400	&	14400	&		&	4.19	&	4.2	&	4.04	&		&	242.196	&	242.075	&	\bf 242.666	\\			
inst2.n6.c30	&	36	&	14400	&	14400	&	14400	&		&	6.92	&	7.22	&	7.54	&		&	\bf 263.736	&	262.608	&	262.416	\\	\midrule		
inst3.n6.c20	&	37	&	494.476	&	598.56	&	438.144	&		&	0	&	0	&	0	&		&	-	&	-	&	-	\\			
inst3.n6.c25	&	38	&	14400	&	14400	&	14400	&		&	4.77	&	4.01	&	3.86	&		&	260.208	&	262.296	&	\bf 262.588	\\			
inst3.n6.c30	&	39	&	3435.58	&	3923.33	&	3488.34	&		&	0	&	0	&	0	&		&	-	&	-	&	-	\\	\midrule		
inst4.n6.c20	&	40	&	14400	&	14400	&	14400	&		&	5.42	&	5.58	&	5.52	&		&	256.478	&	256.525	&	\bf 256.562	\\			
inst4.n6.c25	&	41	&	14400	&	14400	&	14400	&		&	1.69	&	0.75	&	1.85	&		&	297.939	&	\bf 300.623	&	297.437	\\			
inst4.n6.c30	&	42	&	14400	&	14400	&	14400	&		&	5.29	&	5.16	&	5.44	&		&	316.922	&	\bf 317.375	&	316.972	\\	\midrule		
inst5.n6.c20	&	43	&	8976.84	&	10027.6	&	6790.07	&		&	0	&	0	&	0	&		&	-	&	-	&	-	\\			
inst5.n6.c25	&	44	&	14400	&	14400	&	14400	&		&	4.08	&	3.73	&	3.15	&		&	240.476	&	241.379	&	\bf 242.348	\\			
inst5.n6.c30	&	45	&	14400	&	14400	&	14400	&		&	6.24	&	5.74	&	6.39	&		&	265.451	&	\bf 265.822	&	265.289	\\	\midrule		
inst6.n6.c20	&	46	&	14400	&	14400	&	14400	&		&	5.55	&	5.26	&	5.25	&		&	208.502	&	\bf 209.029	&	208.854	\\			
inst6.n6.c25	&	47	&	14400	&	14400	&	14400	&		&	8.18	&	8.01	&	7.98	&		&	231.946	&	231.764	&	\bf 232.202	\\			
inst6.n6.c30	&	48	&	4624.15	&	6740.69	&	3229.44	&		&	0	&	0	&	0	&		&	-	&	-	&	-	\\	\midrule		
inst7.n6.c20	&	49	&	14400	&	14400	&	14400	&		&	2.48	&	3.05	&	2.77	&		&	\bf 263.367	&	261.844	&	262.597	\\			
inst7.n6.c25	&	50	&	3498.99	&	2418.56	&	763.609	&		&	0	&	0	&	0	&		&	-	&	-	&	-	\\			
inst7.n6.c30	&	51	&	14400	&	14400	&	14400	&		&	5.12	&	4.86	&	4.87	&		&	314.271	&	314.777	&	\bf 315.136	\\	\midrule		
inst8.n6.c20	&	52	&	14400	&	14400	&	14400	&		&	3.63	&	4.69	&	4.38	&		&	\bf 309.138	&	305.604	&	306.602	\\			
inst8.n6.c25	&	53	&	14400	&	14400	&	14400	&		&	17.94	&	11.24	&	16.15	&		&	330.667	&	336.038	&	\bf 337.274	\\			
inst8.n6.c30	&	54	&	14400	&	14400	&	14400	&		&	22.36	&	22.46	&	17.94	&		&	374.872	&	374.847	&	\bf 377.081	\\	\midrule		
inst9.n6.c20	&	55	&	196.046	&	268.93	&	397.771	&		&	0	&	0	&	0	&		&	-	&	-	&	-	\\			
inst9.n6.c25	&	56	&	3626.18	&	3109.22	&	3372.04	&		&	0	&	0	&	0	&		&	-	&	-	&	-	\\			
inst9.n6.c30	&	57	&	12707.5	&	12083.7	&	9705.63	&		&	0	&	0	&	0	&		&	-	&	-	&	-	\\	\midrule		
inst10.n6.c20	&	58	&	14400	&	14400	&	14400	&		&	6.21	&	6.46	&	6.19	&		&	276.549	&	275.853	&	\bf 276.596	\\			
inst10.n6.c25	&	59	&	14400	&	14400	&	14400	&		&	3.37	&	3.51	&	3.47	&		&	\bf 316.84	&	315.471	&	316.582	\\			
inst10.n6.c30	&	60	&	14400	&	14400	&	14400	&		&	14.22	&	18.88	&	19.7	&		&	350.919	&	\bf 351.275	&	347.602	\\	\midrule		
\multicolumn{2}{l}{\bf Average}			&	11403.01	&	11427.22	&	11063.51	&		&	5.01	&	4.91	&	4.94	&		&	-	&	-	&	-	\\			
\multicolumn{2}{l}{\bf Improvement (\%) }			&	-	&	-0.21	&	2.97	&		&	-	&	-	&	-	&		&	-	&	-	&	-	\\	\midrule		
\multicolumn{2}{l}{\bf Average*}			&	4410.04	&	4490.72	&	3278.38	&		&	-	&	-	&	-	&		&	-	&	-	&	-	\\			
\multicolumn{2}{l}{\bf Improvement* (\%) }			&	-	&	-1.83	&	25.66	&		&	-	&	-	&	-	&		&	-	&	-	&	-	\\	\midrule		
\multicolumn{2}{l}{\bf \# of inst. ($<$4h)}			&	-	&	-	&	-	&		&	9	&	9	&	9	&		&	-	&	-	&	-	\\			
\multicolumn{2}{l}{\bf \# of  best bounds}			&	-	&	-	&	-	&		&	-	&	-	&	-	&		&	5	&	6	&	10	\\			
\hline
		\end{tabular}
	}
\end{table} 
\section{Dedicate-Merge-and-Mix Algorithm}
\label{sec-DMaM}
As discussed several times in \S~\ref{sec-computational_exp}, MIP formulation requires very long computational time to solve the CSSND problem optimally for medium-size instances even. From the business perspective, MIP formulation may fail to solve the realistic size of instances in relatively short time for a decision-maker. To overcome these drawbacks and to find high-quality solutions for any size of instance instantly, we propose a three-step heuristic, Dedicate-merge-and-mix algorithm (DMaM), in this paper. Details of the algorithm are described in the following subsections, and the flowchart of the algorithm is illustrated in Figure~\ref{fig:DMaMflow} at the end of the section.

\subsection{Description of the DMaM}
\label{subsec-DMaM-steps}

DMaM is developed as a multi-phase solution procedure in which the decisions of the problem including commodity flow, routing of assets, demand shifting, asset leasing, and outsourcing are incorporated step-by-step. The procedure starts with the input processing phase to generate all feasible paths for all TCs from their origin to their destination, that succeeds an initialization phase to read the problem data at the beginning. The algorithm proceeds with the construction phase in which we aim to obtain a feasible initial solution for the CSSND problem. In the third and computationally the most extensive phase, we improve the obtained solution through minimizing the number of commodity paths (number of assets utilized implicitly) by merging paths in a two-step searching process and by mixing paths into another path. At the end, the procedure checks the capacity shortage one more time finally and resolves any shortage by leasing additional asset(s). Last, the procedure is finalized, and all results are reported.

\subsubsection{Phase I: Commodity Path Generation}
\label{subsec-DMaM-pathGen}
In this phase, a subprocedure enumerates all possible paths for all TCs on the service network based on the set of arcs ($A$), commodity's origin node and release date ($O_{k}^{2}$) as well as commodity's destination node and due date ($D_{k}^{2}$). Any generated path, as well as its calculated cost, are recorded in \textit{the generated path list} as an input for the following subprocedures.  

Since we consider a complete service network in the problem setting, there is a service arc between any pair of physical nodes on any period based on the distance between these nodes. Thus, the distinction between commodity paths relies on the fact that how many holding arcs exist on any path. There are three types of commodity paths with six distinct patterns available, which are all listed in Table~\ref{tab:comPaths} where $ S $ and $ H $ stand for \textit{service arc} and \textit{holding arc}, respectively.

\begin{table}[htbp]
	\centering
	\scriptsize
	\caption {Possible commodity path patterns.} \label{tab:comPaths}
	{
		\begin{tabular}{l c c c}
			\midrule			
\bf Commodity path pattern	&		&	\bf \# of holding arcs	&	\bf Service type	\\	\midrule
$ O^{'}_{k} \xrightarrow{} S \xrightarrow{} D^{'}_{k} $	&		&	0	&	Offered \&  outsourced	\\	\midrule
$ O^{'}_{k} \xrightarrow{} S \xrightarrow{} H \xrightarrow{} D^{'}_{k} $	&		&	\multirow{2}{*}{1}	&	\multirow{2}{*}{Offered}	\\	
$ O^{'}_{k} \xrightarrow{} H \xrightarrow{} S \xrightarrow{} D^{'}_{k} $	&		&		&		\\	\midrule
$ O^{'}_{k} \xrightarrow{} S \xrightarrow{} H \xrightarrow{} H \xrightarrow{} D^{'}_{k} $	&		&	\multirow{3}{*}{2}	&	\multirow{3}{*}{Offered}	\\	
$ O^{'}_{k} \xrightarrow{} H \xrightarrow{} S \xrightarrow{} H \xrightarrow{} D^{'}_{k} $	&		&		&		\\	
$ O^{'}_{k} \xrightarrow{} H \xrightarrow{} H \xrightarrow{} S \xrightarrow{} D^{'}_{k} $	&		&		&		\\		
			\midrule	
		\end{tabular}
	}
\end{table}

For a few of the commodities, delivery is not possible on any offered services due to lack of enough time between release and due date. Thus, path generation subprocedure ends up with only generated paths based on outsourced services for those of commodities.

\subsubsection{Phase II: Constructing An Initial Feasible Solution}
\label{subsec-DMaM-construct}

The algorithm proceeds to the second phase where it constructs an initial feasible solution to the problem through a multi-step subprocedure when all possible commodity paths are given as an input. The subprocedure starts with searching the path with the minimum routing cost among the all generated ones in a greedy way and selects it as an active path for each OC in the first step. Thus, one of the problem features, demand shifting, is incorporated through the path selection process in the proposed algorithm with the cost perspective because the subprocedure compares the cost of all paths belong to each TC but selects only one for each OC. 

As stated above, few commodities cannot be delivered through offered services; therefore, there would not be any path generated on offered services. For this reason, the only path with an outsourced service is selected for those of commodities. Hence, the other problem feature, outsourcing, is also considered at this point in the scope of the algorithm. When we select a path either based on offered services or outsourced services for each OC, an initial feasible solution is constructed. Consequently, the number of total paths in the initial solution equals the number of OCs in the problem instance. The number of paths on offered services gives us the number of assets needed to transport these OCs. The basic assumption is that each asset is assigned to each path on offered services. Assets are assumed to be busy from path origin to path destination; however, they make empty moves from path destination to path origin in their cycle during the planning horizon.  

Last, the procedure calculates the total operating cost of the obtained solution based on the number of paths via offered services (fixed cost) and the cost of each commodity path (routing cost). Since the initial solution is a very relaxed solution with many empty moves, the number of assets most likely would be higher than the number of owned assets. Thus, it is assumed that any additional asset is leased as many as the difference between the number of asset cycles and the number of owned assets. Hence, the last feature of the proposed problem, leasing (acquiring) additional assets, is incorporated at this step as well as the leasing cost of additional assets in the cost calculation. The cost of demand shifting (early/tardy delivery) is already included through routing cost of commodity paths since there might be paths of shifted demand selected in the first step of this procedure. Likewise, the final cost element - outsourcing cost - can be attained through routing cost of some commodity paths which are also selected in the first step for those of commodities that cannot be delivered through offered services.

\subsubsection{Phase III: Improving Solution by Merging Paths}
\label{subsec-DMaM-improve}

In the third phase, we aim to improve the initial solution basically by merging two commodity paths into one asset cycle. The number of empty moves decreases and the utilization of assets increases, when an asset cycle is eliminated from the solution in each merge operation. Another multi-step subprocedure is developed to explore and execute any merging option in this phase of the proposed DMaM. 

The subprocedure starts with sorting the selected paths in the initial solution by the number of periods in which the asset assigned to the corresponding path is busy or in other words, how long the commodity path takes from TC's origin to destination. Then, the sorted set of paths is partitioned into two distinct sets as \textit{primary paths} and \textit{secondary paths} such that the paths whose number of busy periods is greater than that of idle periods, are classified as primary paths while others (number of busy periods is less than or equal to that of idle periods) are classified as secondary paths. The reason for sorting and partitioning is that subprocedure should start exploring with primary paths because they are almost completely busy and have less potential to be merged with another path in the solution. The process of searching on ordered paths is demonstrated in Table~\ref{tab:searchOrder} in a generalized way. Given that the planning horizon composes of seven periods, and an arbitrary primary commodity path takes five periods, the merging subprocedure searches on two types of secondary paths which have two or one busy periods by referring to the second row of Table~\ref{tab:searchOrder}. The algorithm proceeds to search for merging a pair of two secondary paths when it completes the searching process on the combination of one primary and one secondary path.

\begin{table}[htbp]
	\centering
	\scriptsize
	\caption {Effective search order in exploring possible merge of paths.} \label{tab:searchOrder}
	{
		\begin{tabular}{l c l c c c l c c}
			\midrule			
\multirow{3}{*}{\bf Step}	&		&	\multicolumn{3}{c}{\bf Commodity Path One}					&		&	\multicolumn{3}{c}{\bf Commodity Path Two}					\\	\cmidrule{3-5}	\cmidrule{7-9}
	&		&	\bf Type	&	\bf Busy Period	&	\bf Idle Period	&		&	\bf Type	&	\bf Busy Period	&	\bf Idle Period	\\	\midrule	
$ 1^{st} $	&		&	Primary	&	$|T| -1 $	&	$ 1 $	&		&	Secondary	&	$ 1 $	&	$|T| -1 $	\\	\midrule	
\multirow{2}{*}{$ 2^{nd} $}	&		&	\multirow{2}{*}{Primary}	&	\multirow{2}{*}{$|T| - 2 $}	&	\multirow{2}{*}{$ 2 $}	&		&	\multirow{2}{*}{Secondary}	&	$ 2 $	&	$|T| -2 $	\\		
	&		&		&		&		&		&		&	$ 1 $	&	$|T| -1 $	\\	\midrule	
\multirow{3}{*}{$ 3^{rd} $}	&		&	\multirow{3}{*}{Primary}	&	\multirow{3}{*}{$|T| - 3 $}	&	\multirow{3}{*}{$ 3 $}	&		&	\multirow{3}{*}{Secondary}	&	$ 3 $	&	$|T| -3 $	\\		
	&		&		&		&		&		&		&	$ 2 $	&	$|T| -2 $	\\		
	&		&		&		&		&		&		&	$ 1 $	&	$|T| -1 $	\\	\midrule	
\ldots	&		&	\ldots	&	\ldots	&	\ldots	&		&	\ldots	&	\ldots	&	\ldots	\\	\midrule	
$ z^{th} :  |T| - z > z $	&		&	Primary	&	$|T| - z $	&	$ z $	&		&	Secondary	&	$ z $	&	$|T| -2 $	\\	\midrule	
$ (z+1)^{th} :  z \leq |T| - z $	&		&	Secondary	&	$|T| -(z+1) $	&	$ (z+1) $	&		&	Secondary	&	$ (z+1) $	&	$|T| -(z+1) $	\\	\midrule	
\ldots	&		&	\ldots	&	\ldots	&	\ldots	&		&	\ldots	&	\ldots	&	\ldots	\\	\midrule	
\multirow{3}{*}{$|T| -3^{th} $}	&		&	\multirow{3}{*}{Secondary}	&	\multirow{3}{*}{$ 3 $}	&	\multirow{3}{*}{$|T| -3 $}	&		&	\multirow{3}{*}{Secondary}	&	$ 3 $	&	$|T| -3 $	\\		
	&		&		&		&		&		&		&	$ 2 $	&	$|T| -2 $	\\		
	&		&		&		&		&		&		&	$ 1 $	&	$|T| -1 $	\\	\midrule	
\multirow{2}{*}{$ |T| -2^{th} $}	&		&	\multirow{2}{*}{Secondary}	&	\multirow{2}{*}{$ 2 $}	&	\multirow{2}{*}{$ |T| -2 $}	&		&	\multirow{2}{*}{Secondary}	&	$ 2 $	&	$|T| -2 $	\\		
	&		&		&		&		&		&		&	$ 1 $	&	$|T| -1 $	\\	\midrule	
$|T| -1^{th} $	&		&	Secondary	&	$ 1 $	&	$|T| -1 $	&		&	Secondary	&	$ 1 $	&	$|T| -1 $	\\		

			\midrule	
		\end{tabular}
	}
\end{table}

To determine that a merge operation can be executed on a pair of commodity paths, we test few basic parameters -physical origin and destination nodes as well as release and due dates- of two TCs belong to those paths both spatially and time-wise in this order whether merging conditions are satisfied or not. Before explaining what merging conditions are, it is better to discuss prerequisites (prior conditions) of the exploration process. To be able to compare the aforementioned parameters properly in exploration, we check time-wise parameters (release and due dates) according to prior conditions and adjust these if prior conditions are not satisfied. These parameters can be used as they are in the exploration process if the conditions are satisfied; otherwise, the ones that exist in the failed condition are adjusted as presented in the third column of Table~\ref{tab:prerequisite}.

\begin{table}[htbp]
	\centering
	\scriptsize
	\caption {Prerequisite for exploration process and proper adjustments on time-wise parameters.} \label{tab:prerequisite}
	{
		\begin{tabular}{l c c}
			\midrule			
\bf No.	&	\bf Prerequisite (Pr)	&	\bf Adjustment if Pr fails	\\	\midrule
1	&	$t_{O_{1}} < t_{D_{1}} $	&	$ t_{D_{1}} =  t_{D_{1}} + |T|$	\\	\midrule
2	&	$t_{O_{2}} < t_{D_{2}} $	&	$ t_{D_{2}} =  t_{D_{2}} + |T|$	\\	\midrule
\multirow{2}{*}{3}	&	\multirow{2}{*}{$t_{O_{1}} < t_{O_{2}} $}	&	$ t_{O_{2}} =  t_{O_{2}} + |T|$	\\	
	&		&	$ t_{D_{2}} =  t_{D_{2}} + |T|$	\\	

			\midrule	
		\end{tabular}
	}
\end{table}

The exploration process consists of two types of searching process for a possible merge operation, \textit{search for regular merging} and \textit{search for shifted merging}. Two commodity paths which are investigated for a merge operation can be directly merged in the former one. However, one of the paths/or both that are in question is/are supposed to be shifted in the latter one. To clarify, shifting basically means switching from one commodity path to another (switching from one type of TC to another, i.e., original TC to early TC or early TC to tardy TC), in other words shifting the demand either earlier or later date is required for executing a merge operation. 

\begin{figure}[htbp]
\centering
\fbox{\begin{varwidth}{1\textwidth}
\subfigure[W/o Repositioning.]{
\includegraphics[clip,width=0.35\textwidth]{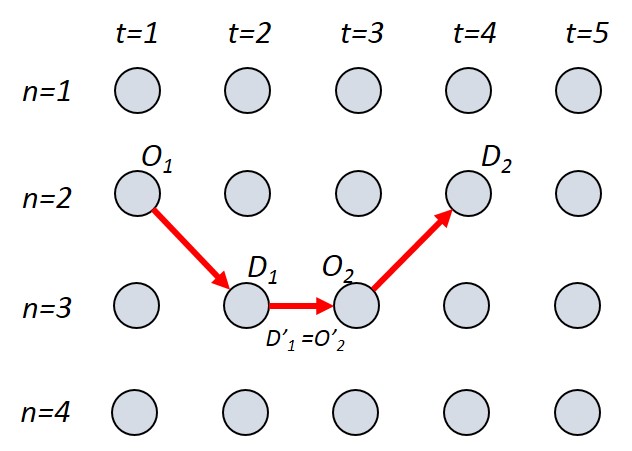}\label{fig:merge1}}
\quad\subfigure[W/ 2-Repositioning.]{
\includegraphics[clip,width=0.35\textwidth]{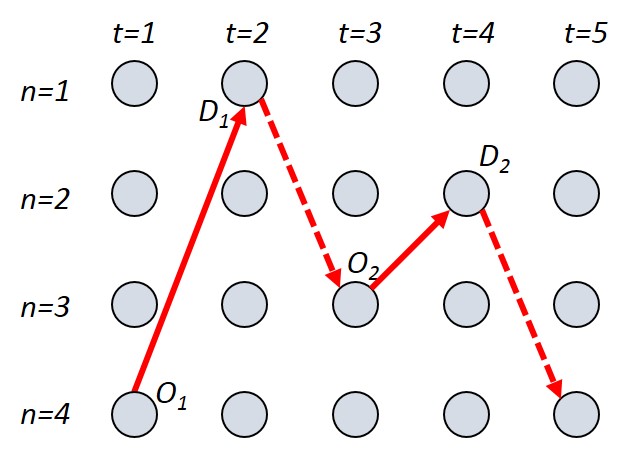}\label{fig:merge2}}\\
\subfigure[W/ 1-Repositioning V1.]{
\includegraphics[clip,width=0.35\textwidth]{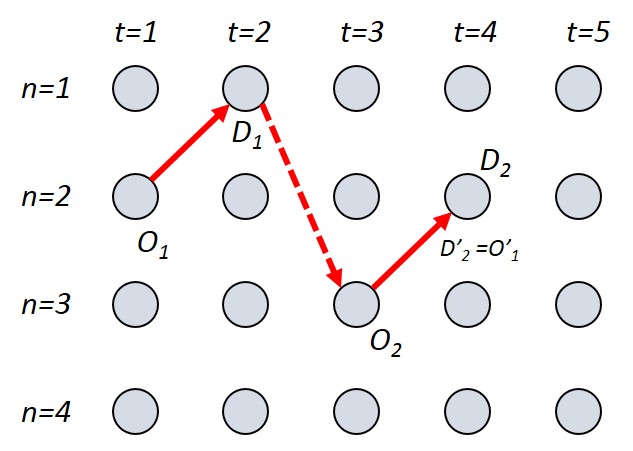}\label{fig:merge3}}
\quad\subfigure[W/ 1-Repositioning V2.]{
\includegraphics[clip,width=0.35\textwidth]{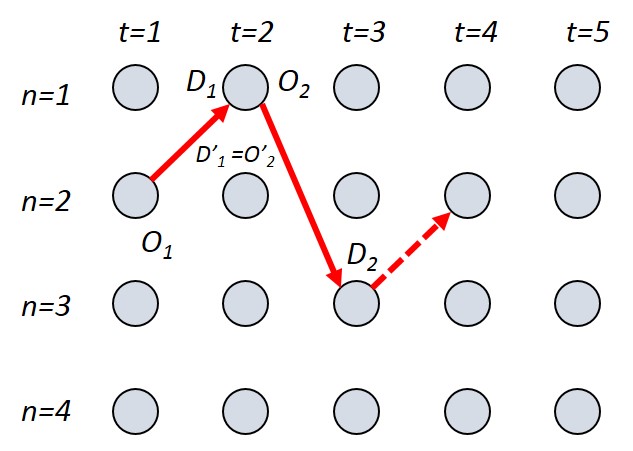}\label{fig:merge4}}
\end{varwidth}}
\caption{Some solutions for the multiple allocation problem.}
\label{fig:mergeTypes}
\end{figure}

There are four types of merge operations which are presented through rows of Table~\ref{tab:mergeConditions} and illustrated in Figure~\ref{fig:mergeTypes} for both merge searching process (regular/shifted). The first type is defined as merging two paths without adding a repositioning arc in the merged path and labeled as `w/o Rep`. The physical origin node of path one is same with the physical destination of path two, and the physical destination node of path one is the same with physical origin node of path two, that means these commodities travel in opposite directions and constitutes a perfect match (see Figure~\ref{fig:merge1}). In the second one that is labeled as `w/ 1-Rep V1`, the physical origin node of path one is same with the physical destination of path two, while the physical destination node of path one is not the same as physical origin node of path two, thus a repositioning arc must be added from the physical destination node of path one to physical origin node of path two (see Figure~\ref{fig:merge3}). The third type (`w/ 1-Rep V2`) is the other way round of merge type two where a repositioning arc must be added from physical destination node of path two to physical origin node of path one (see Figure~\ref{fig:merge4}). Last, two repositioning arcs must be added into the path right after each path to be able to merge two commodity paths in this type that is named as `w/ 2-Rep` (see Figure~\ref{fig:merge2}). Besides, $ t_{[O_{k}]} $ indicates the time index of beginning of the path that belongs to commodity $k$ in next week (i.e, time index of end of any asset cycle) and always equal to $ t_{[O_{1}]}= t_{O_{1}}+ |T| $.

\begin{table}[htbp]
	\centering
	\tiny
	\caption {Merge types and conditions for regular and shifted merging.} \label{tab:mergeConditions}
	{
		\begin{tabular}{l l c l cl}
			\midrule			
\multirow{3}{*}{\bf Merge type}	&	\multirow{3}{*}{\bf Explanation}	&	\multicolumn{2}{c}{\bf Conditions for regular merge}			&	\multicolumn{2}{c}{\bf Conditions for shifted merge}			\\	\cmidrule{3-4}	\cmidrule{5-6}
	&		&	\bf Spatial	&	\bf Time-wise	&	\bf Spatial	&	\bf Time-wise	\\	\midrule	
\multirow{2}{*}{w/o Rep}	&	Merge without repositioning 	&	$ O^{'}_{1} = D^{'}_{2}  $	&	$ t_{D_{1}} \leq t_{O_{2}} $	&	$ O^{'}_{1} = D^{'}_{2}  $	&	$ t_{D_{1}} + a_{m1} \leq t_{O_{2}} + a_{m2} $	\\		
	&	arc (perfect match)	&	$ D^{'}_{1} = O^{'}_{2}  $	&	$ t_{D_{2}} \leq t_{[O_{1}]} $	&	$ D^{'}_{1} = O^{'}_{2}  $	&	$ t_{D_{2}} + a_{m2} \leq t_{[O_{1}]} + a_{m1} $	\\	\midrule	
\multirow{3}{*}{w/ 1-Rep V1}	&	\multirow{1.5}{*}{Merge with one repositioning}	&	$ O^{'}_{1} = D^{'}_{2}  $	&	$ t_{D_{2}} \leq t_{[O_{1}]} $	&	$ O^{'}_{1} = D^{'}_{2}  $	&	$ t_{D_{2}} + a_{m2} \leq t_{[O_{1}]} + a_{m1} $	\\		
	&	\multirow{1.5}{*}{arc (repositioning after first)}	&	$ D^{'}_{1} \ne O^{'}_{2}  $	&	$ t_{D_{1}} < t_{O_{2}} $	&	$ D^{'}_{1} \ne O^{'}_{2}  $	&	$ t_{D_{1}} + a_{m1} < t_{O_{2}} + a_{m2} $	\\		
	&		&	$ D^{'}_{1} \rightarrow O^{'}_{2}  $	&	$ d_{D^{'}_{1}, O^{'}_{2}} \leq t_{O_{2}}-t_{D_{1}} $	&	$ D^{'}_{1} \rightarrow O^{'}_{2}  $	&	$ d_{D^{'}_{1}, O^{'}_{2}} \leq t_{O_{2}}+ a_{m2}-t_{D_{1}} + a_{m1} $	\\	\midrule	
\multirow{3}{*}{w/ 1-Rep V2}	&	\multirow{1.5}{*}{Merge with one repositioning}	&	$ O^{'}_{1} \neq D^{'}_{2}  $	&	$ t_{D_{1}} \leq t_{O_{2}} $	&	$ O^{'}_{1} \neq D^{'}_{2}  $	&	$ t_{D_{1}} + a_{m1} \leq t_{O_{2}} + a_{m2} $	\\		
	&	\multirow{1.5}{*}{arc (repositioning after second)}	&	$ D^{'}_{1} = O^{'}_{2}  $	&	$ t_{D_{2}} < t_{[O_{1}]} $	&	$ D^{'}_{1} = O^{'}_{2}  $	&	$ t_{D_{2}} + a_{m2}< t_{[O_{1}]} + a_{m1}$	\\		
	&		&	$ D^{'}_{2} \rightarrow O^{'}_{1}  $	&	$ d_{D^{'}_{2}, O^{'}_{1}} \leq t_{[O_{1}]}-t_{D_{2}} $	&	$ D^{'}_{2} \rightarrow O^{'}_{1}  $	&	$ d_{D^{'}_{2}, O^{'}_{1}} \leq t_{[O_{1}]}+ a_{m1}-t_{D_{2}} + a_{m2} $	\\	\midrule	
\multirow{4}{*}{w/ 2-Rep}	&	\multirow{3}{*}{Merge with two}	&	$ O^{'}_{1} \neq D^{'}_{2}  $	&	$ t_{D_{1}} < t_{O_{2}} $	&	$ O^{'}_{1} \neq D^{'}_{2}  $	&	$ t_{D_{1}}+ a_{m1} < t_{O_{2}} + a_{m2} $	\\		
	&	\multirow{3}{*}{repositioning arcs}	&	$ D^{'}_{1} \ne O^{'}_{2}  $	&	$ t_{D_{2}} < t_{[O_{1}]} $	&	$ D^{'}_{1} \ne O^{'}_{2}  $	&	$ t_{D_{2}} + a_{m2} < t_{[O_{1}]}+ a_{m1} $	\\		
	&		&	$ D^{'}_{2} \rightarrow O^{'}_{1}  $	&	$ d_{D^{'}_{2}, O^{'}_{1}} \leq t_{[O_{1}]}-t_{D_{2}} $	&	$ D^{'}_{2} \rightarrow O^{'}_{1}  $	&	$ d_{D^{'}_{2}, O^{'}_{1}} \leq t_{[O_{1}]}+ a_{m1}-t_{D_{2}} + a_{m2} $	\\		
	&		&	$ D^{'}_{1} \rightarrow O^{'}_{2}  $	&	$ d_{D^{'}_{1}, O^{'}_{2}} \leq t_{O_{2}}-t_{D_{1}} $	&	$ D^{'}_{1} \rightarrow O^{'}_{2}  $	&	$ d_{D^{'}_{1}, O^{'}_{2}} \leq t_{O_{2}} + a_{m2}-t_{D_{1}}+ a_{m1} $	\\		
	
			\midrule	
		\end{tabular}
	}
\end{table}

The aforementioned spatial and time-wise merging conditions to decide whether two investigated paths can be merged or not are presented in Table~\ref{tab:mergeConditions} in columns three and four for regular merging, and in columns for five and six for shifted merging. Conditions for shifted merging are slightly different from that of regular merging as given in the table. Note that spatial conditions are checked first when searching for a merge operation, time-wise conditions are only tested given that spatial conditions are satisfied for a pair of paths. 

In deciding a regular merge operation, spatial and time-wise conditions given in Table~\ref{tab:mergeConditions} (columns three and four) are sufficient, while we need additional analysis to determine whether two commodity paths can be merged through shifting or not. In this step, we apply a \textit{maximum possible shifting period analysis} for any pair of paths and determine that a shifted merging is possible. We define slack variables for each time-wise condition given in Table~\ref{tab:mergeConditions} such as presented in Table~\ref{tab:slackVar}. Basically, these variables indicate how many periods that commodity paths need to be shifted to be able to merge corresponding paths. In the case of regular merging, these slack variables are within the domain given in column three in Table~\ref{tab:slackVar}, two paths investigated can be merged through regular merging without shifting. In contrast, slack variables take positive values and maximum of available slack variables based on merge type, $s_{max}=Max \left \{s_{1}, .., s_{n} \right \}$, is tested whether that is less than two or not because any value greater than two means that shifted merging operation is not feasible. To remember, one of the problem assumptions states that demand shifting is only possible for one period earlier or later. Thus, the maximum total number of periods to shift is two - when one path is shifted backward, and the other is shifted forward in time - for the problem. The rest of the alternatives in shifting investigated paths is demonstrated in Table~\ref{tab:6alternatives}.

\begin{table}[htbp]
	\centering
	\scriptsize
	\caption {Slack variables as well as domains for shifted merging.} \label{tab:slackVar}
	{
		\begin{tabular}{l l l | c l l l}
			\midrule			
\bf Merge type	&	\bf Slack variables	&	\bf Domain	&		&	\bf Merge type	&	\bf Slack variables	&	\bf Domain	\\	\midrule
\multirow{2}{*}{w/o Rep}	&	$ s_{1}=t_{D_{1}} - t_{O_{2}} $	&	$ s_{1} \leq 0 $	&		&	\multirow{3}{*}{w/ 1-Rep V1}	&	$ s_{1}=t_{D_{2}} - t_{[O_{1}]} $	&	$ s_{1} \leq 0 $	\\	
	&	$ s_{2}=t_{D_{2}} - t_{[O_{1}]} $	&	$ s_{2} \leq 0$	&		&		&	$ s_{2}=t_{D_{1}} - t_{O_{2}} $	&	$ s_{2} +1 \leq 0$	\\	
	&		&		&		&		&	$ s_{3}=d_{D^{'}_{1}, O^{'}_{2}} - (t_{O_{2}}+t_{D_{1}}) $	&	$ s_{3} \leq 0 $	\\	\midrule
\multirow{4}{*}{w/ 2-Rep}	&	$ s_{1}=t_{D_{1}} - t_{O_{2}} $	&	$ s_{1} +1 \leq 0$	&		&	\multirow{3}{*}{w/ 1-Rep V2}	&	$ s_{1}=t_{D_{1}} - t_{O_{2}} $	&	$ s_{1} \leq 0 $	\\	
	&	$ s_{2}=t_{D_{2}}  t_{[O_{1}]} $	&	$ s_{2} +1 \leq 0$	&		&		&	$ s_{2}=t_{D_{2}} - t_{[O_{1}]} $	&	$ s_{2} +1 \leq 0$	\\	
	&	$ s_{3}=d_{D^{'}_{2}, O^{'}_{1}} - (t_{[O_{1}]}+t_{D_{2}}) $	&	$ s_{3} \leq 0 $	&		&		&	$ s_{3}=d_{D^{'}_{2}, O^{'}_{1}} - (t_{[O_{1}]}+t_{D_{2}}) $	&	$ s_{3} \leq 0 $	\\	
	&	$ s_{4}=d_{D^{'}_{1}, O^{'}_{2}} - (t_{O_{2}}+t_{D_{1}}) $	&	$ s_{4} \leq 0 $	&		&		&		&		\\	

			\midrule	
		\end{tabular}
	}
\end{table}

Next step in the shifted merging is to determine which shifting alternative fits paths that are in question. Spatial and time-wise conditions given in column five and six in Table~\ref{tab:mergeConditions} are updated based on determined alternative (see Table~\ref{tab:6alternatives}), then shifted merge operation is executed if updated conditions are satisfied. In case of more than one shifted merging alternative fits with the paths, the merge with the minimum routing cost is chosen for execution. The algorithm proceeds to the next pair of paths when the conditions are not satisfied for any other alternative.

\begin{table}[htbp]
	\centering
	\scriptsize
	\caption {Six shifting alternatives for shifted merging.} \label{tab:6alternatives}
	{
		\begin{tabular}{c c l c l c}
			\midrule			
\bf Alternative (m)	&	\bf $ s_{max}$	&	\bf Path One	&	\bf $ a_{m1} $	&	\bf Path Two	&	\bf $ a_{m2} $	\\	\midrule
1	&	1	&	Forward	&	1	&	-	&	-	\\	
2	&	1	&	Backward	&	-1	&	-	&	-	\\	
3	&	1	&	-	&	-	&	Forward	&	1	\\	
4	&	1	&	-	&	-	&	Backward	&	-1	\\	\midrule
5	&	2	&	Forward	&	1	&	Backward	&	-1	\\	
6	&	2	&	Backward	&	-1	&	Forward	&	1	\\	
			\midrule	
		\end{tabular}
	}
\end{table}

At this step of the DMaM, we propose three distinct configurations for the merging process such as random search, custom search, and advanced search. \textit{The DMaM with random search} ($DMaM-R$) executes any merge operation immediately after it explores any pair of paths, as explained in the previous paragraph. The exploration and the execution are performed subsequently in the $DMaM-R$. However, the other two configurations of DMaM seek for the maximum number of pairs that can be merged in Phase-III because it is most likely to miss some of merging pairs in the random search. As indicated several times in this paper, the fixed cost of assets is dominant in the total operating cost which decreases as the number of assets decreases in this problem. Thus, the possible smallest number of assets can be achieved by merging as many paths as possible in the execution step. \textit{DMaM with custom search} ($DMaM-C$) and \textit{DMaM with advanced search} ($DMaM-A$) are developed to reach the maximum number of pair of paths in exploring step by using alternative methods in the algorithm. 

In \textit{DMaM with custom search} ($DMaM-C$), we introduce \textit{smallest-conflicted-pairs-first (SCoPF)} algorithm to determine maximum possible number of path pairs that can be merged among all pairs explored in the searching process. The SCoPF is thoroughly explained by pseudo code in Algorithm~\ref{alg:SCoPF}. 

\begin{algorithm}
    \small
	\caption{Smallest conflicted pairs first algorithm (SCoPF)} \label{alg:SCoPF}
	\begin{algorithmic}[1]
        \State \textit{Initialization when all pair of paths that can be merged explored.}
        \State \textit{Create a two-dimensional square array with a length of total number of unique paths which exist in explored pairs: "conflict array".}
        \State \textit{Map all pairs of path on the conflict array by putting one at the intersection of paths.}
        \State \textit{Sum up ones in each row and column, record these values as an individual score for each row/column.}
        \State \textit{Replace ones at the intersection of pair of paths with sum of individual scores of paths in any pair.}
        \State \textit{Choose the pair of path with the minimum conflict score to merge and cancel out corresponding row and column from the array.}
        \State \textit{Repeat step-6 until no more score greater than zero is cancelled out on the conflicted array.}
        \State \textit{Move to next step in the DMaM-C: execution of merge operation on chosen pairs.}
	\end{algorithmic}
\end{algorithm}

As an alternative to $DMaM-C$, we also propose another configuration of $DMaM$ in which we aim to guarantee to obtain the maximum number of pairs among the explored pairs. To achieve this, we develop \textit{maximum pairs at minimum routing cost algorithm (MaPMiC)} based on solving Bi-partite graph matching problem ($P_{2}$) optimally in an iterated procedure in \textit{DMaM with advanced search} ($DMaM-A$). The MaPMiC is thoroughly explained by using the pseudo-code in Algorithm~\ref{alg:MaPMiC}.    

\begin{algorithm}
    \small
	\caption{Maximum pairs at minimum routing cost algorithm (MaPMiC)} \label{alg:MaPMiC}
	\begin{algorithmic}[1]
        \State \textit{Initialization when all pair of paths that can be merged explored.}
        \State \textit{Create a two-dimensional square array with a length of total number of unique paths which exist in explored pairs: "conflict array".}        
        \State \textit{Map all pairs of path on the conflict array by putting one at the intersection of paths.}
        \State \textit{Create another two-dimensional square array with a length of total number of unique paths which exist in explored pairs: "routing cost array".}
        \Procedure{Do while: it is infeasible}{}		
		\begin{itemize}
		\item \textit{Solve Bi-partite graph matching problem ($P_{2}$) by feeding the conflict array and the routing cost array as input.} 
		\item \textit{Decrease right hand side of the constraint (\ref{BIP-MatchingNb}) one by one.}
		\end{itemize}
		\EndProcedure
        \State \textit{Move to next step in the DMaM-A: execution of merge operation on chosen pairs.}
	\end{algorithmic}
\end{algorithm}

The Bi-partite graph matching problem is formulated as follows:
\begin{align}
\mbox{\textbf{($P_{2}$)}}
\qquad Min \quad & \sum_{i \in M} \sum_{j \in M: \lambda_{ij}=1} e_{ij} z_{ij}   \label{BIP-Obj} \\
s. \: t. \quad  & \sum_{j \in M: \lambda_{ij}=1} z_{ij} + \sum_{j \in M: \lambda_{ji}=1} z_{ji} \leq 1 \qquad \forall \: i \in M \label{BIP-Originwise} \\
                & \sum_{i \in M: \lambda_{ij}=1} z_{ij} + \sum_{i \in M: \lambda_{ji}=1} z_{ji} \leq 1 \qquad \forall \: j \in M \label{BIP-Destwise} \\
                & \sum_{i \in M} \sum_{j \in M: \lambda_{ij}=1} z_{ij} \geq \frac{|M|}{2} \qquad  \label{BIP-MatchingNb} \\               
                & z_{ij} \in \{0,1\} \qquad \forall \: i \in M, j \in M
                    \label{BIP-Domain-z}
\end{align}
where $M$ is set of unique paths that exist in all explored paths, and $z_{ij}$ is the binary variable which equals to 1 if path $i$ and path $j$ that constitutes a pair is selected to merge; 0, otherwise. The objective function (\ref{BIP-Obj}) minimizes total routing cost of paths that are selected to merge where $e_{ij}$ is the total routing cost of the pair of path $i$ and path $j$. Constraints (\ref{BIP-Originwise}) and (\ref{BIP-Destwise}) restrict to select a path which exists in more than one pair at most once. Constraint (\ref{BIP-MatchingNb}) enforces to select the maximum number of paths to merge from the conflict array. Last, constraints (\ref{BIP-Domain-z}) are the domain constraints for the binary variable $z$. 

\subsubsection{Phase IV: Improving Solution by Mixing Paths}
\label{subsec-DMaM-improve2}

In the fourth phase, another procedure is developed in addition to the third phase to improve the effectiveness of the $DMaM$ in finding high-quality solutions to the problem instances. This procedure focuses on the selected paths that cannot be merged in the third phase (cycles that contain one commodity path). The aim is to decrease the number of selected paths in the current solution by breaking the path into the two parts (at least) and merge these parts with the other selected paths. The overall idea depends on one of the basic assumptions that exist in the problem definition. To remind that holding arcs are assumed to have the infinite capacity that is inherited from many previous studies on service network design problems in the literature.

\begin{figure}[htbp]
	\centering
	\fbox{\includegraphics[width=0.99\textwidth]{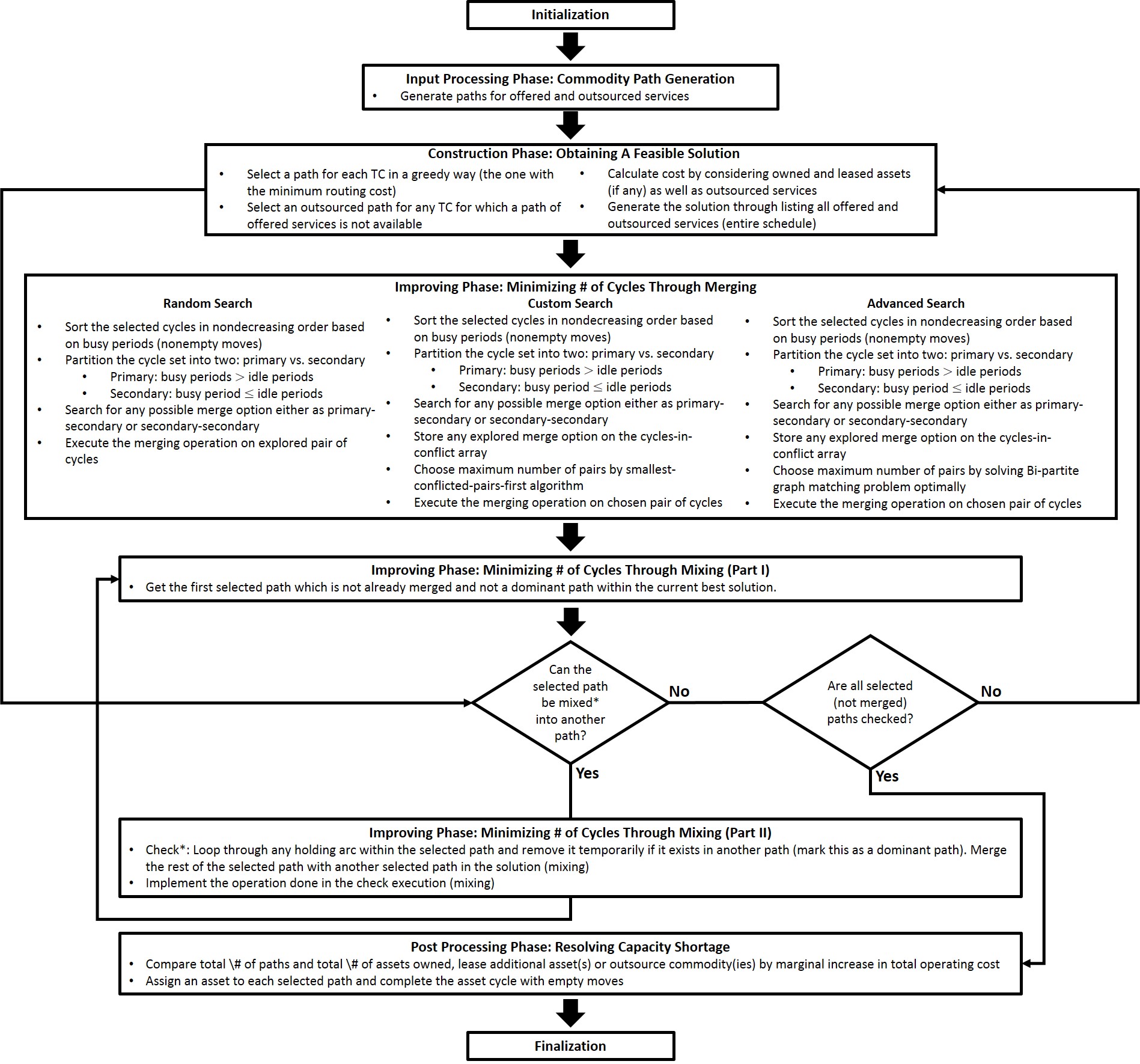}}
	\caption{Flowchart of the DMaM.}
	\label{fig:DMaMflow}
\end{figure}

The procedure searches through unmerged single commodity paths in the current solution. First, it checks every particular holding arc existing in any unmerged single path whether that holding arc also exists in another selected path (called dominant path). If so, that holding path is temporarily removed and the procedure moves to the second step where the rest of the path -also called path particle(s)- is tested for merging with another selected path based on $DMaM-R$ merging search principles. If path particle can be merged with another path, the corresponding holding arc is permanently removed from the path that is just merged. Otherwise, the procedure goes back to the path selection process and replaces the selected path with an alternative commodity path for the same OC (which corresponds to demand shifting decision) and applies the first and second step again. This continues until all alternative paths are checked for the first unmerged single commodity path in the solution. Then this procedure moves to the next unmerged single commodity path and follows the same steps until all unmerged paths are checked. Note that this procedure does not consider any merged and dominant path for mixing operation. The reader can refer to Algorithm~\ref{DMaM_custom} and Figure~\ref{fig:DMaMflow} for illustration.

\subsubsection{Phase V: Resolving Capacity Shortage}
\label{subsec-DMaM-resolve}

In the final phase of $DMaM$, the procedure assigns an asset to each unmerged or merged commodity path and assumes that each asset makes empty moves in its cycle during the rest of the planning horizon. Then the algorithm checks any capacity (asset) shortage through comparing the total number of single paths as well as merged paths (without outsourced paths) and the total number of assets owned in the problem. In case of any asset shortage, the algorithm decides to either lease additional assets or outsource commodities as many as the shortage amount in the capacity based on a marginal increase in total operating cost. 

Furthermore, $DMaM$ creates final asset cycles by completing each cycle with empty moves in addition to service and holding arcs existing in the single/merged commodity paths. In the end, DMaM finalizes the procedure and print the results in terms of asset cycles embedded with TC IDs. The pseudo code of the proposed algorithm ($DMaM$) is presented in Algorithm~\ref{DMaM_custom}.

\begin{algorithm}
    \scriptsize
	\caption{Dedicate-merge-and-mix algorithm (DMaM) w/ random/custom/advanced search} \label{DMaM_custom}
	\begin{algorithmic}[1]
        \State \textit{Initialization}
		\Procedure{Input Processing Phase: Commodity Path Generation}{}
		
		\begin{itemize}
		\item \textit{Generate all possible commodity paths from their origin to destination for all TCs.}
		\item \textit{Generate all paths linked to outsourced services.}
		\end{itemize}
		\EndProcedure
		
		\Procedure{Construction Phase: Obtaining A Feasible Solution}{}		
		\begin{itemize}
		 
		 \item \textit{Choose the path with the minimum cost for any TC among generated commodity paths.}
		 
		 \item \textit{Choose an outsourced path for any TCs for which a path of offered services is not available.}
		 
		 \item \textit{Print selected paths which constitute solution set (entire schedule).}
		\end{itemize}		
		\EndProcedure
		
        \Procedure{Improving Phase: Minimizing \# of Paths Through Merging}{}
		
		\begin{itemize}
		 \item \textit{Sort the paths in the solution set in nonincreasing order based on how many periods in which the asset assigned to corresponding path is busy.}
         \item \textit{Partition the solution set into two distinct sets as primary paths and secondary paths.}
            \begin{itemize}
            \item \textit{Primary: paths on which number of busy periods are greater than that of idle periods.}
            \item \textit{Secondary: paths on which number of busy periods are less than that of idle periods.}
            \end{itemize}

            \For{$\textit{each path within solution set}$}
                \item \textit{Search for any possible merge operation (regular/shifted) on the following combinations: (1) including a primary and a secondary path, and (2) two secondary paths, if any.}  
                
				\item \textit{Jump to last line of this procedure in $DMaM-R$; proceed to next line in $DMaM-C$ and $DMaM-A$.}	\item \textit{Store any explored merge option on the merges-in-conflict array.}		
			\EndFor	                   
        
        \item \textit{Choose maximum number of pairs from the conflict array for merging.}    
            \begin{itemize}
                \item \textit{$DMaM-C$: Identify maximum \# of pairs from conflict array by smallest-conflicted-pairs-first algorithm.}		 
                \item \textit{$DMaM-A$: Identify maximum \# of pairs from conflict array by solving Bi-partite graph matching problem optimally.}	 
            \end{itemize}
            
		\item \textit{Execute the merging operation on chosen pair of paths and calculate total operating cost. Return to the beginning of for loop in $DMaM-R$, proceed to next procedure otherwise.}
		\end{itemize}
		\EndProcedure
		
		\Procedure{Improving Phase: Minimizing \# of Paths Through Mixing}{}
        \State \textit{Loop through all arcs in paths of offered services to check whether a holding arc exists in another path of offered services. If exists, remove it from the current path and merge the rest of path with another path using the principles of random search. If not, continue with next path. Proceed to next procedure when all paths are checked.} 
        \EndProcedure
		
        \Procedure{Post Processing Phase: Resolving Capacity Shortage}{}
        \State \textit{Compare total \# of paths and total \# of assets owned, lease additional asset(s) or outsource commodity(ies) by marginal increase in total operating cost.} 
        
        \EndProcedure
		\State \textit{Create asset cycles by filling each path of offered services with empty moves in addition to existing arcs and assign an asset to each cycle.}
		\State \textit{Finalization and reporting results (asset cycles and TC IDs).}
				
	\end{algorithmic}
\end{algorithm}

\subsection{DMaM Results}
\label{subsec-results_heuristic}

The DMaM is tested on instances larger than small-size in the computational experiment and the results are presented in Tables~\ref{tab:n6_heuristic} for medium-size instances,  Table~\ref{tab:n7_heuristic} for large-size instances, and Table~\ref{tab:n10_heuristic} for very large-size instances with an identical structure. The first column identifies the instance ID, columns two through five presents \textit{best objective value} for the CSSND formulation, DMaM with random search ($DMaM-R$), DMaM with custom search ($DMaM-C$), and DMaM with advanced search ($DMaM-A$), respectively. The \textit{optimality gap} and \textit{CPU time elapsed} in each run are reported on columns seven through ten and eleven through fourteen for each configuration in the same order, respectively. The lower bound (best bound) obtained by the regular formulation is given in column six and used to calculate the optimality gap for each configuration.

As discussed earlier in Section~\S~\ref{subsec-impact-VI}, the proposed formulation is able to solve only nine out of 30 instances within the time limit (four hours) and ends with an optimality gap of 5.01\% on average. The optimality gap fluctuates in a wide range from 1.69\% to 22.36\% for the medium-size instances (see Table~\ref{tab:n6_heuristic}). When it comes to large instances, none of the large instances can be solved within the time limit, and the average gap increases to 10.96\%. The optimality gaps spread over a large interval from 0.61\% to 26.76\% at the end of the time limit on fifteen test problems (see Table~\ref{tab:n7_heuristic}). These outcomes depict the dramatic increase in the complexity when a node is added on the physical network. It also gives a hint about the trade-off between the problem size (complexity) and the computational time overall. This argument is validated transparently by the test of very large instances, OPL CPLEX stopped the solution process of all very large-size instances by throwing \textit{'not responding'} error at early stages based on extensive memory need. Since this error is thrown either in the initialization or presolving period, CPLEX cannot report anything about the root node, even. For this reason, we only denoted the outcomes regarding DMaM heuristic algorithm in Table~\ref{tab:n10_heuristic}. 

When going from DMaM-R to DMaM-A, the average performance of DMaM improves based on the differences in the design of these configurations as expected. Even though a slight improvement is observed on average, the exact gain can be discovered well at the instance level. For instance, the number of assets utilized is decreased by one in instance-36 and 55; the number of commodities outsourced is decreased by one in instance-42 and 43 when going from DMaM-R to DMaM-C. Likewise, the number of assets utilized is decreased by one in instance-34, 46 and 49; the number of commodities outsourced is decreased by one in instance-51 when going from DMaM-C to DMaM-A. In contrast, there are also a few exceptional instances at which an increase in total operating cost is observed as moving from DMaM-R to DMaM-C. In instance-35, 38, 44 and 60 the SCoPF algorithm selects a  subset of path pairs (same number of pairs) that can be merged different from the one chosen in the merging phase of the random search (DMaM-R). This triggers the difference in the mixing phase, either fewer commodity paths or none commodity paths can be mixed into each other. These result in the total operating cost because one more commodity is outsourced in these instances to meet the overall capacity (Table~\ref{tab:n6_heuristic}). Note that, merging and mixing are different modules which aim to maximize their contribution independently rather than the aggregated contribution of DMaM.

All three configurations of DMaM are able to run from the initialization phase to the finalization phase under one second for medium and large-size of instances. The longest computational time is observed in testing very large instances where the maximum time is about 14 seconds. The average optimality gap obtained for medium and large size instances are about 37\% and 43\%, respectively. The performance of DMaM cannot be quantified on problems beyond large size, because OPL CPLEX cannot initialize the solution process for those problems, even.

As stated earlier CSSND formulation is developed based on the existing SND formulations in the literature to be consistent with the previous studies. Earlier versions of this problem had been basically developed for rail transportation and assets has usually corresponded to train cars which have transferred from one block of a train to another en route. Hence, our formulation ($P$) does allow transferring commodities between assets to be consistent with the existing formulations, even though we do not specifically consider rail transportation in this paper. However, DMaM is developed based on the needs of the carriers in practice, which makes the main logic of DMaM different from that of the CSSND formulation. The main idea is to provide the best and the most applicable solution in terms of operating cost and managerial aspect to the decision maker. For this reason, DMaM focuses on dedicated commodity paths -a commodity utilizes the same asset until to reach the destination- at the beginning (construction and merging phases) and allows limited commodity transfer (mixing phase) between assets. This is the main reason behind the difference in optimality gaps between CSSND formulation and the configurations of DMaM algorithm. The bottom line is that DMaM is able to provide an integer solution at an acceptable quality in a reasonable time for the decision maker on any size of the CSSND problem, while the commercial solver cannot obtain any integer solution for the problem, even.  

\begin{table}[htbp]
	\centering
	\tiny
	\caption {Formulation ($P$) vs. $DMaM$ results for instances with $|N|=6$.} \label{tab:n6_heuristic}
	{
		\begin{tabular}{l c r r r r r r r r r r r r r r r r}
			\midrule
\multirow{3}{*}{\bf ID}	&		&	\multicolumn{4}{c}{\bf Best Objective Value}							&		&	\multicolumn{1}{c}{\bf Best}	&		&	\multicolumn{4}{c}{\bf Optimality Gap (\%)}							&		&	\multicolumn{4}{c}{\bf CPU Time (sec.)}							\\	\cmidrule{3-6}	\cmidrule{10-13}	\cmidrule{15-18}
	&		&	\multirow{2}{*}{\bf \textbf{$P$}}	&	\multirow{2}{*}{\bf \textbf{$D^R$}}	&	\multirow{2}{*}{\bf \textbf{$D^C$}}	&	\multirow{2}{*}{\bf \textbf{$D^A$}}	&		&	\multicolumn{1}{c}{\bf Bound}	&		&	\multirow{2}{*}{\bf \textbf{$P$}}	&	\multirow{2}{*}{\bf \textbf{$D^R$}}	&	\multirow{2}{*}{\bf \textbf{$D^C$}}	&	\multirow{2}{*}{\bf \textbf{$D^A$}}	&		&	\multirow{2}{*}{\bf \textbf{$P$}}	&	\multirow{2}{*}{\bf \textbf{$D^R$}}	&	\multirow{2}{*}{\bf \textbf{$D^C$}}	&	\multirow{2}{*}{\bf \textbf{$D^A$}}	\\			
	&		&		&		&		&		&		&		&		&		&		&		&		&		&		&		&		&		\\	\midrule		
31	&		&	270.54	&	395.64	&	395.052	&	395.052	&		&	258.565	&		&	4.43	&	34.65	&	34.55	&	34.55	&		&	14400	&	0.109	&	0.156	&	0.276	\\			
32	&		&	325.624	&	476.356	&	476.662	&	476.164	&		&	305.295	&		&	6.24	&	35.91	&	35.95	&	35.88	&		&	14400	&	0.156	&	0.203	&	0.202	\\			
33	&		&	380.984	&	558.418	&	558.106	&	557.998	&		&	335.929	&		&	11.83	&	39.84	&	39.81	&	39.8	&		&	14400	&	0.141	&	0.265	&	0.402	\\	\midrule		
34	&		&	200.816	&	346.346	&	346.602	&	321.186	&		&	200.816	&		&	0	&	42.02	&	42.06	&	37.48	&		&	2130.6	&	0.172	&	0.266	&	0.233	\\			
35	&		&	252.796	&	376.142	&	401.416	&	375.036	&		&	242.196	&		&	4.19	&	35.61	&	39.66	&	35.42	&		&	14400	&	0.172	&	0.234	&	0.27	\\			
36	&		&	283.342	&	505.62	&	481.3	&	482.336	&		&	263.736	&		&	6.92	&	47.84	&	45.2	&	45.32	&		&	14400	&	0.219	&	0.28	&	0.312	\\	\midrule		
37	&		&	221.202	&	317.876	&	318.478	&	317.574	&		&	221.202	&		&	0	&	30.41	&	30.54	&	30.35	&		&	494.48	&	0.14	&	0.141	&	0.202	\\			
38	&		&	273.25	&	373.94	&	399.304	&	373.518	&		&	260.208	&		&	4.77	&	30.41	&	34.83	&	30.34	&		&	14400	&	0.125	&	0.171	&	0.244	\\			
39	&		&	285.23	&	479.262	&	478.392	&	478.354	&		&	285.23	&		&	0	&	40.49	&	40.38	&	40.37	&		&	3435.58	&	0.172	&	0.218	&	0.286	\\	\midrule		
40	&		&	271.176	&	372.492	&	372.492	&	371.9	&		&	256.478	&		&	5.42	&	31.15	&	31.15	&	31.04	&		&	14400	&	0.109	&	0.141	&	0.177	\\			
41	&		&	303.052	&	452.762	&	453.328	&	452.676	&		&	297.939	&		&	1.69	&	34.2	&	34.28	&	34.18	&		&	14400	&	0.141	&	0.172	&	0.248	\\			
42	&		&	334.638	&	560.384	&	534.858	&	534.764	&		&	316.922	&		&	5.29	&	43.45	&	40.75	&	40.74	&		&	14400	&	0.172	&	0.219	&	0.254	\\	\midrule		
43	&		&	221.326	&	369.43	&	344.138	&	343.742	&		&	221.326	&		&	0	&	40.09	&	35.69	&	35.61	&		&	8976.84	&	0.141	&	0.187	&	0.181	\\			
44	&		&	250.714	&	374.036	&	400.098	&	398.842	&		&	240.476	&		&	4.08	&	35.71	&	39.9	&	39.71	&		&	14400	&	0.187	&	0.234	&	0.27	\\			
45	&		&	283.134	&	480.072	&	479.944	&	478.758	&		&	265.451	&		&	6.25	&	44.71	&	44.69	&	44.55	&		&	14400	&	0.141	&	0.265	&	0.27	\\	\midrule		
46	&		&	220.764	&	319.866	&	319.422	&	294.078	&		&	208.502	&		&	5.55	&	34.82	&	34.73	&	29.1	&		&	14400	&	0.141	&	0.188	&	0.233	\\			
47	&		&	252.612	&	350.848	&	350.506	&	350.628	&		&	231.946	&		&	8.18	&	33.89	&	33.83	&	33.85	&		&	14400	&	0.187	&	0.218	&	0.26	\\			
48	&		&	258.958	&	480.126	&	480.316	&	478.798	&		&	258.958	&		&	0	&	46.06	&	46.09	&	45.91	&		&	4624.15	&	0.203	&	0.296	&	0.27	\\	\midrule		
49	&		&	270.074	&	396.25	&	396.25	&	371.188	&		&	263.367	&		&	2.48	&	33.54	&	33.54	&	29.05	&		&	14400	&	0.156	&	0.172	&	0.192	\\			
50	&		&	279.326	&	451.628	&	451.232	&	450.796	&		&	279.326	&		&	0	&	38.15	&	38.1	&	38.04	&		&	3498.99	&	0.187	&	0.203	&	0.383	\\			
51	&		&	331.228	&	533.434	&	532.548	&	505.852	&		&	314.271	&		&	5.12	&	41.09	&	40.99	&	37.87	&		&	14400	&	0.171	&	0.203	&	0.443	\\	\midrule		
52	&		&	320.778	&	472.114	&	472.114	&	472.114	&		&	309.138	&		&	3.63	&	34.52	&	34.52	&	34.52	&		&	14400	&	0.125	&	0.156	&	0.181	\\			
53	&		&	402.98	&	528.722	&	528.722	&	528.722	&		&	330.667	&		&	17.94	&	37.46	&	37.46	&	37.46	&		&	14400	&	0.141	&	0.171	&	0.197	\\			
54	&		&	482.826	&	558.322	&	558.322	&	558.322	&		&	374.872	&		&	22.36	&	32.86	&	32.86	&	32.86	&		&	14400	&	0.156	&	0.187	&	0.212	\\	\midrule		
55	&		&	193.49	&	292.816	&	268.324	&	266.684	&		&	193.49	&		&	0	&	33.92	&	27.89	&	27.45	&		&	196.05	&	0.125	&	0.172	&	0.222	\\			
56	&		&	223.668	&	349.846	&	350.934	&	346.826	&		&	223.668	&		&	0	&	36.07	&	36.26	&	35.51	&		&	3626.18	&	0.218	&	0.203	&	0.207	\\			
57	&		&	254.224	&	479.002	&	479.488	&	475.914	&		&	254.224	&		&	0	&	46.93	&	46.98	&	46.58	&		&	12707.5	&	0.172	&	0.25	&	0.222	\\	\midrule		
58	&		&	294.856	&	445.272	&	445.272	&	445.272	&		&	276.549	&		&	6.21	&	37.89	&	37.89	&	37.89	&		&	14400	&	0.172	&	0.187	&	0.198	\\			
59	&		&	327.904	&	527.282	&	527.282	&	526.784	&		&	316.84	&		&	3.37	&	39.91	&	39.91	&	39.85	&		&	14400	&	0.156	&	0.172	&	0.232	\\			
60	&		&	409.11	&	584.526	&	610.118	&	583.69	&		&	350.919	&		&	14.22	&	39.97	&	42.48	&	39.88	&		&	14400	&	0.156	&	0.218	&	0.243	\\	\midrule		
\multicolumn{3}{l}{\bf Average}				-	&	-	&	-	&	-	&		&	-	&		&	5.01	&	37.79	&	37.77	&	36.71	&		&	11403.01	&	0.159	&	0.205	&	0.251	\\			
			\hline
		\end{tabular}
	}
\end{table}

\begin{table}[htbp]
	\centering
	\tiny
	\caption {Formulation ($P$) vs. $DMaM$ results for instances with $|N|=7$.} \label{tab:n7_heuristic}
	{
		\begin{tabular}{l c r r r r r r r r r r r r r r r r}
			\midrule
\multirow{3}{*}{\bf ID}	&		&	\multicolumn{4}{c}{\bf Best Objective Value}							&		&	\multicolumn{1}{c}{\bf Best}	&		&	\multicolumn{4}{c}{\bf Optimality Gap (\%)}							&		&	\multicolumn{4}{c}{\bf CPU Time (sec.)}							\\	\cmidrule{3-6}	\cmidrule{10-13}	\cmidrule{15-18}
	&		&	\multirow{2}{*}{\bf \textbf{$P$}}	&	\multirow{2}{*}{\bf \textbf{$D^R$}}	&	\multirow{2}{*}{\bf \textbf{$D^C$}}	&	\multirow{2}{*}{\bf \textbf{$D^A$}}	&		&	\multicolumn{1}{c}{\bf Bound}	&		&	\multirow{2}{*}{\bf \textbf{$P$}}	&	\multirow{2}{*}{\bf \textbf{$D^R$}}	&	\multirow{2}{*}{\bf \textbf{$D^C$}}	&	\multirow{2}{*}{\bf \textbf{$D^A$}}	&		&	\multirow{2}{*}{\bf \textbf{$P$}}	&	\multirow{2}{*}{\bf \textbf{$D^R$}}	&	\multirow{2}{*}{\bf \textbf{$D^C$}}	&	\multirow{2}{*}{\bf \textbf{$D^A$}}	\\			
	&		&		&		&		&		&		&		&		&		&		&		&		&		&		&		&		&		\\	\midrule		
61	&		&	433.19	&	632.37	&	632.37	&	631.86	&		&	386.77	&		&	10.72	&	38.84	&	38.84	&	38.79	&		&	14400	&	0.328	&	0.312	&	0.332	\\			
62	&		&	489.77	&	665.16	&	665.16	&	689.35	&		&	405.19	&		&	17.27	&	39.08	&	39.08	&	41.22	&		&	14400	&	0.281	&	0.344	&	0.363	\\			
63	&		&	572.7	&	749	&	749.25	&	748.15	&		&	455.21	&		&	20.51	&	39.22	&	39.24	&	39.15	&		&	14400	&	0.327	&	0.421	&	0.41	\\	\midrule		
64	&		&	229.06	&	429.67	&	428.33	&	401.22	&		&	212.11	&		&	7.4	&	50.63	&	50.48	&	47.13	&		&	14400	&	0.249	&	0.343	&	0.29	\\			
65	&		&	263.94	&	559.44	&	560.59	&	557.12	&		&	254.72	&		&	3.49	&	54.47	&	54.56	&	54.28	&		&	14400	&	0.359	&	0.421	&	0.363	\\			
66	&		&	297.08	&	716.72	&	718.86	&	712.04	&		&	284.2	&		&	4.34	&	60.35	&	60.46	&	60.09	&		&	14400	&	0.359	&	0.546	&	0.426	\\	\midrule		
67	&		&	334.22	&	482.25	&	507.03	&	455.93	&		&	314.28	&		&	5.97	&	34.83	&	38.02	&	31.07	&		&	14400	&	0.234	&	0.265	&	0.229	\\			
68	&		&	392	&	560.84	&	560.27	&	559.71	&		&	355.07	&		&	9.42	&	36.69	&	36.62	&	36.56	&		&	14400	&	0.203	&	0.296	&	0.285	\\			
69	&		&	425.93	&	744.13	&	718.36	&	716.72	&		&	376	&		&	11.72	&	49.47	&	47.66	&	47.54	&		&	14400	&	0.234	&	0.406	&	0.342	\\	\midrule		
70	&		&	408.36	&	583.92	&	583.92	&	583.65	&		&	347.18	&		&	14.98	&	40.54	&	40.54	&	40.52	&		&	14400	&	0.156	&	0.203	&	0.258	\\			
71	&		&	493.69	&	716.29	&	716.29	&	691.32	&		&	406.81	&		&	17.6	&	43.21	&	43.21	&	41.16	&		&	14400	&	0.234	&	0.265	&	0.289	\\			
72	&		&	627.33	&	825.42	&	825.42	&	825.07	&		&	459.44	&		&	26.76	&	44.34	&	44.34	&	44.32	&		&	14400	&	0.234	&	0.297	&	0.331	\\	\midrule		
73	&		&	306.01	&	479.47	&	478.99	&	453.17	&		&	304.15	&		&	0.61	&	36.57	&	36.5	&	32.88	&		&	14400	&	0.234	&	0.312	&	0.581	\\			
74	&		&	368.25	&	561.84	&	564.2	&	563.15	&		&	353.53	&		&	4.00	&	37.08	&	37.34	&	37.22	&		&	14400	&	0.343	&	0.375	&	0.431	\\			
75	&		&	426.55	&	721.03	&	721.27	&	719.43	&		&	385.33	&		&	9.66	&	46.56	&	46.58	&	46.44	&		&	14400	&	0.344	&	0.499	&	0.415	\\	\midrule		
\multicolumn{3}{l}{\bf Average}				-	&	-	&	-	&	-	&		&	-	&		&	\bf 10.96	&	\bf 43.46	&	\bf 43.56	&	\bf 42.56	&		&	\bf 14400	&	\bf 0.275	&	\bf 0.354	&	\bf 0.356	\\			

			\hline
		\end{tabular}
	}
\end{table}

\begin{table}[htbp]
	\centering
	\tiny
	\caption {Formulation ($P$) vs. $DMaM$ results for instances with $|N|=10+$.} \label{tab:n10_heuristic}
	{
		\begin{tabular}{l c r r r r r r r r r r r r r r r r}
			\midrule
\multirow{3}{*}{\bf ID}	&		&	\multicolumn{4}{c}{\bf Best Objective Value}							&		&	\multicolumn{1}{c}{\bf Best}	&		&	\multicolumn{4}{c}{\bf Optimality Gap (\%)}							&		&	\multicolumn{4}{c}{\bf CPU Time (sec.)}							\\	\cmidrule{3-6}	\cmidrule{10-13}	\cmidrule{15-18}
	&		&	\multirow{2}{*}{\bf \textbf{$P$}}	&	\multirow{2}{*}{\bf \textbf{$D^R$}}	&	\multirow{2}{*}{\bf \textbf{$D^C$}}	&	\multirow{2}{*}{\bf \textbf{$D^A$}}	&		&	\multicolumn{1}{c}{\bf Bound}	&		&	\multirow{2}{*}{\bf \textbf{$P$}}	&	\multirow{2}{*}{\bf \textbf{$D^R$}}	&	\multirow{2}{*}{\bf \textbf{$D^C$}}	&	\multirow{2}{*}{\bf \textbf{$D^A$}}	&		&	\multirow{2}{*}{\bf \textbf{$P$}}	&	\multirow{2}{*}{\bf \textbf{$D^R$}}	&	\multirow{2}{*}{\bf \textbf{$D^C$}}	&	\multirow{2}{*}{\bf \textbf{$D^A$}}	\\			
	&		&		&		&		&		&		&		&		&		&		&		&		&		&		&		&		&		\\	\midrule		
76	&		&	N/A	&	1250.542	&	1250.602	&	1223.66	&		&	N/A	&		&	N/A	&		&		&		&		&	N/A	&	0.858	&	1.124	&	0.878	\\			
77	&		&	N/A	&	1409.866	&	1409.828	&	1382.698	&		&	N/A	&		&	N/A	&		&		&		&		&	N/A	&	0.936	&	1.373	&	1.08	\\			
78	&		&	N/A	&	1469.864	&	1470.162	&	1467.516	&		&	N/A	&		&	N/A	&		&		&		&		&	N/A	&	1.03	&	1.701	&	1.102	\\	\midrule		
79	&		&	N/A	&	1790.28	&	1840.806	&	1812.084	&		&	N/A	&		&	N/A	&		&		&		&		&	N/A	&	1.872	&	3.058	&	2.058	\\			
80	&		&	N/A	&	1979.404	&	1954.002	&	1951.514	&		&	N/A	&		&	N/A	&		&		&		&		&	N/A	&	2.028	&	3.76	&	2.309	\\			
81	&		&	N/A	&	2262.82	&	2219.368	&	2186.47	&		&	N/A	&		&	N/A	&		&		&		&		&	N/A	&	2.215	&	4.539	&	2.715	\\	\midrule		
82	&		&	N/A	&	2942.47	&	3016.514	&	2961.716	&		&	N/A	&		&	N/A	&		&		&		&		&	N/A	&	4.93	&	10.108	&	6.256	\\			
83	&		&	N/A	&	3185.032	&	3210.258	&	3182.158	&		&	N/A	&		&	N/A	&		&		&		&		&	N/A	&	5.304	&	12.121	&	5.807	\\			
84	&		&	N/A	&	3504.906	&	3481.822	&	3525.348	&		&	N/A	&		&	N/A	&		&		&		&		&	N/A	&	5.71	&	14.024	&	6.266	\\			
			\hline
		\end{tabular}
	}
\end{table}

\newpage\section{Conclusion and Future Work}
\label{sec-conclusion}





In this study, we introduce a variant of the service network design problem in which an excessive demand period such as Christmas or Thanksgiving is ahead, and the capacity of the carrier's assets is insufficient to meet that demand. The problem aims to minimize total operating cost and determine the best set of actions to take in responding to excessive demand while satisfying all customers' demand for a consolidation-based freight carrier. This is one of the gaps in the service network design literature because the assumption of sufficient assets' capacity is commonly considered in previous studies. 
We consider three quick actions to respond to peak demand and resolve capacity shortage: (1) shifting demand in time by negotiating the customer and bearing the penalty incurred for an early pick-up or a late delivery, (2) leasing an additional asset (adding asset temporarily) to increase the capacity for the peak demand period, and (3) outsourcing particular amount of capacity from a competitor or other service provider in exchange of a cost.

Our contributions beyond considering a real-world assumption can be summarized as follows: (1) Developing a mixed-integer formulation for the \textit{capacity scaling service network design} (CSSND) problem with a set-covering problem approach for demand shifting feature, (2) proposing valid inequalities for the CSSND formulation, and (3) developing a multi-phase custom \textit{dedicate-merge-and-mix algorithm} (DMaM) to solve the problem for practical purposes.


In the computational experiment of this study, we investigate the performance of MIP formulation as well as valid inequalities on randomly generated test instances. Besides, we conduct additional analysis to clarify possible impacts of cost parameters as well as operational restrictions and limitations that can be observed in practice. The experiment reveals the fact that the studied problem is computationally NP-Hard and solution times are higher even for small and medium size instances. The proposed valid inequalities can achieve time savings for medium-size instances up to 25\% based on optimally solved instances and 3\% overall. The formulations with valid inequalities outperform the regular formulation in obtaining the best bound for medium-size instances, thus shows that valid inequalities make the CSSND formulation tighter. 
Despite these achievements, an efficient solution approach is required to solve real-life instances of CSSND seen in practice. In this sense, the DMaM -a multi-phase procedure- is developed to obtain solutions good enough for CSSND, especially in a reasonable computational time.

The three configurations of proposed DMaM are able to run from the initialization phase to the finalization phase under one second for medium and large-size of instances for which CSSND formulation cannot reach the optimal solution in 36 out of 45 instances within the time limit of four hours. Moreover, the commercial solver, OPL CPLEX, could not produce an initial integer feasible solution for very large size instances whereas DMaM solves these instances between one second and 14 seconds. The proposed DMaM performs well enough and is able to provide an integer solution at an acceptable quality (as low as 29.05\%) in a reasonable time for the decision maker on any size of the CSSND problem. However, there is still a gap to improve DMaM in terms of problem-solving approach. DMaM is developed based on the needs of the carriers in practice; it focuses on dedicated commodity paths in the construction and merging phases but allows a limited commodity transfer between assets in the mixing phase. This approach reduces the additional handling of commodities and increases the managerial impact. From the perspective of integer programming, DMaM initializes at a point which is exceptionally far from the optimum in the solution space comparing to CSSND formulation. This is the main reason why DMaM cannot get closer to the optimality gaps of the CSSND formulation. To overcome this issue, the idea of arc swapping (commodity transfers) between asset cycles should be incorporated into the DMaM. Arc swapping is a computationally excessive operation, which may hurt the efficiency of the problem-solving process. This issue can be further considered as a part of future venues of this study. We aim to commodity transfer
into the DMaM in an efficient way by considering the trade-off between problem-solving effectiveness and computational efficiency.

\bibliography{GenericSND}{}
\bibliographystyle{acm}

\end{document}